\documentclass[12pt,twoside,reqno]{amsart}

\usepackage{amsmath}
\usepackage{amssymb}
\usepackage{latexsym}
\usepackage[all]{xy}

\newtheorem{theorem}{\bf Theorem}
\newtheorem{lemma}[theorem]{\bf Lemma}
\newtheorem{corollary}[theorem]{\bf Corollary}
\newtheorem{proposition}[theorem]{\bf Proposition}
\newtheorem{definition}[theorem]{\bf Definition}
\newtheorem{example}[theorem]{\bf Example}
\newtheorem{remark}[theorem]{\bf Remark}

\def\proof{\noindent {\sc Proof:}\enspace}
\def\endproof{\removelastskip\rightline{$\Box$}\par\bigskip}

\def\conj#1{ \{ #1 \} }

\begin{document}
\setcounter{page}{1}
\title{Derived classification of gentle algebras with two cycles}
\author{Diana Avella-Alaminos}
\address{Diana Avella-Alaminos\newline Departamento de Matem\'aticas \newline Facultad de Ciencias \newline Universidad Nacional Aut\'onoma de M\'exico.
        \newline Ciudad Universitaria\newline M\'exico D.F. C.P. 04510, M\'exico}
\email{avella@matem.unam.mx}
\email{avella@matematicas.unam.mx}
\date{August 2008}

\begin{abstract}
We classify gentle algebras defined by quivers with two cycles
under derived equivalence in a non degenerate case, by using some
combinatorial invariants constructed from the quiver with
relations defining these algebras. We also present a list of
normal forms; any such algebra is derived equivalent to one of the
algebras in the list. The article includes an Appendix presenting
a slightly modified and extended version of a technical result in
the unpublished manuscript \cite{jan01} by Holm, Schr\"{o}er and
Zimmermann, describing some essential elementary transformations
over the quiver with relations defining the algebra.
\end{abstract}
\maketitle

\section{Introduction}\label{intro}
Let $A$ be a finite-dimensional connected ${\mathrm k}$-algebra
$A$ over an algebraically closed field ${\mathrm k}$. Denote by
$D^b(A)$ the bounded derived category of the module category of
finite-dimensional left $A$-modules, $A$-$\operatorname{mod}$. It
is an interesting problem to classify such algebras up to derived
equivalence.

In particular, the family of gentle algebras is closed under
derived equivalence \cite{sz}. The problem of classifying gentle
algebras up to derived equivalence is well understood in the case
where the associated quiver has one cycle, see \cite{ah81},
\cite{as87}, \cite{dv01}, \cite{chj99} and \cite{bg04}. In this
paper we focus our attention on gentle algebras with two cycles.

We use combinatorial invariants $\phi_A:\mathbb{N}^2\to\mathbb{N}$
defined in \cite{ag07} in order to classify them under derived
equivalence. Roughly speaking $\phi_A$ is obtained as follows:
Start with a maximal directed path in $Q$ which contains no
relations. Then continue in opposite direction as long as possible
with zero relations. Repeat this until the first path appears
again, say after $n$ steps. Then we obtain a pair $(n,m)$ where
$m$ is the number of arrows which appeared in a zero relation.
Repeat this procedure until all maximal paths without a zero
relation have been used; $\phi_A$ counts then how often each pair
$(n,m)\in \mathbb{N}^2$ occurred. Recall $\phi_A$ has always a
finite support. Let $\{(n_1,m_1),(n_2,m_2),\dots ,(n_k,m_k)\}$ be
the support of $\phi_A$ , denote $\phi_A$ by
$[(n_1,m_1),(n_2,m_2),\dots ,(n_k,m_k)]$ where each $(n_j,m_j)$ is
written $\phi_A(n_j,m_j)$ times and the order in which they are
written is arbitrary. Define also $\# \phi_A :=\sum _{1\leq j \leq
k } \phi_A(n_j,m_j)$. See \cite[3,5]{ag07} for a precise
description.

We can show by induction over the number of vertices that:

\vspace{0.5 cm}

 {\bfseries Theorem I.} {\em If $A=\mathrm{k}Q/ \left< \mathcal{P} \right>$ is a gentle
 algebra, $Q$ a quiver with two cycles, then $\# \phi_A \in \{ 1,3\}$.}

\vspace{0.5 cm}

In the case $\# \phi_A = 3$, we give the complete classification
of gentle algebras with quivers of two cycles, see \cite{ag07},
under derived equivalence. This main result is also proved by
induction:

\vspace{0.5 cm}

{\bfseries Theorem II.} {\em Let $A=\mathrm{k}Q/ \left<
\mathcal{P} \right>$ and $B=\mathrm{k}Q'/ \left< \mathcal{P}'
\right>$ be gentle algebras so that $Q$ and $Q'$ are quivers with
two cycles. If $\#\phi_A=3$ then,

$A$ and $B$ are derived equivalent if and only if $\phi_A=\phi_B$.
}
 \vspace{0.5 cm}

For the proof of this main result we use strongly some
combinatorial elementary transformations over the quiver with
relations which defines a gentle algebra, introduced by Holm,
Schr\"{o}er and Zimmermann in the unpublished manuscript
\cite{jan01}. As they are constructed by using tilting complexes,
they produce derived equivalent algebras, see Appendix.

In Section \ref{preliminares} we present the basic definitions and
notation about gentle algebras, including some previous results
with respect to derived invariants introduced in \cite{ag07}. In
Section \ref{representantes} we give a list of normal forms with a
representative of each derived equivalence class of gentle
algebras with quivers of two cycles in the non degenerate case
$\#\phi_A=3$. It follows from the proof of Theorem II that any
gentle algebra $A$ with two cycles, $\#\phi_A=3$ and more than
four vertices is not only derived equivalent to one of the
algebras in the list, but also can be modified into it by applying
some transformations from the ones defined in the Appendix. In
Section \ref{transhilos} we give an alternative way to apply such
transformations.

Section \ref{reduc} consists of technical definitions, lemmas and
propositions used for the proofs of Theorems I and II presented in
Section \ref{pruebas}.

The results of this article appeared as part of the Ph.D. thesis
of the author who graduated on August 30, 2006. The complete
digital version of the thesis can be found on the server
http://bidi.unam.mx/ or more precisely on \\
http://132.248.9.9:8080/tesdig/Procesados\_2006/0608299/Index.html.

\section{Preliminaries}\label{preliminares}

Let $Q_0$ be a set of vertices, $Q_1$ a set of arrows and $s,e:Q_1
\rightarrow Q_0$ functions which define the start resp. end point
of each arrow, the tuple $Q=(Q_0,Q_1,s,e)$ is called a {\em
quiver}.  For a finite and connected quiver $Q$ {\em the number of
cycles of $Q$} is the least number of arrows that we have to
remove from $Q$ in order to obtain a tree. It is denoted by $c(Q)$
and can be calculated by the expression $c(Q)=\#Q_1- \#Q_0+1$. A
sequence of arrows $C=\alpha_n \dots \alpha_2 \alpha_1$ with
$s(\alpha_{i+1})=e(\alpha_{i})$ for $1\leq i < n$ is called a {\em
path} of {\em length} $n$, the length is denoted by $l(C):=n$.
Also, for each $v\in Q_0$ we consider a {\em trivial path} $1_v$
of length zero. The functions $s$ and $e$ are extended to paths in
the obvious way. We define inverse of paths as follows:
$1_v^{-1}:=1_v$ for $v\in Q_0$, $\alpha ^{-1} $ with
$s(\alpha^{-1}):=e(\alpha)$, $e(\alpha ^{-1}):=s(\alpha)$ and
$(\alpha ^{-1})^{-1}:=\alpha$ for $\alpha \in Q_1$ and
$C^{-1}:=\alpha_1^{-1} \alpha_{2}^{-1} \dots \alpha_n^{-1}$ for a
path $C=\alpha_n \dots \alpha_2 \alpha_1$. The composition of
paths $C_1$ and $C_2$ in $Q$ is the concatenation of them if
$s(C_2)=e(C_1)$ or $0$ if $s(C_2)\neq e(C_1)$. Let ${\mathrm k}$
be a field and ${\mathrm k}Q$ the {\em path algebra}, with paths
of $Q$ as a basis and multiplication induced by concatenation. A
{\em relation in $Q$} is a non zero linear combination of paths of
length at least two, with the same start point and end point. We
work with algebras ${\mathrm k}Q/ \left< \mathcal{P} \right>$
where $\left< \mathcal{P} \right>$ is the ideal of ${\mathrm k}Q$
generated by a set of relations $\mathcal{P}$. A path in $Q$ is
identified with its corresponding class in ${\mathrm k}Q/ \left<
\mathcal{P} \right>$.

The path algebras we study in this work fulfill very particular
conditions.

\begin{definition}\label{biserial}{\normalfont
We call ${\mathrm k}Q/ \left< \mathcal{P} \right>$ a {\em gentle
algebra} if the following five conditions hold:
\begin{enumerate}
\item For each $v \in Q_0$, $\#\{ \alpha \in Q_1 | s(\alpha) = v
\} \leq 2$ and $\#\{ \alpha \in Q_1 | e(\alpha) = v \} \leq 2$.
\item For each $\beta \in Q_1$, $\#\{ \alpha \in Q_1 | s(\beta) =
e(\alpha) \text{ and } \beta \alpha \notin \mathcal{P} \} \leq 1$
and $\#\{ \gamma \in Q_1 | s(\gamma) = e(\beta) \text{ and }\gamma
\beta \notin \mathcal{P} \} \leq 1$. \item For each $\beta \in
Q_1$ there is a bound $n(\beta)$ such that any path
$\beta_{n(\beta)}\dots \beta_2 \beta_1$ with
$\beta_{n(\beta)}=\beta $ or $\beta_{1}=\beta$ contains a subpath
in $\mathcal{P}$.
\item All relations in $\mathcal{P}$ are
monomials of length $2$. \item For each $\beta \in Q_1$, $\#\{
\alpha \in Q_1 | s(\beta) = e(\alpha) \text{ and } \beta \alpha
\in \mathcal{P} \} \leq 1$ and
\\$\#\{ \gamma \in Q_1 | s(\gamma) = e(\beta) \text{ and } \gamma
\beta \in \mathcal{P} \} \leq 1$.
\end{enumerate}
}\end{definition}

We will make a slight abuse of notation by talking about gentle
algebras referring to the quiver with relations which define those
algebras.

\subsection{Threads of a gentle algebra}
Let $A$ be a gentle algebra. A {\em permitted path} of $A$ is a
path $C=\alpha_n \dots \alpha_2 \alpha_1$ with no zero relations,
it is called a {\em non trivial permitted thread} of $A$ if it is
of maximal length, that is, for all $\beta\in Q_1$, neither
$C\beta$ nor $\beta C$ is a permitted path. If $v$ is a vertex
with $\#\{ \alpha \in Q_1 | s(\alpha)=v  \} \leq 1$, $\#\{ \alpha
\in Q_1 | e(\alpha)=v  \} \leq 1$ and if $\beta ,\gamma \in Q_1$
are such that $s(\gamma)=v=e(\beta)$ then $\gamma \beta \notin
\mathcal{P}$, we consider $1_v$ a {\em trivial permitted thread}
in $v$ and denote it by $h_{v}$. Similarly a {\em forbidden path}
of $A$ is a sequence $\Pi=\alpha_n \dots \alpha_2 \alpha_1$ formed
by pairwise different arrows in $Q$ with $\alpha_{i+1} \alpha_{i}
\in \mathcal{P}$ for all $i \in \{ 1,2,\ldots n-1\}$ and it is
called a {\em non trivial forbidden thread} if for all $\beta\in
Q_1$, neither $\Pi\beta$ nor $\beta \Pi$ is a forbidden path. If
$v$ is a vertex with $\#\{ \alpha \in Q_1 | s(\alpha)=v  \} \leq
1$, $\#\{ \alpha \in Q_1 | e(\alpha)=v  \} \leq 1$ and if $\beta
,\gamma \in Q_1$ are such that $s(\gamma)=v=e(\beta)$ then $\gamma
\beta \in \mathcal{P}$, we consider $1_v$ a {\em trivial forbidden
thread} in $v$ and denote it by $p_{v}$. The existence of non
trivial permitted threads is due to point $(3)$ in Definition
\ref{biserial} and the existence of non trivial forbidden threads
is due to the restriction of considering pairwise different
arrows.

Denote by $\mathcal{H}_A$ the set of all permitted threads of $A$,
trivial and non trivial. This set describes completely the algebra
$A$. Notice that certain paths can be permitted and forbidden
threads at the same time.

For a gentle algebra the relations in its quiver can be described
by using two functions $\sigma ,\varepsilon :Q_1 \to \{ 1,-1\} $
as in \cite[3]{br87}, defined by:
\begin{enumerate}
\item If $\beta _1 \neq \beta _2$ are arrows with $s(\beta
_1)=s(\beta _2)$, then $\sigma (\beta _1)= -\sigma (\beta _2)$.
\item If $\gamma _1 \neq \gamma _2$ are arrows with $e(\gamma
_1)=e(\gamma _2)$, then $\varepsilon (\beta _1)=-\varepsilon
(\beta _2)$. \item If $\beta$, $\gamma$ are arrows with $s(\gamma
)=e(\beta )$ and $\gamma \beta \notin \mathcal{P}$, then $\sigma
(\gamma)= -\varepsilon(\beta)$.
\end{enumerate}

They can be extended to the threads of $A$ as follows. For
$H=\alpha_n \dots \alpha_2 \alpha_1$ a non trivial thread of $A$
define $\sigma (H):=\sigma (\alpha_1)$ and
$\varepsilon(H):=\varepsilon(\alpha_n)$. Consider $Q$ has at least
two vertices, as $Q$ is connected, for $v\in Q_0$ there exists
some $\gamma \in Q_1$ with $s(\gamma)=v$ or $\beta \in Q_1$ with
$e(\beta)=v$. If there is a trivial permitted thread at $v$ let
$\sigma (h_v)=-\varepsilon (h_v):=-\sigma (\gamma)$ in the first
case and $\sigma(h_v)=-\varepsilon(h_v):=\varepsilon (\beta)$ in
the second one. For a trivial forbidden thread $p_{v}$ define
$\sigma (p_v)=\varepsilon (p_v):=-\sigma (\gamma)$ in the first
case and $\sigma(p_v)=\varepsilon(p_v):=-\varepsilon (\beta)$
otherwise.

The considered relations are monomials of length two, so we will
indicate them in the quiver by using doted lines, joining each
pair of arrows which form a relation.

Identify $\alpha \in Q_1$ with $(e(\alpha)^{-\varepsilon
(\alpha)},s(\alpha)^{\sigma (\alpha)})$. More generally, identify
the non permitted trivial thread $\alpha_n \dots \alpha_2
\alpha_1$ with the vector
$$(e(\alpha_n)^{-\varepsilon (\alpha_n)}, s(\alpha_n)^{\sigma (\alpha_n)},
\cdots , s(\alpha_2)^{\sigma (\alpha_2)}, s(\alpha_1)^{\sigma
(\alpha_1)})$$ and its inverse with
$$(s(\alpha_1)^{\sigma (\alpha_1)}, s(\alpha_2)^{\sigma (\alpha_2)}, \cdots , s(\alpha_n)^{\sigma (\alpha_n)}, e(\alpha_n)^{-\varepsilon (\alpha_n)}).$$
 If there is a trivial permitted thread $h_{v}$ for some $v \in Q_0$ we identify it with $(v^{ \sigma (h_v)})$. For simplicity we write $v^{+}$ instead of $v^{+1}$
 and $v^{-}$ instead of $v^{-1}$.
We can describe then $\mathcal{H}_A$ by an array $[c_1 c_2 \cdots
c_r ]$ where the $c_i$ are the inverses of different elements of
$\mathcal{H}_A$ written out as columns with the above notation and
$r=\# \mathcal{H}_A$.


\begin{example}\label{ej}{\normalfont

\[A:\quad
\def\objectstyle{\scriptstyle}
\def\labelstyle{\scriptstyle}
\vcenter{
  \xymatrix@-1.1pc{
          &                                         &  & &v_6\ar[ld]^-{ \alpha_5}&\\
      v_5 & v_4 \ar[l]_-{\alpha_4} &  & v_{7}\ar[ld]^-{}="d"^-{ \alpha_6} & & \\
          & &v_3 \ar[lu]^-{}="e"^-{\alpha_3}\ar[ld]^-{}="g"^-{ \alpha_7} &  & v_{1} \ar[lu]^-{}="c"_-{\phantom{\cdot}}^-{\alpha_9} \ar[ld]^-{}="b"_-{\phantom{\cdot}}^-{\alpha_1}  & v_9 \ar[l]_-{}="a"^-{\alpha_8}  \\
      v_{10}\ar[uu]^-{ \alpha_{10}}\ar[r]^-{\alpha_{11}}   & v_8  & & v_{2} \ar[lu]_-{}="f"^-{\alpha_2}& &       \\
      \ar@{.}"a";"b" \ar@{.} "c";"d" \ar@{.}"e";"d" \ar@{.}"g";"f"
                  }
}
\]

Define $\sigma (\alpha_1)=\sigma (\alpha_2)=\sigma
(\alpha_3)=\sigma (\alpha_4)=\sigma (\alpha_5)=\sigma
(\alpha_6)=\sigma (\alpha_8)=\sigma (\alpha_10)=+1$, $\sigma
(\alpha_7)=\sigma (\alpha_9)=\sigma (\alpha_{11})=-1$,
$\varepsilon (\alpha_6)=\varepsilon (\alpha_8)=\varepsilon
(\alpha_9)=\varepsilon (\alpha_{10})=\varepsilon
(\alpha_{11})=+1$, and $\varepsilon (\alpha_1)=\varepsilon
(\alpha_2)=\varepsilon (\alpha_3)=\varepsilon
(\alpha_4)=\varepsilon (\alpha_5)=\varepsilon (\alpha_7)=-1$.

Then $\mathcal{H}_A$ consists of $ ( v_5^{+} ,v_4^{+} ,v_3^{+}
,v_2^{+} , v_1^{+} )$, $( v_8^{+} , v_3^{-}, v_7^{+} , v_6^{+} )$,
$( v_7^{-} ,v_1 ^{-} , v_9 ^{+} )$, $( v_5 ^{-} , v_{10} ^{+} )$,
$( v_8 ^{-}, v_{10}^{-} )$,$( v_2 ^{-} )$ ,$( v_4 ^{-} )$,$( v_6
^{-} )$y $( v_9 ^{-} )$
 and can be described by the expression

$ \scriptscriptstyle
  \left[
    \begin{matrix}
    v_1^{+} & v_6^{+}  & v_9^{+} & v_{10}^{+} & v_{10}^{-} & v_2^{-} & v_4^{-} & v_6^{-}& v_9^{-}\\
    v_2 ^{+} & v_7^{+} & v_1^{-} & v_5^{-}& v_8^{-} \\
    v_3 ^{+} & v_3 ^{-}& v_7^{-}\\
    v_4 ^{+} & v_8 ^{+}\\
    v_5 ^{+}
    \end{matrix}
  \right]
$

Usually trivial paths will be omitted. In this case we write

$ \scriptscriptstyle
  \left[
    \begin{matrix}
    v_1^{+} & v_6^{+}  & v_9^{+} & v_{10}^{+} & v_{10}^{-} \\
    v_2 ^{+} & v_7^{+} & v_1^{-} & v_5^{-}& v_8^{-} \\
    v_3 ^{+} & v_3 ^{-}& v_7^{-}\\
    v_4 ^{+} & v_8 ^{+}\\
    v_5 ^{+}
    \end{matrix}
  \right]
$

}\end{example}

There are several forms in which we can describe $\mathcal{H}_A$,
according to the order in which we write the elements of this set.

In forthcoming examples and results we use this way of codifying a
gentle algebra. The computer calculations made for this work use
another way of describing a gentle algebra. Consider ${\mathcal
P}(2n)$, the collection of partitions\footnote{A partition of a
natural $m$ is a sequence $\lambda=(\lambda_1,\lambda_2,\dots)$ of
naturals such that $\lambda_1\geq\lambda_2\geq \dots$, with just a
finite number of no zero terms and such that
$\sum_{i}\lambda_i=m$; it is denoted by $\lambda\vdash m$.} of
$2n$, with $n$ the number of vertices of $Q$, that is, ${\mathcal
P}(2n):=\conj{\lambda|\lambda\vdash 2n}$. Denote by ${\mathcal
R}(2n)$ the collection of partitions\footnote{Let $A$ be a set, a
collection $\{A_i\}_i$ of non empty subsets of $A$ is a partition
of $A$ if $A_i\cap A_j=\emptyset $ for all $i\neq j$ and
$\cup_iA_i=A$.} of the set $[2n]=\conj{1,2, \ldots, 2n}$ whose
elements have cardinality $2$, that is
$${\mathcal R}(2n):=\conj{R | R\text{ is a partition of $[2n]$, $\#p=2$  $\forall p \in R$} }$$

We can define an equivalence relation in ${\mathcal P}(2n)\times
{\mathcal R}(2n)$ and a bijection $\Phi :{\mathcal
A}(n)\rightarrow{\mathcal P}(2n)\times {\mathcal R}(2n)/ \sim $
where ${\mathcal A}(n)$ is the collection of gentle algebras whose
associated quiver is connected and has $n$ vertices. For each
$A=\mathrm{k}Q/ \left< \mathcal{P} \right> \in {\mathcal A}(n)$
let $H_1,H_2,\ldots , H_r$ be the different threads of $A$ ordered
in a non growing form according to its length, define $\lambda:=
(l(H_1)+1,l(H_2)+1,\ldots ,l(H_r)+1)$. Let $Q_0=\{ v_1,v_2,\ldots
,v_n \}$, $\sigma$ and $\varepsilon$ be as described before,
define $f:\{ v_i^{\epsilon}| 1\leq i\leq n, \epsilon \in
\{+1,-1\}\}\to [2n]$ as $f(v_i^{\epsilon}):=(l(H_1)+1)+(l(H_2)+1)+
\dots (l(H_j)+1)-m+1$ if $v_i^{\epsilon}$ is the $m$ element of
vector $H_j$. Let $\gamma \in {\mathcal R}(2n)$ be such that $\{
c,d\}\in \gamma$ if and only if there exists $v_i \in Q_0$ such
that $f(v_i^{\epsilon})=c$ and $f(v_i^{-\epsilon})=d$. Then
$\Phi(A)$ is the equivalence class of pair $(\lambda,\gamma)$.

\begin{example}{\normalfont
For $A$ as in Example \ref{ej}, $\lambda=(5,4,3,2,2,1,1,1,1)$.
Notice $f( v_1^{+})=1$, $f( v_2^{+})=2 $, $f( v_3^{+})=3$, $f(
v_4^{+})=4 $, $f( v_5^{+})=5$, $f( v_6^{+})=6 $, $f( v_7^{+})=7$,
$f( v_3^{-})=8 $, $f( v_8^{+})=9$, $f( v_9^{+})=10$, $f(
v_1^{-})=11$, $f( v_7^{-})=12$, $f( v_{10}^{+})=13$,$f(
v_5^{-})=14$, $f( v_{10}^{-})=15$, $f( v_8^{-})=16$, $f(
v_2^{-})=17$, $f( v_4^{-})=18$, $f( v_6^{-})=19$, $f(
v_9^{-})=20$. Then

$$\small{\gamma=\{\{1,11\},\{2,17\},\{3,8\},\{4,18\},\{5,14\},\{6,19\},\{7,12\},\{9,16\},\{10,20\},\{13,15\}\}}.$$

}\end{example}

\subsection{Previous results}

For a gentle algebra $A$ there exists a derived equivalent
invariant $\phi_A:\mathbb{N}^2\to\mathbb{N}$, see \cite{ag07} for
a precise definition, which can be determined easily if $A$ is
given as a quiver with relations ${\mathrm k}Q/ \left< \mathcal{P}
\right>$ as mentioned in Section \ref{intro}. As $\phi_A$
describes the action of the suspension functor
$\Omega^{-1}_{\hat{A}}$ for the triangulated category
$\hat{A}$-$\underline{\operatorname{mod}}$ on the Auslander-Reiten
components which contain string modules and Auslander-Reiten
triangles, see \cite{ha88}, of the form $X\to Y\to
\tau^{-1}_{\hat{A}}X\to \Omega^{-1}_{\hat{A}}X$ with $Y$
indecomposable, we have the following result, see \cite{ag07}.\\

 {\bfseries Theorem A.}\label{soninv} {\em Let $A$ and $B$ be
gentle algebras. If $A$ and $B$ are derived equivalent then
$\phi_A=\phi_B$.}

\vspace{0.5 cm}

Also

\vspace{0.5 cm}

{\bfseries Theorem C.}\label{unciclo} {\em Let $A={\mathrm k}Q/
\left< \mathcal{P} \right>$ and $B={\mathrm k}Q'/ \left<
\mathcal{P}' \right>$ be gentle algebras such that $c(Q),c(Q')\leq
1$. Then $A$ and $B$ are derived equivalent if and only if
$\phi_A=\phi_B$. }

\vspace{0.5 cm}

As defined in Section \ref{intro} let $\# \phi_A :=\sum _{1\leq j
\leq k } \phi_A(n_j,m_j)$ where $[(n_1,m_1),\dots ,(n_k,m_k)]$
describes $\phi_A$: each $(n_j,m_j)$ in the support of $\phi_A$ is
written $\phi_A(n_j,m_j)$ times and the order in which they are
written is arbitrary.

In this article we prove that if $A$ is a gentle algebra with two
cycles, then $\#\phi_A
 \in \{ 1,3\}$ and we will use these invariants to classify those
kind of algebras in the case where $\#\phi_A =3$.

\section{Normal forms}\label{representantes}

\subsection{Three pairs of natural numbers}

If $A$ is an algebra defined by a quiver with two cycles and
$\phi_A=[(a_1,a_2),(b_1,b_2),(c_1,c_2)]$, we now
$a_1+b_1+c_1=\#\mathcal{H}_A$ and $a_2+b_2+c_2=\# Q_1$. Also, the
permitted threads considered as disjoint graphs form a forrest
with $\#\mathcal{H}_A$ trees and $2\#Q_0$ vertices, so
$\#\mathcal{H}_A=2\#Q_0-\#Q_1$. By definition of $c(Q)$, $\#Q_0 =
\#Q_1-c(Q)+1$, so $\#\mathcal{H}_A=2(\#Q_1-c(Q)+1) -
\#Q_1=\#Q_1-2c(Q)+2$, that is $\#\mathcal{H}_A+2c(Q)-2=\#Q_1$,
equivalently $a_1+b_1+c_1+2(2)-2=a_2+b_2+c_2$ or
$a_1+b_1+c_1+2=a_2+b_2+c_2$. Now we prove the converse, if
$(a_1,a_2)$, $(b_1,b_2)$ and $(c_1,c_2)$ are pairs such that
$a_1+b_1+c_1+2=a_2+b_2+c_2$ there is a gentle algebra $A$ for
which $\phi_A=[(a_1,a_2),(b_1,b_2),(c_1,c_2)]$. We present a list
of normal forms of all different derived equivalence classes in
this case.
\begin{theorem}
Let $[(a_1,a_2),(b_1,b_2),(c_1,c_2)]\subset \mathbb{N}^2$ be such
that $a_1+b_1+c_1+2=a_2+b_2+c_2$. There is a gentle algebra $A$
such that $\phi_A=[(a_1,a_2),(b_1,b_2),(c_1,c_2)]$.
\end{theorem}
\proof

We construct a gentle algebra $A$ for each possible collection of
invariants. First we work with pairs of natural numbers $(0,a)$,
$(a,0)$ or $(1,1)$.
\begin{enumerate}
\item A gentle algebra $A$ such that
$\phi_A=[(0,a),(0,b),(a+b-2,0)]$  with $a,b\geq 1$ and $a+b>2$ is

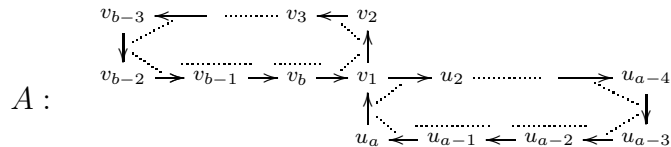
\[A:\quad
\def\objectstyle{\scriptstyle}
\def\labelstyle{\scriptstyle}
\vcenter{
  \xymatrix@-1.1pc{
      v_{b-3} \ar[d]^-{}="a"   &  \ar[l] ^-{}="x"  & v_3  \ar@{.}[l]  & v_{2} \ar[l]^-{}="r"   \\
      v_{b-2} \ar[r]^-{}="b"\ar[r]^-{\phantom{\cdot}}="p" & v_{b-1} \ar[r]^-{\phantom{\cdot}}="c" &    v_{b} \ar[r]^-{\phantom{\cdot}}="d"^-{}="e"  & v_1 \ar[u]_-{}="f" \ar[r]^-{}="i" & u_{2} \ar @{.} [r]  & \ar[r]  ^-{}="l"    & u_{a-4} \ar[d] ^-{}="m" \\
                 &                  &                      & u_{a} \ar[u] ^-{}="h" & u_{a-1} \ar[l]_-{\phantom{\cdot}}="u" ^-{}="v" & u_{a-2} \ar[l]_-{\phantom{\cdot}}="t" & u_{a-3} \ar[l] ^-{}="w"_-{\phantom{\cdot}}="s" \\
      \ar@{.}"a";"b" \ar@{.} "c";"d" \ar@{.}"r";"f" \ar@{.}"h";"i" \ar@{.}"l";"m"\ar@{.}"p";"c" \ar@{.}"u";"v"\ar@{.}"s";"t"\ar@{.}"a";"x" \ar@{.}"e";"f"\ar@{.}"h";"v" \ar@{.}"u";"t" \ar@{.}"m";"w"
                  }
}
\]

\item A gentle algebra $A$ such that
$\phi_A=[(a,0),(b,0),(0,a+b+2)]$ with $a,b\geq 1$ is

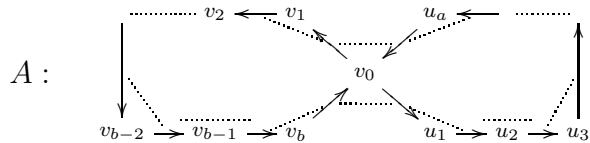
\[A:\quad
\def\objectstyle{\scriptstyle}
\def\labelstyle{\scriptstyle}
\vcenter{
  \xymatrix@-1.1pc{
    \ar[dd]^-{}="p" & v_2 \ar@{.}[l] & v_1 \ar[l]^-{}="j" & & u_{a} \ar[dl]^-{}="l" & \ar[l]^-{}="m" & \ar@{.}[l]  \\
  & & & v_0 \ar[ul]^-{}="k" \ar[dr]^-{}="v"\\
  v_{b-2}\ar[r]^-{}="q"^-{\phantom{\cdot}}="r" & v_{b-1} \ar[r]^-{\phantom{\cdot}}="s"^-{}="t" & v_b \ar[ur]^-{}="u" & & u_1 \ar[r]^-{\phantom{\cdot}}="a"_-{}="d"   & u_2 \ar[r]^-{}="n"^-{\phantom{\cdot}}="b" & u_{3} \ar[uu]^-{}="o"
      \ar@{.}"a";"b"  \ar@{.}"j";"k" \ar@{.}"l";"k" \ar@{.}"m";"l" \ar@{.}"o";"n" \ar@{.}"p";"q" \ar@{.}"r";"s" \ar@{.}"t";"u" \ar@{.}"u";"v"\ar@{.}"d";"v"
                  }
}
\]
\item A gentle algebra $A$ such that
$\phi_A=[(1,1),(b,0),(0,b+2)]$ is
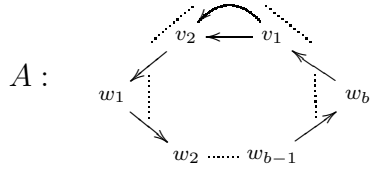
\[A:\quad
\def\objectstyle{\scriptstyle}
\def\labelstyle{\scriptstyle}
\vcenter{
  \xymatrix@-1.1pc{
    & v_2 \ar[dl]_-{\phantom{\cdot}}="b"^-{}="g" &  v_1 \ar[l] \ar@/_1pc/[l]_(.3){\phantom{\cdot}}="d"_(.7){\phantom{\cdot}}="a"  \\
   w_1 \ar[dr]^-{}="e" & & & w_{b} \ar[ul]_-{\phantom{\cdot}}="c"^-{}="h"\\
   & w_2 \ar@{.}[r]  & w_{b-1} \ar[ur]^-{}="f"
      \ar@{.}"a";"b"  \ar@{.}"c";"d" \ar@{.}"g";"e" \ar@{.}"f";"h"
                  }
}
\]
\item A gentle algebra $A$ such that $\phi_A=[(1,1),(1,1),(0,2)]$
is
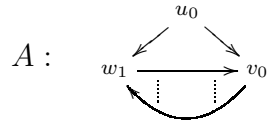
\[A:\quad
\def\objectstyle{\scriptstyle}
\def\labelstyle{\scriptstyle}
\vcenter{
  \xymatrix@-1.1pc{
     & u_0 \ar[dl] \ar[dr] \\
   w_1 \ar[rr]^(.3){}="a"^(.7){}="c" &  & v_0  \ar@/ ^1.5pc/[ll]^(.3){}="d"^(.7){}="b"
      \ar@{.}"a";"b" \ar@{.} "c";"d"
                  }
}
\]

\item A gentle algebra $A$ such that $\phi_A=[(1,1),(0,1),(0,1)]$
is
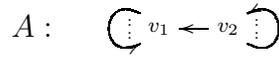
\[A:\quad
\def\objectstyle{\scriptstyle}
\def\labelstyle{\scriptstyle}
\vcenter{
  \xymatrix@-1.1pc{
     & v_1 \ar@(ul,dl)^(.2){}="a"^(.8){}="b" \ar@{}[l]|{_:^:} & v_{2}  \ar[l] \ar@{}[r]|{_:^:} \ar@(dr,ur)^(.2){}="c"^(.3){}="d"  &
                  }
}
\]
For the general cases some trivial permitted threads are added to
the previous algebras:

\item A gentle algebra $A$ such that
$\phi_A=[(k,a+k),(q,b+q),(a+b-2+r,r)]$ with $a\geq b\geq 1$,
$a+b>2$, $k\geq q$ if $a=b$, is
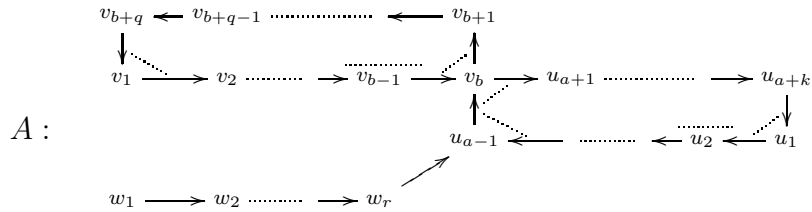
\[A:\quad
\def\objectstyle{\scriptstyle}
\def\labelstyle{\scriptstyle}
\vcenter{
  \xymatrix@-1.1pc{
      v_{b+q} \ar[d]^-{}="a"   & v_{b+q-1} \ar[l]  & & \ar@{.}[ll]  & v_{b+1} \ar[l]   \\
      v_1 \ar[r]^-{}="b" & v_2 \ar@{.}[r] &\ar[r]^-{\phantom{\cdot}}="c" &    v_{b-1} \ar[r]^-{\phantom{\cdot}}="d"^-{}="e"  & v_b \ar[u]_-{}="f" \ar[r]^-{}="i" & u_{a+1} \ar @{.} [rr] &              & \ar[r]     & u_{a+k} \ar[d] ^-{}="j" \\
                 &                  &        &                   & u_{a-1} \ar[u] ^-{}="h"& \ar[l]  ^-{}="g" & \ar @{.} [l] & u_2 \ar[l]_-{\phantom{\cdot}}="m"  & u_1 \ar[l] ^-{}="k"_-{\phantom{\cdot}}="l" \\
      w_1 \ar[r] & w_2 \ar@{.}[r] & \ar[r] & w_r \ar[ur] \\
      \ar@{.}"a";"b" \ar@{.} "c";"d" \ar@{.}"e";"f" \ar@{.}"g";"h" \ar@{.}"h";"i" \ar@{.}"j";"k" \ar@{.}"l";"m"
                  }
}
\]
\item A gentle algebra $A$ such that
$\phi_A=[(a+k,k),(b+q,q),(r,a+b+2+r)]$ with $a\geq b\geq 1$,
$k\geq q$ if $a=b$ is
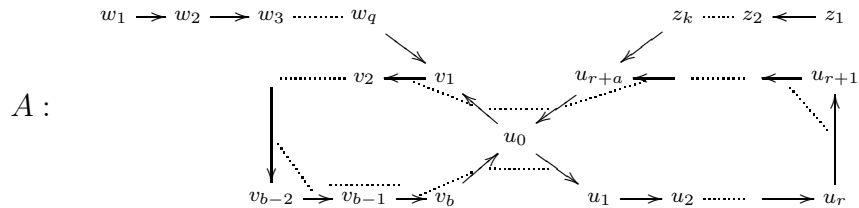
\[A:\quad
\def\objectstyle{\scriptstyle}
\def\labelstyle{\scriptstyle}
\vcenter{
  \xymatrix@-1.1pc{
      w_1 \ar[r]_-{\phantom{\cdot}}="a" & w_2 \ar[r]_-{\phantom{\cdot}}="b" & w_3 \ar@{.}[r] _-{\phantom{\cdot}}="c"^-{}="d" & w_{q} \ar[dr]^-{}="e" & & & & z_k\ar[dl]^-{}="f" & z_2 \ar@{.}[l]^-{}="g"^-{\phantom{\cdot}}="h" & z_1\ar[l]^-{\phantom{\cdot}}="i" \\
      & & \ar[dd]^-{}="p" & v_2 \ar@{.}[l] & v_1 \ar[l]^-{}="j" & & u_{r+a} \ar[dl]^-{}="l" & \ar[l]^-{}="m" & \ar@{.}[l] & u_{r+1} \ar[l]^-{}="n"\\
 & & & & & u_0 \ar[ul]^-{}="k" \ar[dr]^-{}="v"\\
 & & v_{b-2}\ar[r]^-{}="q"^-{\phantom{\cdot}}="r" & v_{b-1} \ar[r]^-{\phantom{\cdot}}="s"^-{}="t" & v_b \ar[ur]^-{}="u" & & u_1 \ar[r] & u_2 \ar@{.}[r] &
\ar[r] & u_{r} \ar[uu]^-{}="o"
       \ar@{.}"j";"k" \ar@{.}"l";"k" \ar@{.}"m";"l" \ar@{.}"o";"n" \ar@{.}"p";"q" \ar@{.}"r";"s" \ar@{.}"t";"u" \ar@{.}"u";"v"
                  }
}
\]
\item A gentle algebra $A$ such that
$\phi_A=[(k,k),(b+q,q),(r,b+2+r)]$ is
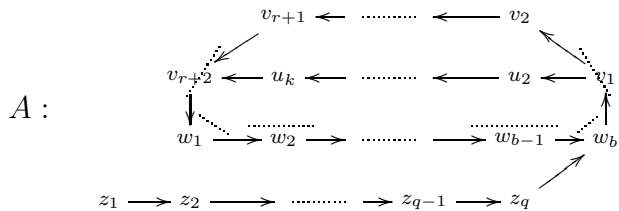
\[A:\quad
\def\objectstyle{\scriptstyle}
\def\labelstyle{\scriptstyle}
\vcenter{
  \xymatrix@-1.1pc{
    &   & v_{r+1} \ar[dl]^-{}="b"_-{\phantom{\cdot}}="k" & \ar[l]^-{}="a" & \ar@{.}[l] & v_2 \ar[l]^-{}="d" \\
    & v_{r+2} \ar[d]^-{}="e"_-{\phantom{\cdot}}="l" & u_k \ar[l]^-{\phantom{\cdot}}="o" & \ar[l]^-{\phantom{\cdot}}="p" & \ar@{.}[l] & u_2 \ar[l]^ -{\phantom{\cdot}}="q" & v_1 \ar[l]^ -{\phantom{\cdot}}="r" \ar[ul]^-{}="c"_-{\phantom{\cdot}}="m" \\
    & w_{1} \ar[r]^-{}="f"^-{\phantom{\cdot}}="s" & w_2 \ar[r]^-{\phantom{\cdot}}="t" & \ar@{.}[r]& \ar[r]^-{\phantom{\cdot}}="u" & w_{b-1} \ar[r]^-{}="g"^-{\phantom{\cdot}}="v" & w_b \ar[u]^-{}="h"_-{\phantom{\cdot}}="n"\\
   z_1 \ar[r]^-{\phantom{\cdot}}="w" & z_2 \ar[r]^-{\phantom{\cdot}}="x" & \ar@{.}[r]   & \ar[r]^-{\phantom{\cdot}}="y" & z_{q-1} \ar[r]^-{}="i"^-{\phantom{\cdot}}="z" & z_q \ar[ur]^-{}="j"
     \ar@{.}"e";"f" \ar@{.}"g";"h" \ar@{.}"k";"l"\ar@{.}"m";"n" \ar@{.}"s";"t"\ar@{.}"u";"v"
                  }
}
\]
\item A gentle algebra $A$ such that
$\phi_A=[(k,k),(q,q),(r,r+2)]$ with $k \leq q$ is
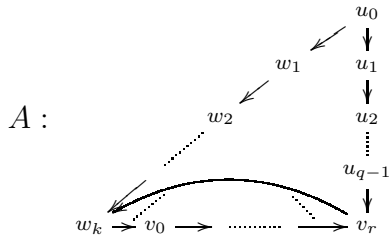
\[A:\;
\def\objectstyle{\scriptstyle}
\def\labelstyle{\scriptstyle}
\vcenter{
  \xymatrix@-1.3pc{
    & &     &     & u_0 \ar[dl] \ar[d]\\
    & &     & w_1 \ar[dl] & u_1 \ar[d] \\
    & & w_2 \ar@{.}[dl] &     & u_2 \ar@{.}[d] \\
    &\ar[dl] &     &     & u_{q-1} \ar[d] \\
 w_{k}\ar[r]^-{}="a" & v_0 \ar[r] & \ar@{.}[r] & \ar[r] ^-{}="c" & v_r \ar@/ _1.5pc/[llll]^(.3){}="d" ^(.7){}="b"
      \ar@{.}"a";"b" \ar@{.}"c";"d"
                  }
}
\]

\item A gentle algebra $A$ such that
$\phi_A=[(k,k),(q,q+1),(r,r+1)]$ with $q\leq r$ is
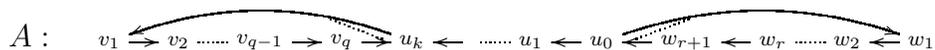
\[A:\quad
\def\objectstyle{\scriptstyle}
\def\labelstyle{\scriptstyle}
\vcenter{
  \xymatrix@-1.1pc{
     v_1 \ar[r] & v_2 \ar@{.}[r] & v_{q-1} \ar[r]_{\phantom{\cdot}}="a" & v_q \ar[r]_{\phantom{\cdot}}="b"^(.7){}="c" & u_{k}  \ar@/_1pc/[llll]^(.3){}="d" &
\ar[l] & u_1 \ar@{.}[l] & u_0 \ar[l] \ar@/^1pc/[rrrr]^(.3){}="f" &
w_{r+1} \ar[l]^(.7){}="e"^-{\phantom{\cdot}}="h" & w_{r}
\ar[l]^-{\phantom{\cdot}}="g" & w_2 \ar@{.}[l] & w_1 \ar[l]
      \ar@{.} "c";"d"\ar@{.}"e";"f"
                  }
}
\]

\end{enumerate}

This covers all possible cases if \# $\phi_A=3$:

Suppose $A$ is a gentle algebra defined by a quiver with two
cycles such that $\phi_A=[(a_1,a_2),(b_1,b_2),(c_1,c_2)]$. If
$a_1=a_2$ and $b_1=b_2$ we have $(9)$. If $a_1=a_2$ but $b_1\ne
b_2$ and $c_1\neq c_2$, one possibility is $b_1> b_2$ and $c_1<
c_2$, or $b_1< b_2$ and $c_1< c_2$, which correspond to $(8)$ and
$(10)$ respectively; as $a_1+b_1+c_1+2=a_2+b_2+c_2$ it is not
possible that $b_1> b_2$ and $c_1> c_2$.

Other possibility is that $a_1\neq a_2$, $b_1\neq b_2$ and
$c_1\neq c_2$ with $a_1< a_2$, $b_1< b_2$ and $c_1> c_2$, or $a_1>
a_2$, $b_1> b_2$ and $c_1< c_2$, which correspond to cases $(6)$
and $(7)$ presented before; as $a_1+b_1+c_1+2=a_2+b_2+c_2$, it is
impossible that $a_1< a_2$, $b_1< b_2$ and $c_1< c_2$, or that
$a_1> a_2$, $b_1> b_2$ and $c_1> c_2$.

\endproof

We prove that any gentle algebra $A$ with two cycles and
$\#\phi_A=3$, is derived equivalent to one of the algebras
presented in this section. In order to do that we will need some
combinatorial transformations over the quiver with relations which
define those algebras.

\section{Description of the elementary transformations in terms
of the threads of $A$}\label{transhilos}

Now, we analyze each one of the combinatorial transformation which
preserve derived equivalence presented in the Appendix by using
the description of a gentle algebra in terms of its threads. These
transformations are very easy to calculate and when described in
terms of its threads we see they are all constructed with a more
elementary combinatorial transformation or its corresponding dual.
The main reason we focus our attention on these transformations is
that for any gentle algebra $A=\mathrm{k}Q/ \left< \mathcal{P}
\right>$, $Q$ with two cycles, $\#\phi_A=3$ and $\#Q_0\geq 5$, all
its derived equivalence class can be obtained by using these
transformations.

\begin{definition}{\normalfont\label{mueveu}
Let $A=\mathrm{k}Q/ \left< \mathcal{P} \right>$ be a gentle
algebra with $Q$ connected and $u,v\in Q_0$ distinct consecutive
vertices. Consider a description of $A$ in terms of its threads:
\[
\mathcal{H}_A:\quad
 \scriptscriptstyle
  \left[
    \begin{matrix}
  &  s(H_{u^+}) & s(H_{u^-}) & s(H_{v^-}) & \\
  &  \vdots &\vdots    &  \vdots & \\
  \dots &  u^+    & u^-       & v^- & \dots \\
  &   v^+     & \vdots   & \vdots \\
  &  \vdots & e(H_{u^-}) & e(H_{v^-})  \\
  &  e(H_{u^+})
    \end{matrix}
  \right]
\]

where $H_{u^-}$ and $H_{v^-}$ are threads of $A$ involving $u$ and
$v$ but not the arrow $(v^+,u^+)$, and $H_{u^+}$ the thread of $A$
which involves such arrow. Denote by $m_{u^+}(\mathcal{H}_A)$ the
corresponding array obtained by removing $u^+$ of its position and
putting it below $v^-$, that is
\[
m_{u^+}(\mathcal{H}_A):\quad
 \scriptscriptstyle
  \left[
    \begin{matrix}
        &  s(H_{u^+}) & s(H_{u^-}) & s(H_{v^-}) & \\
        &  \vdots &\vdots    &  \vdots & \\
  \dots &  v^+      & u^-       & v^- & \dots \\
        &   \vdots& \vdots   & u^+ \\
        &  e(H_{u^+}) & e(H_{u^-}) & \vdots  \\
        &         &          &   e(H_{v^-})
    \end{matrix}
  \right].
\]
We say $u^+$ is moved after $v^-$. The inverse transformation is
denoted by $m_{u^+}^{-1}$, in this case we say $u^+$ is moved
before $v^+$. Define $m_{u^-}$ and its inverse in a similar way.
}\end{definition}

This combinatorial transformations describe the ones presented in
the Appendix as we see next.

\subsection{Transformations over a vertex}
Let $A$ as in Section \ref{vertice}. In all cases,
the transformation over a vertex $i$ corresponds to apply
$m_{i^+}$ followed by $m_{i^-}$.
\begin{enumerate}
\item If $s_1\neq i \neq s_2$. The array associated to
$\mathcal{H}_A$ is
\[
\begin{smallmatrix}
\scriptscriptstyle
  \left[
    \begin{smallmatrix}
         & s(H_{i^+}) & s(H_{i^-}) & s(H_{j_1^-}) & s(H_{j_2^-}) & \\
         & \vdots &\vdots     &  \vdots     & \vdots & \\
   \dots &  p_1^+   & p_2^+       & j_1^-        &  j_2^- & \dots \\
         &  i^+     & i^-        & s_1^+         &  s_2^+\\
         & j_2^+    & j_1^+       & \vdots      &  \vdots\\
         & \vdots & \vdots    & e(H_{j_1^-}) &  e(H_{j_2^-}) \\
         & e(H_{i^+}) & e(H_{i^-})
    \end{smallmatrix}
  \right]
&
\begin{matrix}  m_{i^+} \\ \longmapsto \end{matrix}
&  \scriptscriptstyle
  \left[
   \begin{smallmatrix}
         & s(H_{i^+}) & s(H_{i^-}) & s(H_{j_1^-}) & s(H_{j_2^-}) & \\
         & \vdots &\vdots     &  \vdots     & \vdots & \\
   \dots &  p_1^+   & p_2^+       & j_1^-        &  j_2^- & \dots \\
         &  j_2^+   & i^-        & s_1^+         &  i^+\\
         & \vdots & j_1^+       & \vdots      &  s_2^+\\
         & e(H_{i^+}) & \vdots    & e(H_{j_1^-}) &  \vdots \\
         &        & e(H_{i^-}) &             & e(H_{j_2^-})
    \end{smallmatrix}
  \right]
\end{smallmatrix}
\]
\[
\begin{smallmatrix}
\begin{matrix}  m_{i^-} \\ \longmapsto \end{matrix}
&  \scriptscriptstyle
  \left[
   \begin{smallmatrix}
         & s(H_{i^+}) & s(H_{i^-}) & s(H_{j_1^-}) & s(H_{j_2^-}) & \\
         & \vdots &\vdots     &  \vdots     & \vdots & \\
   \dots &  p_1^+   & p_2^+       & j_1^-        &  j_2^- & \dots \\
         &  j_2^+   & j_1^+       & i^-          &  i^+\\
         & \vdots & \vdots    & s_1^+         &  s_2^+\\
         & e(H_{i^+}) & e(H_{i^-}) & \vdots      &  \vdots \\
         &        &           & e(H_{j_1^-}) & e(H_{j_2^-})
    \end{smallmatrix}
  \right]
\end{smallmatrix}
\]
\item If $s_1=i$ the array $\mathcal{H}_A$ is
\begin{enumerate}
\item \[
\begin{smallmatrix}
\scriptscriptstyle
  \left[
    \begin{smallmatrix}
         & s(H_{i^+}) & s(H_{i^-}) &  s(H_{j_2^-}) & \\
         & \vdots &\vdots     &  \vdots     & \\
   \dots &  j_1^+   & p_2^+       &   j_2^-      & \dots \\
         &  i^+     & i^-        &   s_2^+      \\
         & j_2^+    & j_1^-       &   \vdots   \\
         & \vdots & \vdots    &   e(H_{j_2^-}) \\
         & e(H_{i^+}) & e(H_{i^-})
    \end{smallmatrix}
  \right]
&
\begin{matrix}  m_{i^+} \\ \longmapsto \end{matrix}
&  \scriptscriptstyle
  \left[
   \begin{smallmatrix}
         & s(H_{i^+}) & s(H_{i^-}) &  s(H_{j_2^-}) & \\
         & \vdots &\vdots     &  \vdots & \\
   \dots &  j_1^+   & p_2^+       &  j_2^- & \dots \\
         &  j_2^+   & i^-        &   i^+\\
         & \vdots & j_1^-       &   s_2^+\\
         & e(H_{i^+}) & \vdots    &   \vdots \\
         &        & e(H_{i^-}) &  e(H_{j_2^-})
    \end{smallmatrix}
  \right]
&
\begin{matrix}  m_{i^-} \\ \longmapsto \end{matrix}
&  \scriptscriptstyle
  \left[
   \begin{smallmatrix}
         & s(H_{i^+}) & s(H_{i^-}) &  s(h_{j_2^-}) & \\
         & \vdots &\vdots     &  \vdots & \\
   \dots &  j_1^+  & p_2^+       &  j_2^- & \dots \\
         &  i^-     & j_1^-       &  i^+\\
         & j_2^+    & \vdots    &  s_2^+\\
         & \vdots & e(H_{i^-}) &   \vdots \\
         & e(H_{i^+}) &           &  e(H_{j_2^-})
    \end{smallmatrix}
  \right]
\end{smallmatrix}
\]

or \item \[
\begin{smallmatrix}
\scriptscriptstyle
  \left[
    \begin{smallmatrix}
         & s(H_{i^+}) & s(H_{i^-}) &  s(H_{j_2^-}) & \\
         & \vdots &\vdots     &  \vdots & \\
   \dots &  p_1^+   & j_1^-      &  j_2^- & \dots \\
         &  i^+     & i^-        &   s_2^+\\
         & j_2^+    & j_1^+       &  \vdots\\
         & \vdots & \vdots    &   e(H_{j_2^-}) \\
         & e(H_{i^+}) & e(H_{i^-})
    \end{smallmatrix}
  \right]
&
\begin{matrix}  m_{i^+} \\ \longmapsto \end{matrix}
&  \scriptscriptstyle
  \left[
   \begin{smallmatrix}
         & s(H_{i^+}) & s(H_{i^-}) &  s(H_{j_2^-}) & \\
         & \vdots &\vdots     &  \vdots & \\
   \dots &  p_1^+   & j_1^-      &   j_2^- & \dots \\
         &  j_2^+   & i^-        &   i^+\\
         & \vdots & j_1^+       &   s_2^+\\
         & e(H_{i^+}) & \vdots    &   \vdots \\
         &        & e(H_{i^-}) &  e(H_{j_2^-})
    \end{smallmatrix}
  \right]
\end{smallmatrix}
\]
which remains invariant after applying $m_{i^-}$
\end{enumerate}
\item If $s_1=s_2=i$ the array $\mathcal{H}_A$ is invariant under
$m_{i^+}$ y $m_{i^-}$.
\end{enumerate}
In any case, the resulting arrays correspond to $H_{V_{i}(A)}$.

\subsection{Transformation over an arrow}

Let $A$ be as in Section \ref{flecha}. We prove that the
transformation over an arrow $(j^+,i^+)$ corresponds to the
application of $m_{i^+}$.

The array $\mathcal{H}_A$ is
\[
\begin{smallmatrix}
\scriptscriptstyle
  \left[
    \begin{smallmatrix}
         & s(H_{i^+}) & s(H_{i^-}) & s(H_{j^-})  & \\
         & \vdots &\vdots     &  \vdots    & \\
   \dots &  l^+     & b^+         & j^-         & \dots \\
         &  i^+     & i^-        & x^+         \\
         &  j^+     & c^+         & \vdots    \\
         & \vdots & \vdots    & e(H_{j^-}) \\
         & e(H_{i^+}) & e(H_{i^-})
    \end{smallmatrix}
  \right]
&
\begin{matrix}  m_{i^+} \\ \longmapsto \end{matrix}
&  \scriptscriptstyle
  \left[
   \begin{smallmatrix}
         & s(H_{i^+}) & s(H_{i^-}) &  s(H_{j^-}) & \\
         & \vdots &\vdots     &  \vdots & \\
   \dots &  l^+   & b^+          & j^- & \dots \\
         &  j^+   & i^-        &   i^+\\
         & \vdots & c^+       & x^+ \\
         & e(H_{i^+}) & \vdots    &  \vdots \\
         &        & e(H_{i^-}) &  e(H_{j^-})
    \end{smallmatrix}
  \right]
\end{smallmatrix}
\]
which corresponds to $H_{F_{(j^+,i^+)}(A)}$.

\subsection{Transformation over a loop}

Let $A$ be as in Section \ref{lazo}. We prove that transformation
over a loop $(i^-,i^+)$ corresponds to apply $m_{i^-}$ followed by
$m_{i^+}$. The array $\mathcal{H}_A$ is
\[
\begin{smallmatrix}
\scriptscriptstyle
  \left[
    \begin{smallmatrix}
         & s(H_{i^+}) & s(H_{j^-})  & \\
         & \vdots &  \vdots      & \\
   \dots &  l^+     & j^-         & \dots \\
         &  i^+     & x^+         \\
         &  i^-    & \vdots    \\
         &  j ^+    & e(H_{j^-})  \\
         & \vdots   \\
         & e(H_{i^+})
    \end{smallmatrix}
  \right]
&
\begin{matrix}  m_{i^-} \\ \longmapsto \end{matrix}
&  \scriptscriptstyle
  \left[
   \begin{smallmatrix}
         & s(H_{i^+})  &  s(H_{j^-}) & \\
         & \vdots  &  \vdots    & \\
   \dots &  l^+      &  j^-        & \dots \\
         &  i^+      &  i^- \\
         &  j ^+     &  x^+\\
         & \vdots  &  \vdots \\
         & e(H_{i^+})& e(H_{j^-})
    \end{smallmatrix}
  \right]
&
\begin{matrix}  m_{i^+} \\ \longmapsto \end{matrix}
&  \scriptscriptstyle
  \left[
   \begin{smallmatrix}
         & s(H_{i^+}) &  s(H_{j^-}) & \\
         & \vdots & \vdots & \\
   \dots &  l^+     &  j^- & \dots \\
         &  j ^+    &  i^+\\
         & \vdots &  i^-\\
         & e(H_{i^+}) &  x^+ \\
         &        &  \vdots \\
         &        &  e(H_{j^-})
    \end{smallmatrix}
  \right]
\end{smallmatrix}
\]

which corresponds to $H_{L_{(i^-,i^+)}(A)}$.

\section{Reduction of an algebra to its normal form}\label{reduc}
In this section we study some reductions of gentle algebras with
two cycles which use the elementary transformations mentioned
before, to modify them into some of the normal forms presented in
Section \ref{representantes}. Because of the diversity of the
algebras involved, even for few vertices, we do not present an
algorithm to transform them into its normal form. Else, we develop
induction proofs in which, by knowing the elementary
transformations used to modify a gentle algebra into a normal
form, we can transform the algebra obtained by adding a vertex in
a convenient way also into a normal form. In this section we
present some technical results necessaries in those induction
proofs.

\subsection{Simplification of the branches in a quiver}

In order to present induction proofs over the number of vertices
of the quiver which defines the algebra we need the following:

\begin{definition}{\normalfont\label{grado}
Let $Q$ be a quiver. The {\bfseries degree of a vertex} is the
number of arrows beginning in $v$ plus the number of arrows ending
in $v$. }\end{definition}

\begin{definition}{\normalfont\label{a-x}
Let $A=\mathrm{k}Q/ \left< \mathcal{P} \right>$ be a gentle
algebra, with $Q$ a connected quiver and $x \in Q_0$ a transition
vertex (that is there is just one arrow $\alpha$ which ends in the
vertex, only one arrow $\beta$ which starts in it and
$\beta\alpha\notin \mathcal{P}$), or a vertex of degree one which
is the start point of an arrow; graphically we have one of the
following cases:
\begin{enumerate}
\item
\[
\def\objectstyle{\scriptstyle}
\def\labelstyle{\scriptstyle}
\vcenter{
  \xymatrix@-1.1pc{
 \ar@{-}@/^1pc/[dd] & & & & \ar@{-}@/_1pc/[dd]\\
    & u  \ar[r] ^-{\alpha} & x \ar[r]^-{\beta} & v \\
    & & & &
             }
}
\]
\item
\[
\def\objectstyle{\scriptstyle}
\def\labelstyle{\scriptstyle}
\vcenter{
  \xymatrix@-1.1pc{
  \ar@{-}@/^1pc/[dd] &    &  x \ar[d]_-{\theta} ^-{}="a" & & \ar@{-}@/_1pc/[dd]\\
   & u \ar[r]^-{\alpha}="c" & w \ar[r] _-{\kappa}^-{}="b"  \ar[d]^-{}="d" & z & \\
   &   &  v & &         \\
&  \ar@{-}@/^1pc/[rr] & &     &
   \ar@{.}"a";"b" \ar@{.}"c";"d"
             }
}
\]

\end{enumerate}
Consider the corresponding quiver with relations obtained by
removing vertex $x$, that is, an algebra
$A\setminus\{x\}:=\mathrm{k}Q'/ \left< \mathcal{P}' \right>$ where
$Q_0'=Q_0\setminus\{x\}$, $Q_1'=(Q_1\setminus \{\alpha ,\beta
\})\cup \{\alpha '\}$ with $s(\alpha ')= u$, $e(\alpha ')=v$
$\mathcal{P}'=\mathcal{P}\setminus (\{ \alpha\gamma|\gamma \in
Q_1\}\cup\{\gamma \beta |\gamma \in Q_1\} )\cup \{\alpha
'\gamma|\alpha \gamma\in \mathcal{P}  \}\cup \{\gamma\alpha
'|\gamma\beta\in \mathcal{P}  \}$ for case $(1)$; and
$Q_0'=Q_0\setminus\{x\}$, $Q_1'=(Q_1\setminus \{ \theta \})$
$\mathcal{P}'=\mathcal{P}\setminus\{ \kappa \theta \}$ for case
$(2)$. That is:
\[
\begin{matrix}
\def\objectstyle{\scriptstyle}
\def\labelstyle{\scriptstyle}
\vcenter{
  \xymatrix@-1.1pc{
 \ar@{-}@/^1pc/[dd] & & & & \ar@{-}@/_1pc/[dd]\\
    & u  \ar[r] ^-{\alpha} & x \ar[r]^-{\beta} & v \\
    & & & &
             }
} & \longmapsto &
\def\objectstyle{\scriptstyle}
\def\labelstyle{\scriptstyle}
\vcenter{
  \xymatrix@-1.1pc{
 \ar@{-}@/^1pc/[dd] & & & \ar@{-}@/_1pc/[dd]\\
    & u  \ar[r] ^-{\alpha '} & v \\
    & & & &
             }
}
\end{matrix}
\]

\[
\begin{matrix}
\def\objectstyle{\scriptstyle}
\def\labelstyle{\scriptstyle}
\vcenter{
  \xymatrix@-1.1pc{
  \ar@{-}@/^1pc/[dd] &    &  x \ar[d]_-{\theta} ^-{}="a" & & \ar@{-}@/_1pc/[dd]\\
   & u \ar[r]^-{}="c" & w \ar[r] _-{\kappa}^-{}="b"  \ar[d]^-{}="d" & z & \\
   &   &  v & &         \\
&  \ar@{-}@/^1pc/[rr] & &     &
   \ar@{.}"a";"b" \ar@{.}"c";"d"
             }
} & \longmapsto &
\def\objectstyle{\scriptstyle}
\def\labelstyle{\scriptstyle}
\vcenter{
  \xymatrix@-1.1pc{
  \ar@{-}@/^1pc/[dd] &    &   & & \ar@{-}@/_1pc/[dd]\\
   & u \ar[r]^-{}="c" & w \ar[r] _-{\kappa}^-{}="b"  \ar[d]^-{}="d" & z & \\
   &   &  v & &         \\
&  \ar@{-}@/^1pc/[rr] & &     &
   \ar@{.}"c";"d"
             }
}
\end{matrix}
\]

}\end{definition}

\begin{remark}{\normalfont\label{a-xind}
If $A=\mathrm{k}Q/ \left< \mathcal{P} \right>$ is a gentle
algebra, with $Q$ a quiver with two cycles and $x\in Q_0$ as in
previous definition, $A\setminus \{x \}=\mathrm{k}Q'/ \left<
\mathcal{P}' \right>$ is also a gentle algebra with $Q'$ of two
cycles, because it has exactly one vertex an one arrow less than
$Q$. }\end{remark}

\begin{remark}\label{metorama}{\normalfont
If
\[
\def\objectstyle{\scriptstyle}
\def\labelstyle{\scriptstyle}
\vcenter{
  \xymatrix@-1.1pc{
 \ar@{-}@/^1pc/[dd] &  & x \ar[dl]_-{\delta}_-{}="a" & \ar@{-}@/_1pc/[dd]\\
    & u \ar[r]^-{}="b"& v   &  \\
    & & &
 \ar@{.}"a";"b"
             }
}
\]
after applying the arrow transformation $F_{\delta}$ we get
\[
\def\objectstyle{\scriptstyle}
\def\labelstyle{\scriptstyle}
\vcenter{
  \xymatrix@-1.1pc{
 \ar@{-}@/^1pc/[dd] & & & & \ar@{-}@/_1pc/[dd]\\
    & u \ar[r] & x \ar[r] & v & \\
    & & & &
             }
}
\]
which is a vertex as in Definition \ref{a-x} $(1)$. Similarly in
the dual case. }
\end{remark}

\begin{definition}{\normalfont\label{arbolenraizado}
Consider $Q=(Q_0,Q_1,s,e)$ a finite connected quiver. Let $S$ be
the collection of vertices which belong to a cycle or to a
trajectory which joins two vertices belonging to cycles. We say
$Q'$ a subquiver of $Q$ is a {\bfseries rooted tree in $v\in Q_0$}
if $Q'$ is tree, $v$ is the only vertex of $Q'$ which belongs to
$S$; also $\# Q_1'$ is maximal with respect to this property. The
{\bfseries depth of the rooted tree in $v$} is the number of
arrows in a trajectory which begins in $v$, of maximal length in
the subjacent graph. We say $Q'$ a rooted tree in $v$ is a
{\bfseries branch} if all vertices have degree at most two in
$Q'$. }\end{definition}

Let us study rooted trees of a quiver associated to a gentle
algebra:

\begin{lemma}\label{ramas}
Let $Q$ be a quiver of two cycles and more than five vertices,
associated to a gentle algebra. If there exist rooted trees deeper
than one we have one of the following cases:
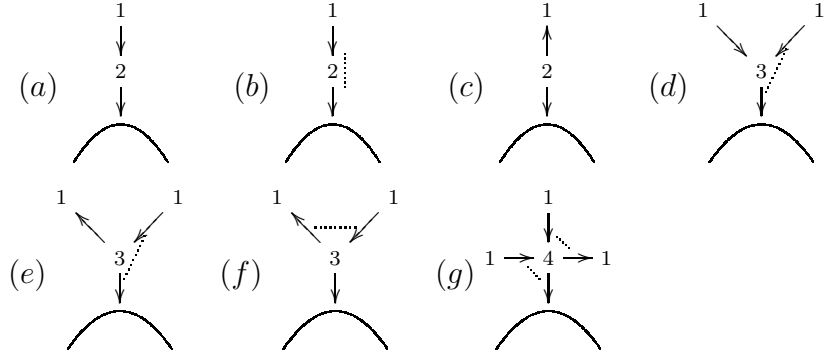
\[
\begin{matrix}
 (a)
\def\objectstyle{\scriptstyle}
\def\labelstyle{\scriptstyle}
\vcenter{
  \xymatrix@-1.1pc{
     & 1 \ar[d]   \\
     & 2 \ar[d]   \\
     &            \\
 \ar@{-}@/^1.5pc/[rr] &      &
                  }
}
 &
 (b)
\def\objectstyle{\scriptstyle}
\def\labelstyle{\scriptstyle}
\vcenter{
  \xymatrix@-1.1pc{
     & 1 \ar[d] ^-{\phantom{\cdot}}="a"  \\
     & 2 \ar[d] ^-{\phantom{\cdot}}="b"  \\
     &            \\
 \ar@{-}@/^1.5pc/[rr] &      &
   \ar@{.}"a";"b"               }
} &
 (c)
\def\objectstyle{\scriptstyle}
\def\labelstyle{\scriptstyle}
\vcenter{
  \xymatrix@-1.1pc{
     & 1    \\
     & 2 \ar[u] \ar[d]   \\
     &            \\
 \ar@{-}@/^1.5pc/[rr] &      &
                  }
} &
 (d)
\def\objectstyle{\scriptstyle}
\def\labelstyle{\scriptstyle}
\vcenter{
  \xymatrix@-1.1pc{
    1 \ar[dr] &   & 1 \ar[dl] ^-{}="a" \\
     & 3 \ar[d] ^-{}="b" & \\
     &            \\
 \ar@{-}@/^1.5pc/[rr] &      &
   \ar@{.}"a";"b"               }
}
\\
 (e)
\def\objectstyle{\scriptstyle}
\def\labelstyle{\scriptstyle}
\vcenter{
  \xymatrix@-1.1pc{
    1  &   & 1 \ar[dl] ^-{}="a" \\
     & 3 \ar[ul] \ar[d] ^-{}="b" & \\
     &            \\
 \ar@{-}@/^1.5pc/[rr] &      &
   \ar@{.}"a";"b"               }
} &
 (f)
\def\objectstyle{\scriptstyle}
\def\labelstyle{\scriptstyle}
\vcenter{
  \xymatrix@-1.1pc{
    1  &   & 1 \ar[dl] ^-{}="a" \\
     & 3 \ar[ul] ^-{}="b"  \ar[d]& \\
     &            \\
 \ar@{-}@/^1.5pc/[rr] &      &
   \ar@{.}"a";"b"               }
} &
 (g)
\def\objectstyle{\scriptstyle}
\def\labelstyle{\scriptstyle}
\vcenter{
  \xymatrix@-1.1pc{
      &  1 \ar[d] ^-{}="a" &\\
   1 \ar[r]^-{}="c" & 4 \ar[r] ^-{}="b"  \ar[d]^-{}="d" & 1 \\
     &            \\
 \ar@{-}@/^1.5pc/[rr] &      &
   \ar@{.}"a";"b" \ar@{.}"c";"d"
             }
} & &
\end{matrix}
\]
or some of their dual. In the previous diagrams the numbers
indicate the degree of each vertex.

\end{lemma}

\proof

Let $u\in Q_0$ be of degree $1$ such that the only neighbor of
$u$, $v$, is of minimal degree. As the algebra is gentle, the
quiver connected and with at least $6$ vertices, the degree of $v$
is $2$, $3$ or $4$. If $v$ has degree $2$, the only possible cases
are $(a)$, $(b)$ or $(c)$, or their respective dual, if it has
degree $3$ the possibilities are $(d)$, $(e)$ or $(f)$ or one of
their dual, and if it has degree $4$, the only possibility is
$(g)$ or its dual.

\endproof

More specifically

\begin{lemma}\label{ramas2}
Let $Q$ be a two cycle quiver with more than five vertices,
associated to a gentle algebra $A$. If there are rooted trees
deeper than one, $A$ is derived equivalent under elementary
transformations to a gentle algebra $B=\mathrm{k}Q'/ \left<
\mathcal{P'} \right>$ such that there is an $x\in Q_0'$ as in
Definition \ref{a-x} $(1)$.
\end{lemma}

\proof We analyze each one of the cases presented in Lemma
\ref{ramas}. For $(a)$ there is nothing to prove and for
$(b)$,$(d)$, $(e)$, $(f)$ and $(g)$ the result follows from Remark
\ref{metorama}. For $(c)$ we have:

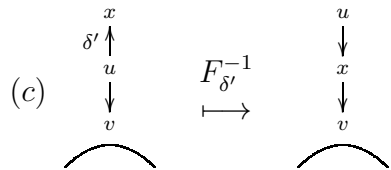
\[
\begin{matrix}
(c)
\def\objectstyle{\scriptstyle}
\def\labelstyle{\scriptstyle}
\vcenter{
  \xymatrix@-1.1pc{
     & x    \\
     & u \ar[u]^-{\delta '} \ar[d]   \\
     & v           \\
 \ar@{-}@/^1pc/[rr] &      &
                  }
} &
\begin{matrix} F_{\delta '}^{-1} \\ \longmapsto \end{matrix}
&
\def\objectstyle{\scriptstyle}
\def\labelstyle{\scriptstyle}
\vcenter{
  \xymatrix@-1.1pc{
     & u \ar[d]   \\
     & x \ar[d]   \\
     & v           \\
 \ar@{-}@/^1pc/[rr] &      &
                  }
}
\end{matrix}
\]
and we apply the dual transformations for the dual cases.
\endproof

Now we analyze what happens with rooted trees of degree one.

\begin{lemma}\label{ramas1}
Let $A=\mathrm{k}Q/ \left< \mathcal{P} \right>$ be a gentle
algebra, with $Q$ a quiver of two cycles. If there is a rooted
tree in $Q$ of degree one, there is an $x \in Q_0$ such that $A$
is derived equivalent under elementary transformations to
$B=\mathrm{k}Q'/ \left< \mathcal{P}' \right>$ where the
corresponding $x$ is presented in one of the following situations:

\begin{enumerate}
\item
\[
\def\objectstyle{\scriptstyle}
\def\labelstyle{\scriptstyle}
\vcenter{
  \xymatrix@-1.1pc{
 \ar@{-}@/^1pc/[dd] & & & & \ar@{-}@/_1pc/[dd]\\
    & u  \ar[r] ^-{\alpha} & x \ar[r]^-{\beta} & v \\
    & & & &
             }
}
\] with $x$ of degree $2$, or
\item
\[
\def\objectstyle{\scriptstyle}
\def\labelstyle{\scriptstyle}
\vcenter{
  \xymatrix@-1.1pc{
 \ar@{-}@/^1pc/[dd] & & x \ar[d] & & \ar@{-}@/_1pc/[dd]\\
    & u  \ar[r] ^-{\alpha}_-{\phantom{\cdot}}="a" & w \ar[r]^-{\beta}_-{\phantom{\cdot}}="b" & v \\
    & & & &
   \ar@{.}"a";"b"
             }
}
\] with $w$ of degree $3$ and $\alpha$ and $\beta$ in a cycle.
\end{enumerate}
\end{lemma}

\proof

We analyze the different cases in which the rooted tree of degree
one can be presented. If the situation is as in Remark
\ref{metorama} there is nothing to prove, we study now what
happens with branches of depth one which look like:

\[
\def\objectstyle{\scriptstyle}
\def\labelstyle{\scriptstyle}
\vcenter{
  \xymatrix@-1.1pc{
 \ar@{-}@/^1pc/[dd] & & x \ar[d] & & \ar@{-}@/_1pc/[dd]\\
    & u  \ar[r] ^-{\alpha}_-{\phantom{\cdot}}="a" & w \ar[r]^-{\beta}_-{\phantom{\cdot}}="b" & v \\
    & & & &
   \ar@{.}"a";"b"
             }
}
\]

with $w$ of degree $3$. If $\alpha$ and $\beta$ belong to a cycle
we have the result. If not, we apply $F_{\beta}$ to obtain

\[
\def\objectstyle{\scriptstyle}
\def\labelstyle{\scriptstyle}
\vcenter{
  \xymatrix@-1.1pc{
 \ar@{-}@/^1pc/[dd] & & x \ar[dr]_-{}="a" & & \ar@{-}@/_1pc/[dd]\\
    & u \ar[r] & w  & v \ar[l]_-{}="b" & \\
    & & & &
 \ar@{.}"a";"b"
             }
}
\]

and then we use Remark \ref{metorama}. In the dual case we apply
the corresponding inverse transformation.

\endproof

Under certain conditions we can assure the existence of a vertex
$x$ for which it makes sense to consider the algebra $A\setminus
\{ x \}$ as in Definition \ref{a-x}.

\begin{proposition}\label{verticex}
Let $A=\mathrm{k}Q/ \left< \mathcal{P} \right>$ be a gentle
algebra, with $Q$ a connected quiver of two cycles and more than
$5$ vertices. Then $A$ is derived equivalent under elementary
transformations to $B=\mathrm{k}Q'/ \left< \mathcal{P}' \right>$
such that there is an $x\in Q_0'$ of the following type:
\begin{enumerate}
\item
\[
\def\objectstyle{\scriptstyle}
\def\labelstyle{\scriptstyle}
\vcenter{
  \xymatrix@-1.1pc{
 \ar@{-}@/^1pc/[dd] & & & & \ar@{-}@/_1pc/[dd]\\
    & u  \ar[r] ^-{\alpha} & x \ar[r]^-{\beta} & v \\
    & & & &
             }
}
\] with $x$ of degree $2$
\item
\[
\def\objectstyle{\scriptstyle}
\def\labelstyle{\scriptstyle}
\vcenter{
  \xymatrix@-1.1pc{
 \ar@{-}@/^1pc/[dd] & & x \ar[d] & & \ar@{-}@/_1pc/[dd]\\
    & u  \ar[r] ^-{\alpha}_-{\phantom{\cdot}}="a" & w \ar[r]^-{\beta}_-{\phantom{\cdot}}="b" & v \\
    & & & &
   \ar@{.}"a";"b"
             }
}
\] with $w$ of degree $3$ and $\alpha$ and $\beta$ in a cycle
\item
\[
\def\objectstyle{\scriptstyle}
\def\labelstyle{\scriptstyle}
\vcenter{
  \xymatrix@-1.1pc{
 \ar@{-}@/^1pc/[dd] & & & & \ar@{-}@/_1pc/[dd]\\
    & u  \ar[r] ^-{\alpha}_-{\phantom{\cdot}}="a" & x \ar[r]^-{\beta}_-{\phantom{\cdot}}="b" & v \\
    & & & &
   \ar@{.}"a";"b"
             }
}
\] with $x$ of degree $2$ and $\alpha$ and $\beta$ in a cycle.
\end{enumerate}
\end{proposition}

\proof

If there are rooted trees the result follows from Lemmas
\ref{ramas2} or \ref{ramas1}.

If there are no rooted trees all vertices have degree $2$, $3$ or
$4$. We know $$\sum _{v\in Q_0} gr(v)=2\#Q_1=2(\#
Q_0+1)=2\#Q_0+2.$$

Then all vertices must have degree $2$, except one of degree $4$
or two of degree $3$. By hypothesis there are at least $6$
vertices, so there are at least two consecutive vertices of degree
$2$ denote them by $u$ and $v$. If situation $(1)$ is not
presented in any of these vertices we have one of the following
situations:
$$
\begin{matrix}
\def\objectstyle{\scriptstyle}
\def\labelstyle{\scriptstyle}
\vcenter{
  \xymatrix@-1.1pc{
 \ar@{-}@/^1pc/[dd] & & & & & \ar@{-}@/_1pc/[dd]\\
    & & u  \ar[l] \ar[r] & v  & \ar[l] \\
    & & & & &
             }
}
&
,
&
\def\objectstyle{\scriptstyle}
\def\labelstyle{\scriptstyle}
\vcenter{
  \xymatrix@-1.1pc{
 \ar@{-}@/^1pc/[dd] & & & & & \ar@{-}@/_1pc/[dd]\\
    &\ar[r]_-{\phantom{\cdot}}="c" & u  \ar[r] ^-{\delta}_-{\phantom{\cdot}}="a" & v \ar[r]_-{\phantom{\cdot}}="b" &  \\
    & & & & &
   \ar@{.}"a";"b" \ar@{.}"c";"a"
             ,}
}
&
,
&
\def\objectstyle{\scriptstyle}
\def\labelstyle{\scriptstyle}
\vcenter{
  \xymatrix@-1.1pc{
 \ar@{-}@/^1pc/[dd] & & & & & \ar@{-}@/_1pc/[dd]\\
    &\ar[r]_-{\phantom{\cdot}}="c" & u  \ar[r]^-{\delta} _-{\phantom{\cdot}}="a" & v  & \ar[l]  \\
    & & & & &
   \ar@{.}"c";"a"
             }
}
\end{matrix}
$$
or the dual of the last one.
In the first situation we apply a vertex transformation
 \[
\begin{matrix}
\def\objectstyle{\scriptstyle}
\def\labelstyle{\scriptstyle}
\vcenter{
  \xymatrix@-1.1pc{
 \ar@{-}@/^1pc/[dd] & & & & & \ar@{-}@/_1pc/[dd]\\
    & & u  \ar[l] \ar[r] & v  & \ar[l] \\
    & & & & &
             }
} &
\begin{matrix} V_{u} \\ \longmapsto \end{matrix}
&
\def\objectstyle{\scriptstyle}
\def\labelstyle{\scriptstyle}
\vcenter{
  \xymatrix@-1.1pc{
 \ar@{-}@/^1pc/[dd] & & & & \ar@{-}@/_1pc/[dd]\\
    & u  & v  \ar[l] &  \ar[l]\\
    & & & &
             }
}
\end{matrix}
\]
where the degree of vertex $u$ can change but the degree of $v$ is
$2$, so we have $(1)$. For the remaining cases, if the arrow
$\delta$ belongs to a cycle we have $(3)$; if not, we can apply
the transformation over that arrow:
\[
\begin{matrix}
\def\objectstyle{\scriptstyle}
\def\labelstyle{\scriptstyle}
\vcenter{
  \xymatrix@-1.1pc{
 \ar@{-}@/^1pc/[dd] & & & & & \ar@{-}@/_1pc/[dd]\\
    &\ar[r]_-{\phantom{\cdot}}="c" & u  \ar[r]^-{\delta} _-{\phantom{\cdot}}="a" & v \ar[r]_-{\phantom{\cdot}}="b" &  \\
    & & & & &
   \ar@{.}"a";"b" \ar@{.}"c";"a"
             }
} &
\begin{matrix} F_{\delta} \\ \longmapsto \end{matrix}
&

\def\objectstyle{\scriptstyle}
\def\labelstyle{\scriptstyle}
\vcenter{
  \xymatrix@-1.1pc{
 \ar@{-}@/^1pc/[dd] & & v \ar[d] & & \ar@{-}@/_1pc/[dd]\\
    &   \ar[r] ^-{\alpha}_-{\phantom{\cdot}}="a" & u \ar[r]^-{\beta}_-{\phantom{\cdot}}="b" &  \\
    & & & &
   \ar@{.}"a";"b"
             }
}
\end{matrix}
\]
and result follows by Lemma \ref{ramas1}, else
\[
\begin{matrix}
\def\objectstyle{\scriptstyle}
\def\labelstyle{\scriptstyle}
\vcenter{
  \xymatrix@-1.1pc{
 \ar@{-}@/^1pc/[dd] & & & & & \ar@{-}@/_1pc/[dd]\\
    &\ar[r]_-{\phantom{\cdot}}="c" & u  \ar[r]^-{\delta} _-{\phantom{\cdot}}="a" & v  & \ar[l]  \\
    & & & & &
   \ar@{.}"c";"a"
             }
} &
\begin{matrix} F_{\delta} \\ \longmapsto \end{matrix}
&
\def\objectstyle{\scriptstyle}
\def\labelstyle{\scriptstyle}
\vcenter{
  \xymatrix@-1.1pc{
 \ar@{-}@/^1pc/[dd] & & & & & \ar@{-}@/_1pc/[dd]\\
    & \ar[r] & u  & v \ar[l] & \ar[l] \\
    & & & & & .
             }
}
\end{matrix}
\]

For the dual case we apply the corresponding inverse
transformation to produce one of the situations already analyzed.

\endproof

Much more:

\begin{proposition}\label{irreducibles}
Let $A=\mathrm{k}Q/ \left< \mathcal{P} \right>$ be a gentle
algebra, with $Q$ a quiver of two cycles and more than $5$
vertices. If $A$ is not derived equivalent to an algebra
$B=\mathrm{k}Q'/ \left< \mathcal{P}' \right>$ with an $x\in Q_0'$
as follows:
\begin{enumerate}
\item
\[
\def\objectstyle{\scriptstyle}
\def\labelstyle{\scriptstyle}
\vcenter{
  \xymatrix@-1.1pc{
 \ar@{-}@/^1pc/[dd] & & & & \ar@{-}@/_1pc/[dd]\\
    & u  \ar[r] ^-{\alpha} & x \ar[r]^-{\beta} & v \\
    & & & &
             }
}
\] $x$ of degree $2$
\item
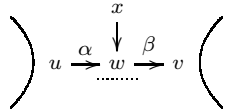
\[
\def\objectstyle{\scriptstyle}
\def\labelstyle{\scriptstyle}
\vcenter{
  \xymatrix@-1.1pc{
 \ar@{-}@/^1pc/[dd] & & x \ar[d] & & \ar@{-}@/_1pc/[dd]\\
    & u  \ar[r] ^-{\alpha}_-{\phantom{\cdot}}="a" & w \ar[r]^-{\beta}_-{\phantom{\cdot}}="b" & v \\
    & & & &
   \ar@{.}"a";"b"
             }
}
\] $w$ of degree $3$ and $\alpha$ and $\beta$ in a cycle
\end{enumerate}
then $\# \phi_A =1$ or $A$ is derived equivalent to one of the
representatives $(1)$, $(2)$ or $(3)$ of Section
\ref{representantes}.
\end{proposition}

\proof

Suppose $A$ is not derived equivalent to an algebra whose
associated quiver has a vertex $x$ as described, we know then by
Lemmas \ref{ramas2} and \ref{ramas1} that $Q$ has no branches.
Also $\# Q_0 \geq 6$ and by the proof of Proposition
\ref{verticex} there must exist a sequence of $3$ arrows of type
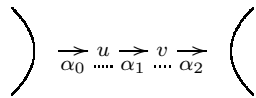
\[
\def\objectstyle{\scriptstyle}
\def\labelstyle{\scriptstyle}
\vcenter{
  \xymatrix@-1.1pc{
 \ar@{-}@/^1pc/[dd] & & & & & \ar@{-}@/_1pc/[dd]\\
    &\ar[r]_-{\alpha_0}="c" & u  \ar[r] _-{\alpha_1}="a" & v \ar[r]_-{\alpha_2}="b" &  \\
    & & & & &
   \ar@{.}"a";"b" \ar@{.}"c";"a"
             }
}
\]
with $u$ and $v$ of degree $2$, $\alpha_0$, $\alpha_1$ and
$\alpha_2$ belonging to the same cycle; moreover, any sequence
consisting of consecutive vertices of degree two $u_1$, $u_2$,
\ldots $u_n$, $n\geq 2$ must be like:
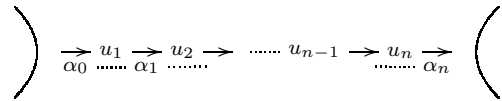
\[
\def\objectstyle{\scriptstyle}
\def\labelstyle{\scriptstyle}
\vcenter{
  \xymatrix@-1.1pc{
 \ar@{-}@/^1pc/[dd] & & & & & & & & \ar@{-}@/_1pc/[dd]\\
    &\ar[r]_-{\alpha_0}="c" & u_1  \ar[r] _-{\alpha_1}="a" & u_2 \ar[r]_-{\phantom{\cdot}}="b" & \ar@{.}[r] & u_{n-1} \ar[r]_-{\phantom{\cdot}}="d" & u_n \ar[r]_-{\alpha_n}="e"& &\\
    & & & & & & & &
   \ar@{.}"a";"b" \ar@{.}"c";"a" \ar@{.}"d";"e"
             }
}
\]

and is constructed by arrows which belong to the same cycle, if
not, by the proof of Proposition \ref{verticex} we would be able
to transform the quiver with relations into another with a vertex
$x$ as described. We also know every vertex in $A$ has degree two,
except one of degree $4$ or two of degree $3$. If there is a $4$
degree vertex $v$, it can be presented in one of the following
ways:

$
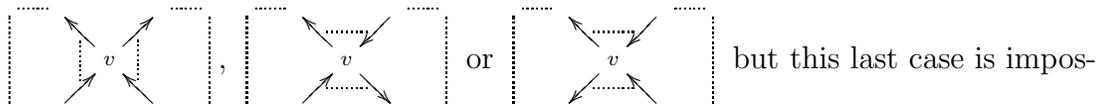
\def\objectstyle{\scriptstyle}
\def\labelstyle{\scriptstyle}
\vcenter{
  \xymatrix@-1.1pc{
  \ar@{.} [dd]\ar@{.}\ar@{.}[r] & & & \ar@{.}[r] & \ar@{.}[dd] \\
    & & v \ar[ul]_-{}="b" \ar[ur] _-{}="c"& & \\
    \ar@{.}[r]& \ar[ur] _-{}="a"& & \ar[ul] _-{}="d"&  \ar@{.}[l]
   \ar@{.}"a";"b" \ar@{.}"c";"d"
             }
} $, $
\def\objectstyle{\scriptstyle}
\def\labelstyle{\scriptstyle}
\vcenter{
  \xymatrix@-1.1pc{
  \ar@{.} [dd]\ar@{.}\ar@{.}[r] & & & \ar[dl] _-{}="a" \ar@{.}[r] & \ar@{.}[dd] \\
    & & v \ar[ul]_-{}="b" \ar[dr] _-{}="c"& & \\
    \ar@{.}[r]& \ar[ur] _-{}="d"& & &  \ar@{.}[l]
   \ar@{.}"a";"b" \ar@{.}"c";"d"
             }
} $ or $
\def\objectstyle{\scriptstyle}
\def\labelstyle{\scriptstyle}
\vcenter{
  \xymatrix@-1.1pc{
  \ar@{.} [dd]\ar@{.}\ar@{.}[r] & & & \ar[dl] _-{}="a" \ar@{.}[r] & \ar@{.}[dd] \\
    & & v \ar[ul]_-{}="b" \ar[dl] _-{}="c"& & \\
    \ar@{.}[r]& & &\ar[ul] _-{}="d" &  \ar@{.}[l]
   \ar@{.}"a";"b" \ar@{.}"c";"d"
             }
} $ but this last case is impossible because in one of the two
cycles there should be more than three arrows and by their
orientation there would appear a sequence different from the ones
mentioned before. In the first case we have a representant of type
$(1)$ and in the second case one of type $(2)$. We study now what
happens if there are two vertices of degree $3$, $a$ and $b$:

By previous analysis all arrows but one or two are part of a
cycle. In any case, by the conditions of the quiver with
relations, $\# \phi_A =1$ or $A$ is derived equivalent to one of
representant $(1)$, $(2)$ or $(3)$ of Section
\ref{representantes}:

Case 1.- All arrows belong to some cycle.

We know there must be a sequence of three arrows at least as the
one mentioned before and its end terms must be $a$ and $b$. Then
one of the next situations is presented:

$ i)
\def\objectstyle{\scriptstyle}
\def\labelstyle{\scriptstyle}
\vcenter
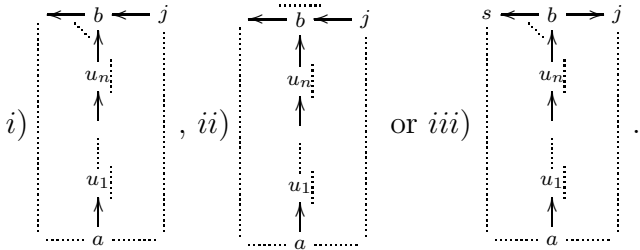
{
  \xymatrix@-1.1pc{
 \ar@{.}[dddd]  & b \ar[l] _-{}="a" & j \ar[l]\ar@{.}[dddd] \\
    & u_n\ar[u] _-{}="b"_-{\phantom{\cdot}}="d" &   \\
   & \ar[u] _-{\phantom{\cdot}}="c" &         \\
   & u_1\ar@{.}[u]_-{\phantom{\cdot}}="f"  & \\
   & a \ar[u]_-{\phantom{\cdot}}="e"   \ar@{.}[l] \ar@{.}[r]    &
   \ar@{.}"a";"b" \ar@{.}"c";"d"  \ar@{.}"e";"f"
             }
} $, $ ii)
\def\objectstyle{\scriptstyle}
\def\labelstyle{\scriptstyle}
\vcenter{
  \xymatrix@-1.1pc{
 \ar@{.}[dddd]   & b \ar[l] _-{\phantom{\cdot}}="a" & j \ar[l]_-{\phantom{\cdot}}="b"\ar@{.}[dddd] \\
   & u_n\ar[u] _-{\phantom{\cdot}}="d" &   \\
   & \ar[u] _-{\phantom{\cdot}}="c" &         \\
   & u_1\ar@{.}[u]_-{\phantom{\cdot}}="f"  & \\
   & a \ar[u]_-{\phantom{\cdot}}="e"   \ar@{.}[l] \ar@{.}[r]    &
   \ar@{.}"a";"b" \ar@{.}"c";"d" \ar@{.}"e";"f"
             }
} $ or $ iii)
\def\objectstyle{\scriptstyle}
\def\labelstyle{\scriptstyle}
\vcenter{
  \xymatrix@-1.1pc{
 s \ar@{.}[dddd]  & b \ar[l] _-{}="a"\ar[r] & j \ar@{.}[dddd] \\
    & u_n\ar[u] _-{}="b"_-{\phantom{\cdot}}="d" &   \\
   & \ar[u] _-{\phantom{\cdot}}="c" &         \\
   & u_1\ar@{.}[u]_-{\phantom{\cdot}}="e"  & \\
   & a \ar[u]_-{\phantom{\cdot}}="f"   \ar@{.}[l] \ar@{.}[r]    &
   \ar@{.}"a";"b" \ar@{.}"c";"d"  \ar@{.}"e";"f"
             }
} $.

$i)$ we get a trivial permitted thread in $u_n$ by applying
$V_b^{-1}$.

$ii)$ If $j\neq a$, we obtain a trivial permitted thread in $j$
after applying $V_b^{-1}$; also, if $j=a$ we get

$
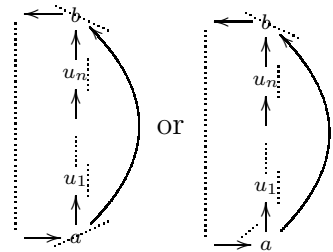
\def\objectstyle{\scriptstyle}
\def\labelstyle{\scriptstyle}
\vcenter{
  \xymatrix@-1.1pc{
 \ar@{.}[dddd]  & b \ar[l] _-{\phantom{\cdot}}="a" &  \\
    & u_n\ar[u]_-{\phantom{\cdot}}="d" &   \\
   & \ar[u] _-{\phantom{\cdot}}="c" &         \\
   & u_1\ar@{.}[u]_-{\phantom{\cdot}}="g"  & \\
    \ar[r] _-{\phantom{\cdot}}="e"& a \ar[u]_-{\phantom{\cdot}}="h"   \ar@{-<}@/_2pc/[uuuu]_(.85){\phantom{\cdot}}="b"_(.15){\phantom{\cdot}}="f" &
   \ar@{.}"a";"b" \ar@{.}"c";"d"\ar@{.}"e";"f" \ar@{.}"g";"h"
             }
} $ or $
\def\objectstyle{\scriptstyle}
\def\labelstyle{\scriptstyle}
\vcenter{
  \xymatrix@-1.1pc{
 \ar@{.}[dddd]  & b \ar[l] _-{\phantom{\cdot}}="a" &  \\
    & u_n\ar[u]_-{\phantom{\cdot}}="d" &   \\
   & \ar[u] _-{\phantom{\cdot}}="c" &         \\
   & u_1\ar@{.}[u]_-{\phantom{\cdot}}="g"  & \\
    \ar[r] _-{}="e"& a \ar[u]_-{\phantom{\cdot}}="h" _-{}="f"   \ar@{-<}@/_2pc/[uuuu]_(.85){\phantom{\cdot}}="b" &
   \ar@{.}"a";"b" \ar@{.}"c";"d"\ar@{.}"e";"f" \ar@{.}"g";"h"
             }
} $ in the first case, we produce a trivial permitted thread in
$a$ by applying $V_b^{-1}$, and in the second case, by the
algorithm to calculate $\phi_A$ we get $\# \phi_A=1$.

$iii)$ If $s\neq a$, we produce a trivial permitted thread in $s$
after applying $V_b$. Also, if $s=a$ we have one of the following
cases:

$
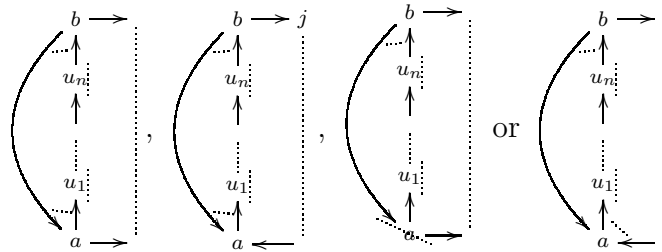
\def\objectstyle{\scriptstyle}
\def\labelstyle{\scriptstyle}
\vcenter{
  \xymatrix@-1.1pc{
& b \ar[r] \ar@{->}@/_2pc/[dddd]_(.85){}="b"_(.15){}="f" &  \ar@{.}[dddd] \\
    & u_n\ar[u]_-{\phantom{\cdot}}="d"_-{}="e" &   \\
   & \ar[u] _-{\phantom{\cdot}}="c" &         \\
   & u_1\ar@{.}[u]_-{\phantom{\cdot}}="g"  & \\
   & a \ar[u]_-{\phantom{\cdot}}="h" _-{\phantom{\cdot}}="a"   \ar[r]  &
   \ar@{.}"a";"b" \ar@{.}"c";"d"\ar@{.}"e";"f" \ar@{.}"g";"h"
             }
} $, $
\def\objectstyle{\scriptstyle}
\def\labelstyle{\scriptstyle}
\vcenter{
  \xymatrix@-1.1pc{
& b \ar[r] \ar@{->}@/_2pc/[dddd]_(.85){}="b"_(.15){}="f" & j \ar@{.}[dddd] \\
    & u_n\ar[u]_-{\phantom{\cdot}}="d"_-{}="e" &   \\
   & \ar[u] _-{\phantom{\cdot}}="c" &         \\
   & u_1\ar@{.}[u]_-{\phantom{\cdot}}="g"  & \\
   & a \ar[u] _-{\phantom{\cdot}}="h"_-{\phantom{\cdot}}="a"    &  \ar[l]
   \ar@{.}"a";"b" \ar@{.}"c";"d"\ar@{.}"e";"f"  \ar@{.}"g";"h"
             }
} $,
 $
\def\objectstyle{\scriptstyle}
\def\labelstyle{\scriptstyle}
\vcenter{
  \xymatrix@-1.1pc{
& b \ar[r] \ar@{->}@/_2pc/[dddd]_(.85){\phantom{\cdot}}="b"_(.15){}="f" &  \ar@{.}[dddd] \\
    & u_n\ar[u]_-{\phantom{\cdot}}="d"_-{}="e" &   \\
   & \ar[u] _-{\phantom{\cdot}}="c" &         \\
   & u_1\ar@{.}[u]_-{\phantom{\cdot}}="g"   & \\
   & a \ar[u]_-{\phantom{\cdot}}="h"   \ar[r]  _-{\phantom{\cdot}}="a"  &
   \ar@{.}"a";"b" \ar@{.}"c";"d"\ar@{.}"e";"f" \ar@{.}"g";"h"
             }
} $ or $
\def\objectstyle{\scriptstyle}
\def\labelstyle{\scriptstyle}
\vcenter{
  \xymatrix@-1.1pc{
& b \ar[r] \ar@{->}@/_2pc/[dddd]_(.15){}="f" &  \ar@{.}[dddd] \\
    & u_n\ar[u]_-{\phantom{\cdot}}="d"_-{}="e" &   \\
   & \ar[u] _-{\phantom{\cdot}}="c" &         \\
   & u_1\ar@{.}[u]_-{\phantom{\cdot}}="g"   & \\
   & a \ar[u]_-{\phantom{\cdot}}="h" _-{}="b"  & \ar[l]  _-{}="a"
   \ar@{.}"a";"b" \ar@{.}"c";"d"\ar@{.}"e";"f" \ar@{.}"g";"h"
             }
} $. In the first one, we obtain a trivial permitted thread in
$u_1$ applying $V_a$, in the second case, if $j=a$ we have one of
representatives of type $(3)$ and if $j\neq a$ we get a trivial
permitted thread in $b$ after applying $V_a^{-1}$. The last two
cases correspond to algebras with $\# \phi_A=1$.

Case 2.- All arrows but one, defined by vertices $b$ and $a$, are
in a cycle.

We have one of the following situations

$i)
\def\objectstyle{\scriptstyle}
\def\labelstyle{\scriptstyle}
\vcenter
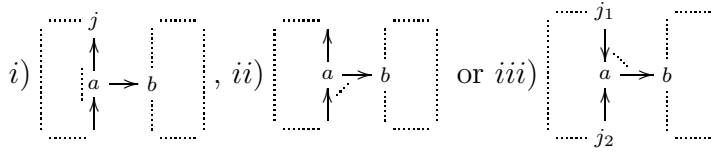
{
  \xymatrix@-1.1pc{
  \ar@{.}[r] \ar@{.}[dd] & j & \ar@{.}[r]\ar@{.}[d]     & \ar@{.}[dd] \\
                         & a \ar[u]^-{\phantom{\cdot}}="b" \ar[r] & b &   \\
  \ar@{.}[r]             & \ar[u] ^-{\phantom{\cdot}}="a" & \ar@{.}[u] \ar@{.}[r] &
   \ar@{.}"a";"b"
             }
} $, $ii)
\def\objectstyle{\scriptstyle}
\def\labelstyle{\scriptstyle}
\vcenter{
  \xymatrix@-1.1pc{
  \ar@{.}[r] \ar@{.}[dd] &  & \ar@{.}[r]\ar@{.}[d]     & \ar@{.}[dd] \\
                         & a \ar[u] \ar[r]^-{}="b" & b &   \\
  \ar@{.}[r]             & \ar[u] ^-{}="a" & \ar@{.}[u] \ar@{.}[r] &
   \ar@{.}"a";"b"
             }
} $ or $iii)
\def\objectstyle{\scriptstyle}
\def\labelstyle{\scriptstyle}
\vcenter{
  \xymatrix@-1.1pc{
  \ar@{.}[r] \ar@{.}[dd] & j_1 \ar[d]  ^-{}="a"  & \ar@{.}[r]\ar@{.}[d]     & \ar@{.}[dd] \\
                         & a \ar[r]^-{}="b" & b &   \\
  \ar@{.}[r]             & j_2 \ar[u]& \ar@{.}[u] \ar@{.}[r] &
   \ar@{.}"a";"b"
             }
} $

$i)$ If $j\neq a$ we get a trivial permitted thread in $j$ after
applying $V_a$. If $j=a$, that is, if there is a loop $\lambda$ in
$a$ there are two possibilities

$
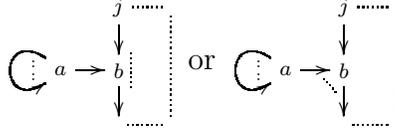
\def\objectstyle{\scriptstyle}
\def\labelstyle{\scriptstyle}
\vcenter{
  \xymatrix@-1.1pc{
   &  & j \ar@{.}[r]\ar[d]^-{\phantom{\cdot}}="c"     & \ar@{.}[dd] \\
   & a \ar@(ul,dl)^(.2){}="a"^(.8){}="b" \ar@{}[l]|{_:^:} \ar[r] & b \ar[d]^-{\phantom{\cdot}}="d" &   \\
   &  & \ar@{.}[r] &
    \ar@{.}"c";"d"
             }
} $ or $
\def\objectstyle{\scriptstyle}
\def\labelstyle{\scriptstyle}
\vcenter{
  \xymatrix@-1.1pc{
   &  & j \ar@{.}[r]\ar[d]    & \ar@{.}[dd] \\
   & a \ar@(ul,dl)^(.2){}="a"^(.8){}="b" \ar@{}[l]|{_:^:} \ar[r]^-{}="c"  & b \ar[d]^-{}="d" &   \\
   &  & \ar@{.}[r] &
    \ar@{.}"c";"d"
             }
} $ in the first one we obtain a trivial permitted thread in $j$
by applying $V_b^{-1}$ and in the second one we get a trivial
permitted thread in $b$ by applying $L_{\lambda}$.

$ii)$ We know there can not be a loop in $a$ because it is a
finite-dimensional algebra, and then we get a trivial permitted
thread in $b$ after applying $V_a$.

$iii)$ If $j_1 \neq j_2$, we produce a trivial permitted thread in
$j_1$ by applying $V_a^{-1}$. If $j_1=j_2$ there are only two
possibilities:

$
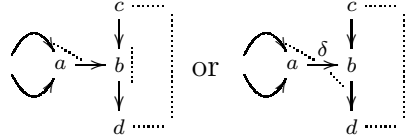
\def\objectstyle{\scriptstyle}
\def\labelstyle{\scriptstyle}
\vcenter{
  \xymatrix@-1.1pc{
 &  & c \ar@{.}[r]\ar[d]^-{\phantom{\cdot}}="c"     & \ar@{.}[dd] \\
 \ar@/_1pc/[r]\ar@/^1pc/[r]_(.7){}="a"     & a \ar[r]^-{}="b" & b \ar[d]^-{\phantom{\cdot}}="d" &   \\
        & & d \ar@{.}[r] &
   \ar@{.}"a";"b" \ar@{.}"c";"d"
             }
} $ or $
\def\objectstyle{\scriptstyle}
\def\labelstyle{\scriptstyle}
\vcenter{
  \xymatrix@-1.1pc{
 &  & c \ar@{.}[r]\ar[d]     & \ar@{.}[dd] \\
 \ar@/_1pc/[r]\ar@/^1pc/[r]_(.7){}="a"     & a \ar[r]^-{\delta}^-{}="b" & b \ar[d]^-{}="d" &   \\
        & & d \ar@{.}[r] &
   \ar@{.}"a";"b" \ar@{.}"b";"d"
             }
} $ in the first one,  we get a trivial permitted thread in $c$ by
using $V_b^{-1}$ and in the second one, after applying $F_{\delta
'}^{-1}$ we obtain a quiver in which every arrow belongs to a
cycle, so it is reduced to case 1.

Case 3.- There are exactly two arrows not belonging to a cycle.
There is one of the following options:

$i)
\def\objectstyle{\scriptstyle}
\def\labelstyle{\scriptstyle}
\vcenter
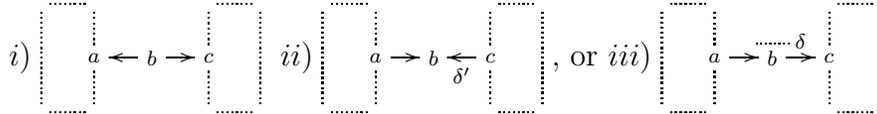
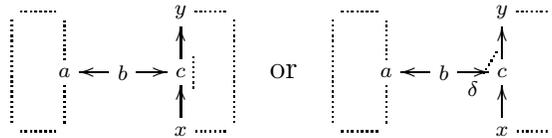
{
  \xymatrix@-1.1pc{
  \ar@{.}[r] \ar@{.}[dd] &   & & \ar@{.}[r]\ar@{.}[d]     & \ar@{.}[dd] \\
                         & a \ar@{.}[u]^-{\phantom{\cdot}}="b"  & b\ar[l]\ar[r]& c  &   \\
  \ar@{.}[r]             & \ar@{.}[u] ^-{\phantom{\cdot}}="a" & & \ar@{.}[u] \ar@{.}[r] &
             }
} $ $ii)
\def\objectstyle{\scriptstyle}
\def\labelstyle{\scriptstyle}
\vcenter{
  \xymatrix@-1.1pc{
  \ar@{.}[r] \ar@{.}[dd] &   & & \ar@{.}[r]\ar@{.}[d]     & \ar@{.}[dd] \\
                         & a \ar@{.}[u]^-{\phantom{\cdot}}="b" \ar[r] & b& c \ar[l]^-{\delta '} &   \\
  \ar@{.}[r]             & \ar@{.}[u] ^-{\phantom{\cdot}}="a" & & \ar@{.}[u] \ar@{.}[r] &
             }
} $, or $iii)
\def\objectstyle{\scriptstyle}
\def\labelstyle{\scriptstyle}
\vcenter{
  \xymatrix@-1.1pc{
  \ar@{.}[r] \ar@{.}[dd] &   & & \ar@{.}[r]\ar@{.}[d]     & \ar@{.}[dd] \\
                         & a \ar@{.}[u] \ar[r]^-{\phantom{\cdot}}="a" & b\ar[r]^-{\delta}^-{\phantom{\cdot}}="b"& c  &   \\
  \ar@{.}[r]             & \ar@{.}[u] & & \ar@{.}[u] \ar@{.}[r] &
\ar@{.}"a";"b"
             }
} $

$i)$ We can have
$$
\begin{matrix}
\def\objectstyle{\scriptstyle}
\def\labelstyle{\scriptstyle}
\vcenter{
  \xymatrix@-1.1pc{
  \ar@{.}[r] \ar@{.}[dd] &   & & y \ar@{.}[r]     & \ar@{.}[dd] \\
                         & a \ar@{.}[u]^-{\phantom{\cdot}}="b"  & b\ar[l]\ar[r]& c\ar[u]_-{\phantom{\cdot}}="c"   &   \\
  \ar@{.}[r]             & \ar@{.}[u] ^-{\phantom{\cdot}}="a" &                &x \ar[u]_-{\phantom{\cdot}}="d"  \ar@{.}[r] &
\ar@{.}"c";"d"
             }
}
&
\text{or}
&
\def\objectstyle{\scriptstyle}
\def\labelstyle{\scriptstyle}
\vcenter{
  \xymatrix@-1.1pc{
  \ar@{.}[r] \ar@{.}[dd] &   & & y \ar@{.}[r]     & \ar@{.}[dd] \\
                         & a \ar@{.}[u]^-{\phantom{\cdot}}="b"  & b\ar[l] \ar[r]_-{\delta}="e"& c\ar[u]^-{}="f"   &   \\
  \ar@{.}[r]             & \ar@{.}[u] ^-{\phantom{\cdot}}="a" &                &x \ar[u]^-{\phantom{\cdot}}="d"  \ar@{.}[r] &
\ar@{.}"e";"f"
             }
}
\end{matrix}
$$

so applying $V_c^{-1}$ we get a trivial permitted thread in $x$ or
applying $F_{\delta }$ it is reduced to case 2, respectively.

$ii)$ Is dual to $i)$.

$iii)$ Applying $F_{\delta}$ we obtain an algebra like in $i)$ or
it is reduced to case 2.

Finally, recall that trivial permitted threads correspond to
quivers with relations presenting transition vertices or branches,
and by the conditions of $A$ and Lemmas \ref{ramas2} and
\ref{ramas1} this is impossible. The only remaining options are
$\# \phi_A =1$, or $A$ derived equivalent to one of the
representatives $(1)$, $(2)$ or $(3)$ of Section
\ref{representantes}.

\endproof

Using previous results we can study and compare what happens after
applying an elementary transformation to a gentle algebra $A$ and
to an algebra obtained by removing a vertex to the quiver which
defines $A$ in a proper way:
\begin{lemma}\label{ta-x} Let $A=\mathrm{k}Q/ \left<
\mathcal{P} \right>$ be a gentle algebra, with $Q$ a quiver of two
cycles and $\#Q_0\geq 6$. Let $x\in Q_0$ as in Definition
\ref{a-x} and $T$ an elementary transformation which can be
applied to $A\setminus \{ x\}$. Then $A$ is derived equivalent
under elementary transformations to an algebra $B=\mathrm{k}Q'/
\left< \mathcal{P}' \right>$ such that the corresponding $x\in
Q_0'$ is also a vertex as in Definition \ref{a-x} and such that
$T(A\setminus \{ x\})=B\setminus\{ x \}$.
\end{lemma}

\proof

Let us analyze each one of the possible cases for vertex $x$.
Observe that any elementary transformation applied to
$A\setminus\{ x\}$ can be also be applied to $A$. Let $H_{u^{-}}$
and $H_{v^{-}}$ be the permitted threads in $A$ which involve $u$
and $v$ and which do not involve $\alpha$ and $\beta$, $H_{u^{+}}$
and $H_{v^{+}}$ be the permitted threads in $A$ which involve
$\alpha$ and $\beta$ respectively. In the corresponding
presentations of $A\setminus\{ x\}$ and $A$ in terms of their
permitted threads we have, in each case:
\begin{enumerate}
\item The presentations of $A\setminus\{ x\}$ and $A$ are

$$
\begin{matrix}
 \scriptscriptstyle
  \left[
    \begin{smallmatrix}
  &  s(H_{u^{+}}) & s(H_{u^{-}}) & s(H_{v^{-}}) & \\
  &  \vdots &\vdots    &  \vdots & \\
  \dots &  u^{+}     & u^{-}       & v^{} & \dots \\
  &   v^{+}     & \vdots   & \vdots \\
  &  \vdots & e(H_{u^{-}}) & e(H_{v^{-}})  \\
  &  e(H_{u^{+}})
    \end{smallmatrix}
  \right]
&
\text{and}
&
 \scriptscriptstyle
  \left[
    \begin{smallmatrix}
  &  s(H_{u^{+}}) & s(H_{u^{-}}) & s(C_{v^{-}})  & x^{-} & \\
  &  \vdots &\vdots    &  \vdots & &\\
  \dots &   u^{+}     & u^{-}       & v^{-} & & \dots \\
  &   x^{+}     & \vdots   & \vdots \\
  &   v^{+}     & e(H_{u^{-}}) & e(H_{v^{-}})  \\
  &  \vdots \\
  &  e(H_{u^{+}})
    \end{smallmatrix}
  \right]
\end{matrix}
$$

When the transformation $T$ does not move $u^{+}$ after $v^{-}$,
or $v^{+}$ before $u^{-}$, after applying $T$ to $A$ the trivial
permitted thread $x^{-}$ does not disappear, $x$ is also a vertex
as in the first part of Definition \ref{a-x} and we obtain the
result defining $B=T(A)$. Now, if $u^{+}$ is moved after $v^{-}$:

\[
\begin{matrix}
\scriptscriptstyle
  \left[
    \begin{smallmatrix}
   & s(H_{u^{+}}) & s(H_{u^{-}}) & s(H_{v^{-}}) & \\
   & \vdots &\vdots    &  \vdots & \\
   \dots &  u^{+}     & u^{-}       & v^{-} & \dots \\
   &  v^{+}     & \vdots   & \vdots \\
   & \vdots & e(H_{u^{-}}) & e(H_{v^{-}})  \\
   & e(H_{u^{+}})
    \end{smallmatrix}
  \right]
&
\begin{matrix} T \\ \longmapsto \end{matrix}
&  \scriptscriptstyle
  \left[
    \begin{smallmatrix}
   & s(H_{u^{+}}) & \vdots & s(H_{v^{-}}) &  \\
   & \vdots & u^{-}   &  \vdots &\\
   \dots&  v^{+}   &  \vdots       & v^{-} & \dots  \\
   &  \vdots     &    & u^{+} \\
   & e(H_{u^{+}}) &  & \vdots \\
   &        &           &  e(H_{v^{-}})
    \end{smallmatrix}
  \right]
\end{matrix}
\]

In $A$ we do the following, if $ e(H_{v^{-}})\neq v^{-}$

\[
\begin{matrix}
 \scriptscriptstyle
  \left[
   \begin{smallmatrix}
    & s(H_{u^{+}}) & s(H_{u^{-}}) & s(H_{v^{-}})  & x^{-}& \\
    &\vdots &\vdots    &  \vdots & &\\
    \dots & u^{+}     & u^{-}       & v^{-} & & \dots  \\
    & x^{+}     & \vdots   & \vdots \\
    & v^{+}     & e(H_{u^{-}}) & e(H_{v^{-}})  \\
    & \vdots \\
    & e(H_{u^{+}})
    \end{smallmatrix}
  \right]
&
\begin{matrix} F_{(v^{+},x^{+})} \\ \longmapsto \end{matrix}
&
 \scriptscriptstyle
  \left[
    \begin{smallmatrix}
    & s(H_{u^{+}}) & s(H_{u^{-}}) & s(H_{v^{-}})  & x^{-} & \\
    & \vdots &\vdots    &  \vdots & &\\
    \dots & u^{+}     & u^{-}       & v^{-} & & \dots \\
    & v^{+}     & \vdots   & x^{+} \\
    & \vdots     & e(H_{u^{-}}) & \vdots  \\
    & e(H_{u^{+}})     &           & e(H_{v^{-}})
    \end{smallmatrix}
  \right]
\end{matrix}
\]

\[
\begin{matrix}
\begin{matrix} T \\ \longmapsto \end{matrix}
&
 \scriptscriptstyle
  \left[
    \begin{smallmatrix}
  &  s(H_{u^{+}}) & \vdots & s(H_{v^{-}}) & x^{-} & \\
  &  \vdots & u^{-}   &  \vdots & &\\
  \dots &   v^{+}   &  \vdots       & v^{-} & & \dots  \\
  &   \vdots     &    & u^{+} \\
  &  e(H_{u^{+}}) &  & x^{+} \\
  &         &           &  \vdots \\
  &        &           &  e(H_{v^{-}})
    \end{smallmatrix}
  \right]
\end{matrix}
\]

else

\[
\begin{matrix}
 \scriptscriptstyle
  \left[
   \begin{smallmatrix}
    & s(H_{u^{+}}) & s(H_{u^{-}}) & s(H_{v^{-}})  & x^{-}& \\
    &\vdots &\vdots    &  \vdots & &\\
    \dots & u^{+}     & u^{-}       & v^{-} & & \dots  \\
    & x^{+}           & \vdots    \\
    & v^{+}     & e(H_{u^{-}}) \\
    & \vdots \\
    & e(H_{u^{+}})
    \end{smallmatrix}
  \right]
&
\begin{matrix} F_{(v^{+},x^{+})} \\ \longmapsto \end{matrix}
&
 \scriptscriptstyle
  \left[
    \begin{smallmatrix}
    & s(H_{u^{+}}) & s(H_{u^{-}}) & s(H_{v^{-}})  & x^{-} & \\
    & \vdots &\vdots    &  \vdots & &\\
    \dots & u^{+}     & u^{-}       & v^{-} & & \dots \\
    & v^{+}     & \vdots   & x^{+} \\
    & \vdots     & e(H_{u^{-}})  \\
    & e(H_{u^{+}})     &
    \end{smallmatrix}
  \right]
\end{matrix}
\]

\[
\begin{matrix}
\begin{matrix} T \\ \longmapsto \end{matrix}
 \scriptscriptstyle
  \left[
    \begin{smallmatrix}
  &  s(H_{u^{+}}) & \vdots & s(H_{v^{-}}) & x^{-} & \\
  &  \vdots & u^{-}   &  \vdots & &\\
  \dots &   v^{+}   &  \vdots       & v^{-} & & \dots  \\
  &   \vdots     &    & u^{+} \\
  &  e(H_{u^{+}}) &  & x^{+} \\
  &         &            \\
  &        &
    \end{smallmatrix}
  \right]
&
\begin{matrix} F^{-1}_{(x^{+},u^{+})} \\ \longmapsto \end{matrix}
&
 \scriptscriptstyle
  \left[
    \begin{smallmatrix}
  &  s(H_{u^{+}}) & \vdots & s(H_{v^{-}}) & x^{-} & \\
        &  \vdots &  x^{+}   &  \vdots & &\\
  \dots &   v^{+} &  u^{-}      & v^{-} & & \dots  \\
  &   \vdots     &  \vdots   & u^{+} \\
  &  e(H_{u^{+}}) &   \\
  &         &            \\
  &        &
    \end{smallmatrix}
  \right]
\end{matrix}
\]

In any case, the final array describes an algebra $B$ which is
obtained from the corresponding array of $T(A\setminus \{x\})$
adding a permitted thread $x^{-}$ corresponding to a vertex $x\in
Q_0 '$ as the one described in Definition \ref{a-x}. If $T$
involves the movement of $v^{+}$ before $u^{-}$, we do something
similar but using transformation $F_{(x^{+},u^{+})}^{-1}$.

\item Presentations of $A\setminus\{ x\}$ and $A$ are

$$
\begin{matrix}
 \scriptscriptstyle
  \left[
    \begin{smallmatrix}
    & s(H_{u^{+}}) & s(H_{u^{-}}) & w^{-} & \\
    & \vdots &\vdots    &  v^{-} & \\
    \dots & u^{+}      & u^{-}       & \vdots & \dots  \\
    & w^{+}      & \vdots   & e(H_{v^{-}}) \\
    & z^{+}       & e(H_{u^{-}}) \\
    &\vdots \\
    & e(H_{u^{+}})
    \end{smallmatrix}
  \right]
&
\text{and}
&
 \scriptscriptstyle
  \left[
   \begin{smallmatrix}
    & s(H_{u^{+}}) & s(H_{u^{-}}) & x^{+}  & x^{-} & \\
    & \vdots &\vdots    &  w^{-} & & \\
    \dots & u^{+}      & u^{-}       & v^{-} & &\dots  \\
    & w^{+}      & \vdots   & \vdots \\
    & z^{+}       & e(H_{u^{-}})& e(H_{v^{-}}) \\
    & \vdots \\
    & e(H_{u^{+}})
    \end{smallmatrix}
  \right].
\end{matrix}
$$

The only way which could make the permitted thread in $x$
disappear would be the movement of $w^{-}$ before $x^{-}$, but
this does not happens because $w$ is the beginning of a permitted
thread in $A\setminus \{x \}$. In this case $B:=T(A)$ solves the
problem.
\end{enumerate}

\endproof

Applying previous lemma several times we have:

\begin{corollary}\label{cta-x}
Let $A=\mathrm{k}Q/ \left< \mathcal{P} \right>$ be a gentle
algebra, with $Q$ a quiver with two cycles and $\# Q_0\geq 6$. Let
$x\in Q_0$ as in Definition \ref{a-x} and $T$ a composition of
elementary transformations such that $T(A\setminus\{x\})$ makes
sense. Then $A$ is derived equivalent under elementary
transformations to an algebra $B=\mathrm{k}Q'/ \left< \mathcal{P}'
\right>$ such that the corresponding $x\in Q_0'$ is also as in
Definition \ref{a-x} and $T(A\setminus \{ x\})=B\setminus\{ x \}$.
\end{corollary}

\subsection{Adding a vertex to the representatives}

We analyze now the different ways in which it is possible to add a
vertex as the ones described in Definition \ref{a-x} to the
representatives of Section \ref{representantes} and we prove that
in any case there exist some elementary transformations which turn
these algebras into one of the normal forms. First let us see what
happens if the added vertex is a transition vertex. We need the
following:

\begin{remark}\label{transicion}{\normalfont

\[
\def\objectstyle{\scriptstyle}
\def\labelstyle{\scriptstyle}
\vcenter{
  \xymatrix@-1.1pc{
 \ar@{-}@/^1pc/[dd] & & \ar@{-}@/_1pc/[rr] & & & \ar@{-}@/_1pc/[dd]\\
    & u  \ar[r] _-{\phantom{\cdot}}="a" & w \ar[r]_-{\delta}="b" & z \ar[r]_-{\phantom{\cdot}}="c" & v \\
   & & & & &
   \ar@{.}"b";"c"
             }
}
\]
is derived equivalent to
\[
\def\objectstyle{\scriptstyle}
\def\labelstyle{\scriptstyle}
\vcenter{
  \xymatrix@-1.1pc{
 \ar@{-}@/^1pc/[dd] & \ar@{-}@/_1pc/[rr] &  & & & \ar@{-}@/_1pc/[dd]\\
    & u  \ar[r] ^-{}_-{\phantom{\cdot}}="a" & w \ar[r]^-{}_-{\phantom{\cdot}}="b" & z \ar[r]_-{\phantom{\cdot}}="c" & v \\
   & & & & &
   \ar@{.}"a";"b"
             }
}
\]
under transformation $F_{\delta}$. }\end{remark}

Under certain conditions, adding a vertex which defines a trivial
permitted thread to a branch and removing it from another produces
derived equivalent algebras.

\begin{lemma}\label{lema1}
The gentle algebra defined by the next quiver with relations
\[
\def\objectstyle{\scriptstyle}
\def\labelstyle{\scriptstyle}
\vcenter{
  \xymatrix@-1.1pc{
 & w_1\ar[r] & w_2\ar[r] & \ar@{.}[r] & \ar[r] & w_r \ar[ddr] & & & z_m \ar[ddl] & \ar[l] & \ar@{.}[l] & z_1 \ar[l] & x \ar[l] & \\
&  &  &  & \ar@{-}@/^1pc/[dd]  &  & &  & & \ar@{-}@/_1pc/[dd]\\
    & & & & &  u  \ar[r] _-{\phantom{\cdot}}="a" & w \ar[r]_-{\phantom{\cdot}}="b" & z \ar[r]_-{\phantom{\cdot}}="c" & v \\
   & & & & & & & & &
   \ar@{.}"a";"b"\ar@{.}"b";"c"
             }
}
\]
is derived equivalent under elementary transformations to a gentle
algebra defined by
\[
\def\objectstyle{\scriptstyle}
\def\labelstyle{\scriptstyle}
\vcenter{
  \xymatrix@-1.1pc{
  x \ar[r] & w_1\ar[r] & w_2\ar[r] & \ar@{.}[r] & \ar[r] & w_r \ar[ddr] & & & z_m \ar[ddl] & \ar[l] & \ar@{.}[l] & z_1 \ar[l] &  & \\
&  &  &  & \ar@{-}@/^1pc/[dd]  &  & &  & & \ar@{-}@/_1pc/[dd]\\
    & & & & &  u  \ar[r] _-{\phantom{\cdot}}="a" & w \ar[r]_-{\phantom{\cdot}}="b" & z \ar[r]_-{\phantom{\cdot}}="c" & v \\
   & & & & & & & & &
   \ar@{.}"a";"b"\ar@{.}"b";"c"
             }
}
\]
\end{lemma}

\proof

First we apply $V_z^{-1}$ and after renaming the vertices we get
\[
\def\objectstyle{\scriptstyle}
\def\labelstyle{\scriptstyle}
\vcenter{
  \xymatrix@-1.1pc{
  x \ar[r] & w_1\ar[r] & w_2\ar[r] & \ar@{.}[r] & \ar[r] & w_r & & & z_m \ar[ddl] & \ar[l] & \ar@{.}[l] & z_1 \ar[l] &  & \\
&  &  &  & \ar@{-}@/^1pc/[dd]  &  & &  & & \ar@{-}@/_1pc/[dd]\\
    & & & & &  u  \ar[r]_-{\phantom{\cdot}}="a" & w  \ar[uul]_-{\delta^{'}} \ar[r]_-{\phantom{\cdot}}="b" & z \ar[r]_-{\phantom{\cdot}}="c" & v \\
   & & & & & & & & &
   \ar@{.}"a";"b"\ar@{.}"b";"c"
             }
}
\]
using $F_{\delta^{'}}^{-1}$ we obtain
\[
\def\objectstyle{\scriptstyle}
\def\labelstyle{\scriptstyle}
\vcenter{
  \xymatrix@-1.1pc{
  x \ar[r] & w_1\ar[r] & w_2\ar[r] & \ar@{.}[r] & \ar[r]_-{}="d" & w_r \ar[ddr]_-{}="e"_-{\delta}& & & z_m \ar[ddl] & \ar[l] & \ar@{.}[l] & z_1 \ar[l] &  & \\
&  &  &  & \ar@{-}@/^1pc/[dd]  &  & &  & & \ar@{-}@/_1pc/[dd]\\
    & & & & &  u  \ar[r] _-{\phantom{\cdot}}="a" & w   \ar[r]_-{\phantom{\cdot}}="b" & z \ar[r]_-{\phantom{\cdot}}="c" & v \\
   & & & & & & & & &
   \ar@{.}"a";"b"\ar@{.}"b";"c"\ar@{.}"d";"e"
             }
}
\]
and by Remark \ref{transicion} it is derived equivalent to
\[
\def\objectstyle{\scriptstyle}
\def\labelstyle{\scriptstyle}
\vcenter{
  \xymatrix@-1.1pc{
  x \ar[r]^-{\delta} _-{\phantom{\cdot}}="d"& w_1\ar[r]_-{\phantom{\cdot}}="e" & w_2\ar[r] & \ar@{.}[r] & \ar[r] & w_r \ar[ddr]& & & z_m \ar[ddl] & \ar[l] & \ar@{.}[l] & z_1 \ar[l] &  & \\
&  &  &  & \ar@{-}@/^1pc/[dd]  &  & &  & & \ar@{-}@/_1pc/[dd]\\
    & & & & &  u  \ar[r] _-{\phantom{\cdot}}="a" & w   \ar[r]_-{\phantom{\cdot}}="b" & z \ar[r]_-{\phantom{\cdot}}="c" & v \\
   & & & & & & & & &
   \ar@{.}"a";"b"\ar@{.}"b";"c"\ar@{.}"d";"e"
             }
}
\]
and using $F_{\delta}$ we get the result.

\endproof

\begin{remark}\label{movrama}{\normalfont
As a consequence
\[
\def\objectstyle{\scriptstyle}
\def\labelstyle{\scriptstyle}
\vcenter{
  \xymatrix@-1.1pc{
 \ar@{-}@/^1pc/[dd] & & & x \ar[d] & & \ar@{-}@/_1pc/[dd]\\
    & u  \ar[r] _-{\phantom{\cdot}}="a" & w \ar[r]_-{\phantom{\cdot}}="b" & z \ar[r]_-{\phantom{\cdot}}="c" & v \\
   & & & & &
   \ar@{.}"a";"b"\ar@{.}"b";"c"
             }
}
\]
is derived equivalent to
\[
\def\objectstyle{\scriptstyle}
\def\labelstyle{\scriptstyle}
\vcenter{
  \xymatrix@-1.1pc{
 \ar@{-}@/^1pc/[dd] & & x \ar[d] & & & \ar@{-}@/_1pc/[dd]\\
    & u  \ar[r] ^-{}_-{\phantom{\cdot}}="a" & w \ar[r]^-{}_-{\phantom{\cdot}}="b" & z \ar[r]_-{\phantom{\cdot}}="c" & v \\
   & & & & &
   \ar@{.}"a";"b"\ar@{.}"b";"c"
             }
}
\]
}\end{remark}

We work now with the normal forms and we analyze all different
ways in which we can add a vertex as described in Definition
\ref{a-x}. If it is a transition vertex, in most of the cases
using Remark \ref{transicion} and the symmetry of the normal
forms, we can transform the algebra into one of the
representatives by using elementary transformations, like in the
case of representative $(7)$. However, the conditions over the
parameters $a$, $b$, $k$ and $q$ could change, we need then the
following results:

\begin{proposition}
The algebras associated to the next quiver with relations are
derived equivalent.

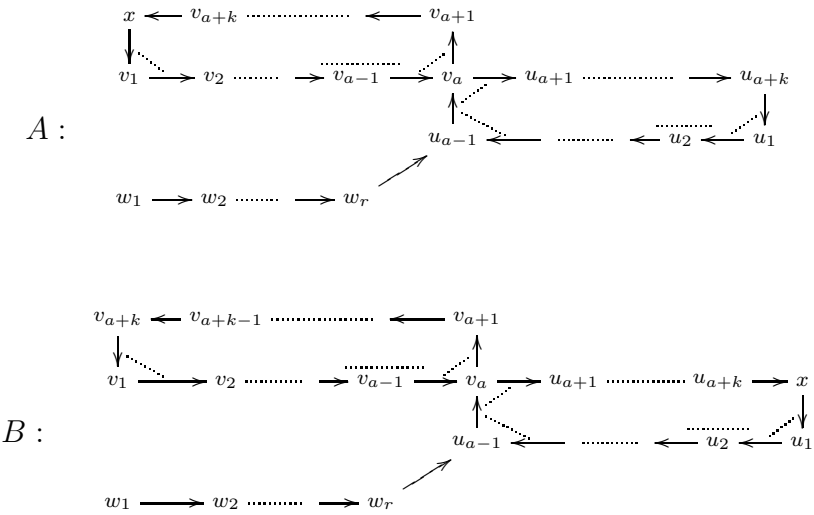
\[A:\quad
\def\objectstyle{\scriptstyle}
\def\labelstyle{\scriptstyle}
\vcenter{
  \xymatrix@-1.1pc{
      x \ar[d]^-{}="a"   & v_{a+k} \ar[l]  & & \ar@{.}[ll]  & v_{a+1} \ar[l]   \\
      v_1 \ar[r]^-{}="b" & v_2 \ar@{.}[r] &\ar[r]^-{\phantom{\cdot}}="c" &    v_{a-1} \ar[r]^-{\phantom{\cdot}}="d"^-{}="e"  & v_a \ar[u]_-{}="f" \ar[r]^-{}="i" & u_{a+1} \ar @{.} [rr] &              & \ar[r]     & u_{a+k} \ar[d] ^-{}="j" \\
                 &                  &        &                   & u_{a-1} \ar[u] ^-{}="h"& \ar[l]  ^-{}="g" & \ar @{.} [l] & u_2 \ar[l]_-{\phantom{\cdot}}="m"  & u_1 \ar[l] ^-{}="k"_-{\phantom{\cdot}}="l" \\
      w_1 \ar[r] & w_2 \ar@{.}[r] & \ar[r] & w_r \ar[ur] \\
      \ar@{.}"a";"b" \ar@{.} "c";"d" \ar@{.}"e";"f" \ar@{.}"g";"h" \ar@{.}"h";"i" \ar@{.}"j";"k" \ar@{.}"l";"m"
                  }
}
\]

\[B:\quad
\def\objectstyle{\scriptstyle}
\def\labelstyle{\scriptstyle}
\vcenter{
  \xymatrix@-1.1pc{
      v_{a+k} \ar[d]^-{}="a"   & v_{a+k-1} \ar[l]  & & \ar@{.}[ll]  & v_{a+1} \ar[l]   \\
      v_1 \ar[r]^-{}="b" & v_2 \ar@{.}[r] &\ar[r]^-{\phantom{\cdot}}="c" &    v_{a-1} \ar[r]^-{\phantom{\cdot}}="d"^-{}="e"  & v_a \ar[u]_-{}="f" \ar[r]^-{}="i" & u_{a+1} \ar @{.} [rr] &              & u_{a+k} \ar[r]     & x \ar[d] ^-{}="j" \\
                 &                  &        &                   & u_{a-1} \ar[u] ^-{}="h"& \ar[l]  ^-{}="g" & \ar @{.} [l] & u_2 \ar[l]_-{\phantom{\cdot}}="m"  & u_1 \ar[l] ^-{}="k"_-{\phantom{\cdot}}="l" \\
      w_1 \ar[r] & w_2 \ar@{.}[r] & \ar[r] & w_r \ar[ur] \\
      \ar@{.}"a";"b" \ar@{.} "c";"d" \ar@{.}"e";"f" \ar@{.}"g";"h" \ar@{.}"h";"i" \ar@{.}"j";"k" \ar@{.}"l";"m"
                  }
}
\]
\end{proposition}

\begin{remark}{\normalfont
$\phi_A =\phi_B=[(k,a+k),(k+1,a+k+1),(2a-2+r,r)]$. }\end{remark}

\proof

Applying $V_{v_a}^{-1}$ to $A$ we get:
%
%
%
%
%
%

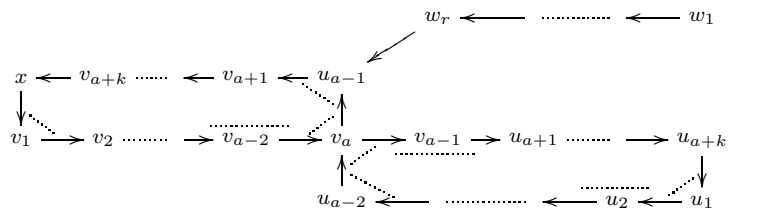
\[
\def\objectstyle{\scriptstyle}
\def\labelstyle{\scriptstyle}
\vcenter{
  \xymatrix@-1.1pc{
      & & & & &  w_r \ar[dl] & \ar[l] & \ar@{.}[l] & w_1 \ar[l] \\
      x \ar[d]^-{}="x"   & v_{a+k} \ar[l]  & \ar@{.}[l]&  v_{a+1}\ar[l] & u_{a-1} \ar[l]^-{}="b"   \\
      v_1 \ar[r]^-{}="y" & v_2 \ar@{.}[r] &\ar[r]^-{\phantom{\cdot}}="c" &    v_{a-2} \ar[r]^-{\phantom{\cdot}}="d"^-{}="e"  & v_a \ar[u]_-{}="f" \ar[r]^-{}="i"_-{\phantom{\cdot}}="z" & v_{a-1} \ar[r]_-{\phantom{\cdot}}="w" & u_{a+1} \ar@{.}[r]         & \ar[r]     & u_{a+k} \ar[d] ^-{}="j" \\
                 &                  &        &                   & u_{a-2} \ar[u] ^-{}="h"& \ar[l]  ^-{}="g" & \ar @{.} [l] & u_2 \ar[l]_-{\phantom{\cdot}}="m"  & u_1 \ar[l] ^-{}="k"_-{\phantom{\cdot}}="l" &.\\
      \ar@{.}"f";"b" \ar@{.} "c";"d" \ar@{.}"e";"f" \ar@{.}"g";"h" \ar@{.}"h";"i" \ar@{.}"j";"k" \ar@{.}"l";"m" \ar@{.}"x";"y" \ar@{.}"z";"w"
                  }
}
\]

Using Lemma \ref{lema1} and Remarks \ref{transicion} and
\ref{movrama} it is derived equivalent to

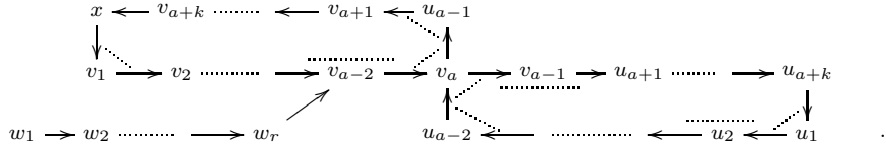
\[
\def\objectstyle{\scriptstyle}
\def\labelstyle{\scriptstyle}
\vcenter{
  \xymatrix@-1.1pc{
     &  x \ar[d]^-{}="x"   & v_{a+k} \ar[l]  & \ar@{.}[l] &v_{a+1}\ar[l]  & u_{a-1} \ar[l]^-{}="b"   \\
     & v_1 \ar[r]^-{}="y" & v_2 \ar@{.}[r] &\ar[r]^-{\phantom{\cdot}}="c" &    v_{a-2} \ar[r]^-{\phantom{\cdot}}="d"^-{}="e"  & v_a \ar[u]_-{}="f" \ar[r]^-{}="i"_-{\phantom{\cdot}}="z" & v_{a-1} \ar[r]_-{\phantom{\cdot}}="w" & u_{a+1} \ar@{.}[r]         & \ar[r]     & u_{a+k} \ar[d] ^-{}="j" \\
      w_1 \ar[r] & w_2 \ar@{.}[r] & \ar[r] & w_r \ar[ur] &  & u_{a-2} \ar[u] ^-{}="h"& \ar[l]  ^-{}="g" & \ar @{.} [l] & u_2 \ar[l]_-{\phantom{\cdot}}="m"  & u_1 \ar[l] ^-{}="k"_-{\phantom{\cdot}}="l" & .\\
      \ar@{.}"f";"b" \ar@{.} "c";"d" \ar@{.}"e";"f" \ar@{.}"g";"h" \ar@{.}"h";"i" \ar@{.}"j";"k" \ar@{.}"l";"m" \ar@{.}"x";"y" \ar@{.}"z";"w"
                  }
}
\]
Renaming the vertices, it is the quiver with relations associated
to $B$.

\endproof

\begin{proposition}
The algebras associated to the following quivers with relations
are derived equivalent.

$\begin{matrix} A:\;
\def\objectstyle{\scriptstyle}
\def\labelstyle{\scriptstyle}
\vcenter
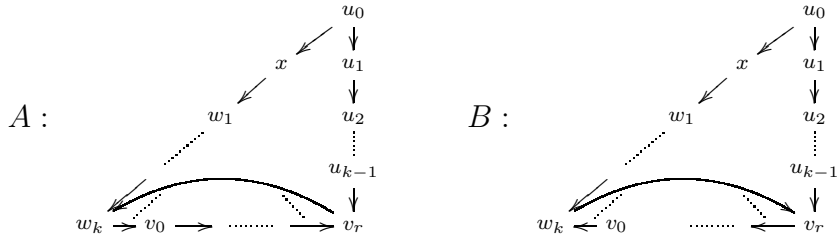
{
  \xymatrix@-1.3pc{
    & &     &     & u_0 \ar[dl] \ar[d]\\
    & &     & x \ar[dl] & u_1 \ar[d] \\
    & & w_1 \ar@{.}[dl] &     & u_2 \ar@{.}[d] \\
    &\ar[dl] &     &     & u_{k-1} \ar[d] \\
 w_{k}\ar[r]^-{}="a" & v_0 \ar[r] & \ar@{.}[r] & \ar[r] ^-{}="c" & v_r \ar@/ _1.5pc/[llll]^(.3){}="d" ^(.7){}="b"
      \ar@{.}"a";"b" \ar@{.}"c";"d"
                  }
} & & & B:\;
\def\objectstyle{\scriptstyle}
\def\labelstyle{\scriptstyle}
\vcenter{
  \xymatrix@-1.3pc{
    & &     &     & u_0 \ar[dl] \ar[d]\\
    & &     & x \ar[dl] & u_1 \ar[d] \\
    & & w_1 \ar@{.}[dl] &     & u_2 \ar@{.}[d] \\
    &\ar[dl] &     &     & u_{k-1} \ar[d] \\
 w_{k}\ar@/ ^1.5pc/[rrrr]^(.7){}="d" ^(.3){}="b" & v_0 \ar[l]^-{}="a" & \ar@{.}[r] &  & v_r \ar[l] ^-{}="c"
      \ar@{.}"a";"b" \ar@{.}"c";"d"
                  }
}
\end{matrix}$
\end{proposition}

\begin{remark}{\normalfont
$\phi_A =\phi_B =[(k,k),(k+1,k+1),(r,r+2)]$. }\end{remark}

\proof

Using Remark \ref{transicion}, $A$ is derived equivalent to
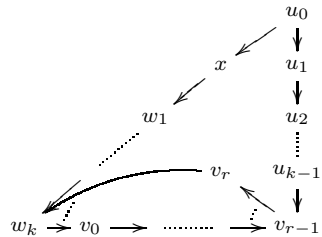
\[
\def\objectstyle{\scriptstyle}
\def\labelstyle{\scriptstyle}
\vcenter{
  \xymatrix@-1.3pc{
    & &     &     & u_0 \ar[dl] \ar[d]\\
    & &     & x \ar[dl] & u_1 \ar[d] \\
    & & w_1 \ar@{.}[dl] &     & u_2 \ar@{.}[d] \\
    &\ar[dl] &     &  v_r  \ar@/ _.7pc/[dlll]^(.7){}="b"    & u_{k-1} \ar[d] \\
 w_{k}\ar[r]^-{}="a" & v_0 \ar[r] & \ar@{.}[r] & \ar[r] ^-{}="c" & v_{r-1} \ar[ul]^{}="d"
      \ar@{.}"a";"b" \ar@{.}"c";"d"
                  }
}
\]

Applying the result $r$ times and renaming vertices we obtain $B$.
\endproof

\begin{proposition}
The algebras associated to the next quiver with relations are
derived equivalent.

\[A:\quad
\def\objectstyle{\scriptstyle}
\def\labelstyle{\scriptstyle}
\vcenter{
  \xymatrix@-1.1pc{
     x \ar[r] & v_1 \ar@{.}[r] & v_{q-1} \ar[r]_{\phantom{\cdot}}="a" & v_q \ar[r]_{\phantom{\cdot}}="b"^(.7){}="c" & u_{k}  \ar@/_1pc/[llll]^(.3){}="d" &
\ar[l] & u_1 \ar@{.}[l] & u_0 \ar[l] \ar@/^1pc/[rrrr]^(.3){}="f" &
w_q \ar[l]^(.7){}="e"^-{\phantom{\cdot}}="h" & w_{q-1}
\ar[l]^-{\phantom{\cdot}}="g" & w_2 \ar@{.}[l] & w_1 \ar[l]
       \ar@{.} "c";"d"\ar@{.}"e";"f"
                  }
}
\]

\[B:\quad
\def\objectstyle{\scriptstyle}
\def\labelstyle{\scriptstyle}
\vcenter{ \xymatrix@-1.1pc{
     x \ar[r] & v_1 \ar@{.}[r] & v_{q-1} \ar[r]_{\phantom{\cdot}}="a" & v_q \ar[r]_{\phantom{\cdot}}="b"^(.7){}="c" & u_{k}  \ar@/_1pc/[llll]^(.3){}="d"\ar[r] &  & u_1 \ar@{.}[l] \ar[r] & u_0 \ar@/^1pc/[rrrr]^(.3){}="f" & w_q \ar[l]^(.7){}="e"^-{\phantom{\cdot}}="h" & w_{q-1} \ar[l]^-{\phantom{\cdot}}="g" & w_2 \ar@{.}[l] & w_1 \ar[l]
       \ar@{.} "c";"d"\ar@{.}"e";"f"
                  }
}
\]
\end{proposition}

\begin{remark}{\normalfont
$\phi_A = \phi_B =[(k,k),(q+1,q+2),(q,q+1)]$. }\end{remark}

\proof

Using Remark \ref{metorama} several times, $A$ is derived
equivalent to
\[
\def\objectstyle{\scriptstyle}
\def\labelstyle{\scriptstyle}
\vcenter{
  \xymatrix@-1.1pc{
     x \ar[r] & v_1 \ar@{.}[r] & v_{q-1} \ar[r]_{\phantom{\cdot}}="a" & v_q \ar[r]_{\phantom{\cdot}}="b"^(.7){}="c" & u_{k}  \ar@/_1pc/[l]^(.3){}="d" &
\ar[l] & u_1 \ar@{.}[l] & u_0 \ar[l] \ar@/^1pc/[r]^(.3){}="f" &
w_q \ar[l]^(.7){}="e"^-{\phantom{\cdot}}="h" & w_{q-1}
\ar[l]^-{\phantom{\cdot}}="g" & w_2 \ar@{.}[l] & w_1 \ar[l]
      \ar@{.}"a";"b" \ar@{.} "c";"d"\ar@{.}"e";"f"\ar@{.}"g";"h"
                  }
}
\]

%

After applying $V_{u_0}$ and $V_{w_q}$ we have:

\[
\def\objectstyle{\scriptstyle}
\def\labelstyle{\scriptstyle}
\vcenter{
  \xymatrix@-1.1pc{
     x \ar[r] & v_1 \ar@{.}[r] & v_{q-1} \ar[r]_{\phantom{\cdot}}="a" & v_q \ar[r]_{\phantom{\cdot}}="b"^(.7){}="c" & u_{k}  \ar@/_1pc/[l]^(.3){}="d" &
\ar[l] & \ar@{.}[l] & u_2 \ar[l]  & u_1 \ar[l] \ar[r]& w_q
\ar@/^1pc/[r]^(.3){}="f" & u_0
\ar[l]^(.7){}="e"^-{\phantom{\cdot}}="h" & w_{q-1}
\ar[l]^-{\phantom{\cdot}}="g" & w_2 \ar@{.}[l] & w_1 \ar[l]
      \ar@{.}"a";"b" \ar@{.} "c";"d"\ar@{.}"e";"f"\ar@{.}"g";"h"
                  }
}
\]

and applying $V_{u_1}$ it is derived equivalent to
\[
\def\objectstyle{\scriptstyle}
\def\labelstyle{\scriptstyle}
\vcenter{
  \xymatrix@-1.1pc{
     x \ar[r] & v_1 \ar@{.}[r] & v_{q-1} \ar[r]_{\phantom{\cdot}}="a" & v_q \ar[r]_{\phantom{\cdot}}="b"^(.7){}="c" & u_{k}  \ar@/_1pc/[l]^(.3){}="d" &
\ar[l] & \ar@{.}[l] & u_2 \ar[l] \ar[r] & u_1  & w_q \ar[l]
\ar@/^1pc/[r]^(.3){}="f" & u_0
\ar[l]^(.7){}="e"^-{\phantom{\cdot}}="h" & w_{q-1}
\ar[l]^-{\phantom{\cdot}}="g" & w_2 \ar@{.}[l] & w_1 \ar[l]
      \ar@{.}"a";"b" \ar@{.} "c";"d"\ar@{.}"e";"f"\ar@{.}"g";"h"
                  }
}
\]

after renaming the vertices we get

\[
\def\objectstyle{\scriptstyle}
\def\labelstyle{\scriptstyle}
\vcenter{
  \xymatrix@-1.1pc{
     x \ar[r] & v_1 \ar@{.}[r] & v_{q-1} \ar[r]_{\phantom{\cdot}}="a" & v_q \ar[r]_{\phantom{\cdot}}="b"^(.7){}="c" & u_{k}  \ar@/_1pc/[l]^(.3){}="d" &
\ar[l] & \ar@{.}[l] & u_2 \ar[l] \ar[r] & u_1 & u_0
\ar[l]\ar@/^1pc/[r]^(.3){}="f" & w_q
\ar[l]^(.7){}="e"^-{\phantom{\cdot}}="h" & w_{q-1}
\ar[l]^-{\phantom{\cdot}}="g" & w_2 \ar@{.}[l] & w_1 \ar[l]
      \ar@{.}"a";"b" \ar@{.} "c";"d"\ar@{.}"e";"f"\ar@{.}"g";"h"
                  }
}
\]

Applying $V_{u_0}$ and $V_{w_q}$ we get

\[
\def\objectstyle{\scriptstyle}
\def\labelstyle{\scriptstyle}
\vcenter{
  \xymatrix@-1.1pc{
     x \ar[r] & v_1 \ar@{.}[r] & v_{b-1} \ar[r]_{\phantom{\cdot}}="a" & v_b \ar[r]_{\phantom{\cdot}}="b"^(.7){}="c" & u_{k}  \ar@/_1pc/[l]^(.3){}="d" &
\ar[l] & \ar@{.}[l] & u_2 \ar[l] \ar[r] & u_1 \ar[r]& u_0
\ar@/^1pc/[r]^(.3){}="f" & w_q
\ar[l]^(.7){}="e"^-{\phantom{\cdot}}="h" & w_{q-1}
\ar[l]^-{\phantom{\cdot}}="g" & w_2 \ar@{.}[l] & w_1 \ar[l]
      \ar@{.}"a";"b" \ar@{.} "c";"d"\ar@{.}"e";"f"\ar@{.}"g";"h"
                  }
}
\]

After $k$ steps we obtain

\[
\def\objectstyle{\scriptstyle}
\def\labelstyle{\scriptstyle}
\vcenter{
  \xymatrix@-1.1pc{
     x \ar[r] & v_1 \ar@{.}[r] & v_{q-1} \ar[r]_{\phantom{\cdot}}="a" & v_q \ar[r]_{\phantom{\cdot}}="b"^(.7){}="c" &u_{k}\ar[r]\ar@/_1pc/[l]^(.3){}="d" &
 & u_1 \ar@{.}[l]\ar[r] & u_0  \ar@/^1pc/[r]^(.3){}="f" & w_q \ar[l]^(.7){}="e"^-{\phantom{\cdot}}="h" & w_{q-1} \ar[l]^-{\phantom{\cdot}}="g" & w_2 \ar@{.}[l] & w_1 \ar[l]
      \ar@{.}"a";"b" \ar@{.} "c";"d"\ar@{.}"e";"f"\ar@{.}"g";"h"
                  }
}
\]

which is derived equivalent to $B$ by Remark \ref{metorama}.

\endproof

If the added vertex is as in Definition \ref{a-x} $(2)$, and there
is in fact an arrow $\kappa$ as the one described in the mentioned
definition, we can reduce the problem to the case of a transition
vertex using Remark \ref{metorama} which has been just analyzed.
If this is not the case, by Remark \ref{movrama} and the type of
normal forms which we are working with, it is enough to prove the
next result:

\begin{proposition}
The algebras associated to the following quivers with relations
are derived equivalent.

\[A:\quad
\def\objectstyle{\scriptstyle}
\def\labelstyle{\scriptstyle}
\vcenter{
  \xymatrix@-1.1pc{
      v_{b+q} \ar[d]^-{}="a"   & v_{b+q-1} \ar[l]  & & \ar@{.}[ll]  & v_{b+1} \ar[l]   \\
      v_1 \ar[r]^-{}="b" & v_2 \ar@{.}[r] &\ar[r]^-{\phantom{\cdot}}="c" &    v_{b-1} \ar[r]^-{\phantom{\cdot}}="d"^-{}="e"  & v_b \ar[u]_-{}="f" \ar[r]^-{}="i" & u_{a+1} \ar @{.} [rr] &              & \ar[r]     & u_{a+k} \ar[d] ^-{}="j" \\
                 &                  &        &     x\ar[u]     & u_{a-1} \ar[u] ^-{}="h"& \ar[l]  ^-{}="g" & \ar @{.} [l] & u_2 \ar[l]_-{\phantom{\cdot}}="m"  & u_1 \ar[l] ^-{}="k"_-{\phantom{\cdot}}="l" \\
      w_1 \ar[r] & w_2 \ar@{.}[r] & \ar[r] & w_r \ar[ur] \\
      \ar@{.}"a";"b" \ar@{.} "c";"d" \ar@{.}"e";"f" \ar@{.}"g";"h" \ar@{.}"h";"i" \ar@{.}"j";"k" \ar@{.}"l";"m"
                  }
}
\]

\[B:\quad
\def\objectstyle{\scriptstyle}
\def\labelstyle{\scriptstyle}
\vcenter{
  \xymatrix@-1.1pc{
      v_{b+q} \ar[d]^-{}="a"   & v_{b+q-1} \ar[l]  & & \ar@{.}[ll]  & v_{b+1} \ar[l]   \\
      v_1 \ar[r]^-{}="b" & v_2 \ar@{.}[r] &\ar[r]^-{\phantom{\cdot}}="c" &    v_{b-1} \ar[r]^-{\phantom{\cdot}}="d"^-{}="e"  & v_b \ar[u]_-{}="f" \ar[r]^-{}="i" & u_{a+1} \ar @{.} [rr] &              &  \ar[r]     & u_{a+k} \ar[d] ^-{}="j" \\
                 &                  &        &                   & u_{a-1} \ar[u] ^-{}="h"& \ar[l]  ^-{}="g" & \ar @{.} [l] & u_2 \ar[l]_-{\phantom{\cdot}}="m"  & u_1 \ar[l] ^-{}="k"_-{\phantom{\cdot}}="l" \\
      x \ar[r] & w_1 \ar@{.}[r] & \ar[r] & w_r \ar[ur] \\
      \ar@{.}"a";"b" \ar@{.} "c";"d" \ar@{.}"e";"f" \ar@{.}"g";"h" \ar@{.}"h";"i" \ar@{.}"j";"k" \ar@{.}"l";"m"
                  }
}
\]
\end{proposition}

\begin{remark}{\normalfont
$\phi_A =\phi_B=[(k,a+k),(q,b+q),(a+b-2+r+1,r)]$. }\end{remark}

\proof

Let $\delta$ be the only arrow of $A$ such that $s(\delta)=x$ and
$e(\delta)=v_{b-1}$. After applying $F_{\delta}$ and $V_{v_{b-1}}$
we get
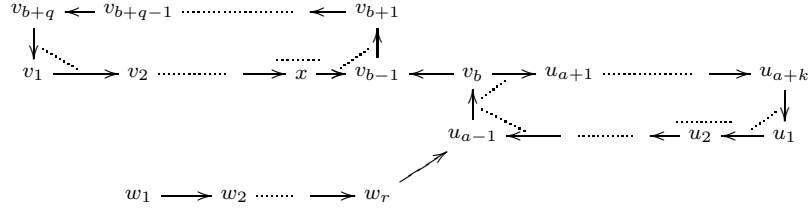
\[
\def\objectstyle{\scriptstyle}
\def\labelstyle{\scriptstyle}
\vcenter{
  \xymatrix@-1.1pc{
      v_{b+q} \ar[d]^-{}="a"   & v_{b+q-1} \ar[l]  &  & \ar@{.}[ll]  & v_{b+1} \ar[l]   \\
      v_1 \ar[r]^-{}="b" & v_2 \ar@{.}[r] &\ar[r]^-{\phantom{\cdot}}="z" & x \ar[r]^-{}="c" \ar[r]^-{\phantom{\cdot}}="y" &    v_{b-1} \ar[u]^-{}="d"  & v_b \ar[l] \ar[r]^-{}="i" & u_{a+1} \ar @{.} [rr] &              & \ar[r]     & u_{a+k} \ar[d] ^-{}="j" \\
                 &                  &        &     &    & u_{a-1} \ar[u] ^-{}="h"& \ar[l]  ^-{}="g" & \ar @{.} [l] & u_2 \ar[l]_-{\phantom{\cdot}}="m"  & u_1 \ar[l] ^-{}="k"_-{\phantom{\cdot}}="l" \\
     &  w_1 \ar[r] & w_2 \ar@{.}[r] & \ar[r] & w_r \ar[ur] \\
      \ar@{.}"a";"b" \ar@{.} "c";"d"  \ar@{.}"g";"h" \ar@{.}"h";"i" \ar@{.}"j";"k" \ar@{.}"l";"m"\ar@{.}"y";"z"
                  }
}
\]

and by Lemma \ref{lema1} and Remark \ref{movrama} it is derived
equivalent to

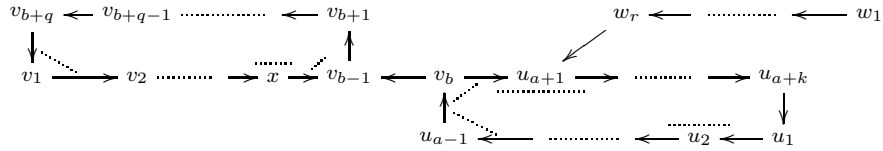
\[
\def\objectstyle{\scriptstyle}
\def\labelstyle{\scriptstyle}
\vcenter{
  \xymatrix@-1.1pc{
      v_{b+q} \ar[d]^-{}="a"   & v_{b+q-1} \ar[l]  &  & \ar@{.}[ll]  & v_{b+1} \ar[l] & & &  w_r \ar[dl] &  \ar[l] & \ar@{.}[l] & w_1 \ar[l] \\
      v_1 \ar[r]^-{}="b" & v_2 \ar@{.}[r] &\ar[r]^-{\phantom{\cdot}}="z" & x \ar[r]^-{\phantom{\cdot}}="y"^-{}="c" &    v_{b-1} \ar[u]^-{\phantom{\cdot}}="d"  & v_b \ar[l] \ar[r]^-{}="i"_-{\phantom{\cdot}}="k" & u_{a+1}\ar[r]_ -{\phantom{\cdot}}="j" &      \ar @{.} [r]      & \ar[r]     & u_{a+k} \ar[d]  \\
                 &                  &        &     &    & u_{a-1} \ar[u] ^-{}="h"& \ar[l]  ^-{}="g" & \ar @{.} [l] & u_2 \ar[l]_-{\phantom{\cdot}}="m"  & u_1 \ar[l] _-{\phantom{\cdot}}="l" \\
      \ar@{.}"a";"b" \ar@{.} "c";"d"  \ar@{.}"g";"h" \ar@{.}"h";"i" \ar@{.}"j";"k" \ar@{.}"l";"m"\ar@{.}"y";"z"
                  }
}
\]

%

After applying $V_{v_{b}}$ we get:

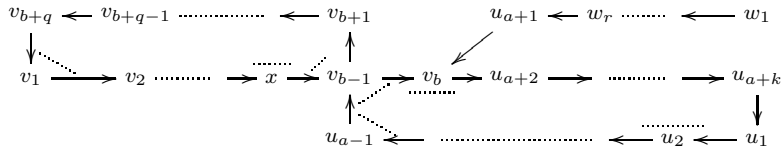
\[
\def\objectstyle{\scriptstyle}
\def\labelstyle{\scriptstyle}
\vcenter{
  \xymatrix@-1.1pc{
      v_{b+q} \ar[d]^-{}="a"   & v_{b+q-1} \ar[l]  &  & \ar@{.}[ll]  & v_{b+1} \ar[l] &  &  u_{a+1} \ar[dl] &  w_r \ar[l] & \ar@{.}[l] & w_1 \ar[l] \\
      v_1 \ar[r]^-{}="b" & v_2 \ar@{.}[r] &\ar[r]^-{\phantom{\cdot}}="z" & x \ar[r]^-{\phantom{\cdot}}="y"^-{}="c" &    v_{b-1} \ar[u]^-{\phantom{\cdot}}="d"\ar[r]^-{}="i" _-{\phantom{\cdot}}="k" & v_b  \ar[r]_ -{\phantom{\cdot}}="j" & u_{a+2}\ar[r] &      \ar @{.} [r]      & \ar[r]     & u_{a+k} \ar[d]  \\
  &   &     &    & u_{a-1} \ar[u] ^-{}="h"& \ar[l]  ^-{}="g" & & \ar @{.} [ll] & u_2 \ar[l]_-{\phantom{\cdot}}="m"  & u_1 \ar[l] _-{\phantom{\cdot}}="l" \\
      \ar@{.}"a";"b" \ar@{.} "c";"d"  \ar@{.}"g";"h" \ar@{.}"h";"i" \ar@{.}"j";"k" \ar@{.}"l";"m"\ar@{.}"y";"z"
                  }
}
\]
and by Lemma \ref{lema1} and Remark \ref{movrama} it is derived
equivalent to $B$.

\endproof

\section{Proof of the main result}\label{pruebas}

In this section we present induction proofs of Theorems I and II,
which use lemmas and propositions of Section \ref{reduc}. For the
induction basis a computer program in Groups Algorithms and
Programming (GAP) version 18-may-2000 was developed in order to
calculate all derived equivalent classes of gentle algebras $A$
with two cycles, $\phi_A=3$ and at most five vertices. The program
and the calculations are presented in www.matem.unam.mx/avella

From calculations we see all derived equivalent classes under
elementary transformations have distinct invariants except in case
of algebras defined by quivers with 4 vertices and $\phi_A= [ [ 1,
1 ], [ 1, 1 ], [ 1, 3 ] ]$, were there are two distinct classes
with the same invariant:
$$[ [ [ 1, 2, 3, 4 ], [ 5, 6 ], [ 7, 8 ] ], [ [ 1, 5 ], [ 2, 6 ],
[ 3, 7 ], [ 4, 8 ] ], [ [ 1, 1 ], [ 1, 1 ], [ 1, 3 ] ], 1 ]$$

appears in the first one and

$$[ [ [ 1, 2, 3, 4 ], [ 5, 6 ], [ 7, 8 ] ], [ [ 1, 5 ], [ 2, 7 ], [ 3, 8 ], [ 4, 6 ] ], [ [ 1, 1 ], [ 1, 1 ], [ 1, 3 ] ], 2 ]$$

in the second one. These algebras correspond to the following
quiver with relations:

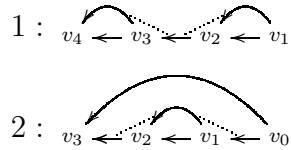
\[
1:
\def\objectstyle{\scriptstyle}
\def\labelstyle{\scriptstyle}
\vcenter{
  \xymatrix@-1.1pc{
     v_4  & v_3 \ar[l]\ar@/_1pc/[l]^(.3){}="c" & v_2 \ar[l]_-{}="b" & v_1\ar[l]  \ar@/_1pc/[l]^(.7){}="a"
      \ar@{.}"a";"b" \ar@{.} "c";"b"
                  }
}
\]

\[
2:
\def\objectstyle{\scriptstyle}
\def\labelstyle{\scriptstyle}
\vcenter{
  \xymatrix@-1.1pc{
     v_3  & v_2 \ar[l]\ar[l]_-{}="b" & v_1\ar[l] \ar@/_1pc/[l]^(.3){}="c"^(.7){}="a"  & v_0\ar[l]_-{}="d"  \ar@/_2pc/[lll]
      \ar@{.}"a";"b" \ar@{.} "c";"d"
                  }
}
\]

they are a coextension and an extension in a vertex of the algebra
defined by the quiver with relations

\[
\def\objectstyle{\scriptstyle}
\def\labelstyle{\scriptstyle}
\vcenter{
  \xymatrix@-1.1pc{
     v_3  & v_2 \ar[l]\ar[l]_-{}="b" & v_1\ar[l] \ar@/_1pc/[l]^(.3){}="c"^(.7){}="a"
      \ar@{.}"a";"b"
                  }
}
\]

the corresponding repetitive algebras, see \cite{ri97}, are
isomorphic and they have finite global dimension so they are
derived equivalent, see \cite{hw83}. However we can not transform
one into the other by using the elementary transformations.

Now we prove Theorem I by induction over the number of vertices.

Let $A=\mathrm{k}Q/ \left< \mathcal{P} \right>$ be a gentle
algebra, with $Q$ a connected quiver of two cycles. For $\#Q_0\leq
5$ we present the calculations of $\phi_A$ for each such algebra,
see www.matem.unam.mx/avella, in all cases $\#\phi_A\in\{1,3\}$.
Consider now $\#Q_0\geq$ 6. By Proposition \ref{verticex} there
exists $x\in Q_0'$ of the following type:
\begin{enumerate}
\item
\[
\def\objectstyle{\scriptstyle}
\def\labelstyle{\scriptstyle}
\vcenter{
  \xymatrix@-1.1pc{
 \ar@{-}@/^1pc/[dd] & & & & \ar@{-}@/_1pc/[dd]\\
    & u  \ar[r] ^-{\alpha} & x \ar[r]^-{\beta} & v \\
    & & & &
             }
}
\] $x$ of degree $2$
\item
\[
\def\objectstyle{\scriptstyle}
\def\labelstyle{\scriptstyle}
\vcenter{
  \xymatrix@-1.1pc{
 \ar@{-}@/^1pc/[dd] & & x \ar[d]^-{\eta} & & \ar@{-}@/_1pc/[dd]\\
    & u  \ar[r] _-{\alpha}_-{\phantom{\cdot}}="a" & w \ar[r]_-{\beta}_-{\phantom{\cdot}}="b" & v \\
    & & & &
   \ar@{.}"a";"b"
             }
}
\] with $\alpha$ and $\beta$ belonging to a cycle and $w$ of degree $3$
\item
\[
\def\objectstyle{\scriptstyle}
\def\labelstyle{\scriptstyle}
\vcenter{
  \xymatrix@-1.1pc{
 \ar@{-}@/^1pc/[dd] & & & & \ar@{-}@/_1pc/[dd]\\
    & u  \ar[r] ^-{\alpha}_-{\phantom{\cdot}}="a" & x \ar[r]^-{\beta}_-{\phantom{\cdot}}="b" & v \\
    & & & &
   \ar@{.}"a";"b"
             }
}
\] with $\alpha$ and $\beta$ belonging to a cycle and $x$ of degree $2$
\end{enumerate}
if there is no vertex $x$ of type $(1)$ or $(2)$, by Proposition
\ref{irreducibles} $\# \phi_A =1$ or $A$ is derived equivalent to
one of the representatives $(1)$, $(2)$ or $(3)$ and then
$\#\phi_{A}\in\{1,3\}$.

Consider now the case where there is a vertex $x$ of type $(1)$ or
$(2)$. Consider the algebra $A\setminus \{ x\}$. By Remark
\ref{a-xind} we know that $A\setminus \{ x\}$ is a gentle algebra
associated to a quiver of two cycles with one vertex less than the
one defining $A$; then, using induction hypothesis
$\#\phi_{A\setminus \{ x\}}\in\{1,3\}$, which means that the
algorithm to calculate the invariants of $A\setminus \{ x\}$
produces three pairs of natural numbers. Now we analyze the
differences presented in the corresponding algorithm applied to
$A$, see \cite[3]{ag07}.

Case 1.-

Observe that the process is different only when we go through
$\alpha$ and $\beta$ backwards, that is, when those arrows are
considered as part of forbidden threads. Let
$\Pi_{\alpha}=\alpha\rho_s\dots\rho_1$ and
$\Pi_{\beta}=\gamma_r\dots\gamma_1\beta$ be the forbidden threads
of $A$ involving $\alpha$ and $\beta$ resp., $\hat{H}$ be the
permitted thread such that $s(\hat{H})=s(\Pi_{\alpha})$ and
$\sigma(\hat{H})=-\sigma(\Pi_{\alpha})$, and let $H$ be the
permitted thread such that $e(H)=e(\Pi_{\beta})$ and
$\eta(H)=-\eta(\Pi_{\beta})$. In the algorithm for $A\setminus
\{x\}$ we have:

\[
\begin{matrix}
H^{A\setminus \{x\}}_i & = & H & & (\Pi^{A\setminus \{x\}}_i)^{-1} & = & (\gamma_r\dots\gamma_1\alpha^{'}\rho_s\dots\rho_1)^{-1} \\
H^{A\setminus \{x\}}_{i+1} & = & \hat{H}
\end{matrix}
\]
for some natural $i$, where $\alpha '$ being an arrow as in
Definition \ref{a-x}. This step of the algorithm produces a pair
$(n,m)$. In the corresponding algorithm for $A$ we get
\[
\begin{matrix}
H^A_{i\phantom{+1}} & = & H      & & (\Pi^A_{i\phantom{+1}})^{-1}     & = & (\Pi_{\beta})^{-1} & = & (\gamma_r\dots\gamma_1\beta)^{-1} \\
H^A_{i+1}           & = & 1_x    & & (\Pi^A_{i+1})^{-1} & = & (\Pi_{\alpha})^{-1}& = & (\alpha\rho_s\dots\rho_1)^{-1}\\
H^A_{i+2}           & = & \hat{H}
\end{matrix}
\]
and we obtain the pair $(n+1,m+1)$.

Case 2.-

Observe the process only changes when we go through $\alpha$ and
$\beta$ in the sense of the arrows, that is, when they are
considered as part of permitted threads. Let
$H_{\alpha}=\alpha\rho_s\dots\rho_1$ and
$H_{\beta}=\gamma_r\dots\gamma_1\beta$ be the permitted threads of
$A$ which involve $\alpha$ and $\beta$ resp. In the algorithm
$A\setminus \{x\}$ we have:

\[
\begin{matrix}
H^{A\setminus \{x\}}_i     & = & H_{ \alpha} & = & \alpha\rho_s\dots\rho_1 & & (\Pi_{i})^{-1} & = & 1_w \\
H^{A\setminus \{x\}}_{i+1} & = & H_{\beta}   & = &
\gamma_r\dots\gamma_1\beta
\end{matrix}
\]
for some natural $i$, and this part of the algorithm produces a
pair $(n,m)$. In the corresponding algorithm for $A$ we have
\[
\begin{matrix}
H^A_{i\phantom{+1}} & = & H_{ \alpha} & = & \alpha\rho_s\dots\rho_1 & &(\Pi^A_{i\phantom{+1}})^{-1}     & = & \eta^{-1} \\
H^A_{i+1}           & = & 1_x         &   &                         & &(\Pi^A_{i+1})^{-1} & = & 1_x\\
H^A_{i+2}           & = & \gamma_r\dots\gamma_1\beta\eta
\end{matrix}
\]
and then we get the pair $(n+1,m+1)$.

Then, $\#\phi_{A\setminus \{x\}}=\#\phi_A$, in fact the algorithm
to calculate $\phi_A$ produces exactly the same pairs but one,
which is different from the one in $\phi_{A\setminus \{x\}}$ only
by a summand $(1,1)$. So the result is completed.

Now we prove Theorem II which gives the classification of gentle
algebras with quivers of two cycles and three series of
characteristic components under derived equivalence:

For $\#Q_0\leq 5$ we have the complete classification of gentle
algebras $A=\mathrm{k}Q/ \left< \mathcal{P} \right>$, with quivers
of two cycles under derived equivalence, presented in
\\
www.matem.unam.mx/avella. We observe the result fulfills and in
fact, for $\#Q_0 = 5$ two derived equivalent algebras of that kind
can be transform one into other by using a composition of
elementary transformations.

Now consider a gentle algebra $A=\mathrm{k}Q/ \left< \mathcal{P}
\right>$, $Q$ a quiver with two cycles, $\#Q_0\geq 6$ and $\#
\phi_A=3$. We prove it is derived equivalent to one of the
representatives of the ones described in Section
\ref{representantes}. If there is $x \in Q_0$ as in $(1)$ or $(2)$
of Proposition \ref{verticex}, consider the algebra $A\setminus \{
x\}$, see Definition \ref{a-x}. We know $A\setminus\{x\}$ is a
gentle algebra with a quiver of two cycles, one vertex less than
$A$ and, by the proof of Theorem I, with $\phi_{A\setminus
\{x\}}=3$. We apply then the induction hypothesis to $A\setminus
\{ x\}$, and conclude that it is derived equivalent to one of the
representatives described in Section \ref{representantes}, call
this representative $R_{A\setminus \{ x\}}$; moreover there is a
composition of elementary transformations $T$, such that
$T(A\setminus\{x\}) = R_{A\setminus \{ x\}}$. By Corollary
\ref{cta-x}, we know $A$ is derived equivalent to an algebra
$\hat{R_A}=\mathrm{k}Q'/ \left< \mathcal{P}' \right>$ such that
the corresponding $x\in Q_0'$ is as the one described in
Definition \ref{a-x} and such that $T(A\setminus \{
x\})=\hat{R_A}\setminus\{ x \}$, so $R_{A\setminus \{
x\}}=\hat{R_A}\setminus\{ x \}$. Then $\hat{R_A}$ is obtain from
the representative $R_{A\setminus \{ x\}}$ by adding a vertex as
the one of Definition \ref{a-x}. By the results of previous
section, $\hat{R_A}$ is a representative of the ones described in
Section \ref{representantes} or can be transform into one by using
elementary transformations and this concludes the proof.

\section*{APPENDIX}\label{apendice}
\section{Elementary transformations of gentle algebras(after a manuscript by T. Holm, J. Schr\"oer, A. Zimmermann)}\label{transformaciones}

This Appendix is a slightly modified version of a technical result
in the unpublished manuscript \cite{jan01}. We show here that the
combinatorial transformations which are essential for the proof of
Theorem II are indeed derived equivalences.

\subsection{Transformations over a vertex}\label{vertice}

Let $A=\mathrm{k}Q/ \left< \mathcal{P} \right>$ be a gentle
algebra with $Q$ connected such that there exist $\alpha_1,
\alpha_2\in Q_1$ $\alpha_1\neq \alpha_2$ and
$s(\alpha_1)=s(\alpha_2)$. Denote $s(\alpha_1)$ by $i$,
$j_1:=e(\alpha_1)$ and $j_2:=e(\alpha_2)$. Suppose also that
$j_1\neq i \neq j_2$ ($j_1$ and $j_2$ could be equal). The
description of the general situation which can be found around the
vertex $i$ in the quiver with relations is the following:

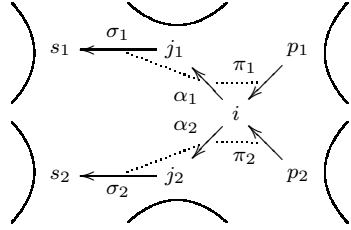
\[
\def\objectstyle{\scriptstyle}
\def\labelstyle{\scriptstyle}
\vcenter{
  \xymatrix@-1.1pc{
       \ar@{-}@/^1pc/[dd] &     & \ar@{-}@/_1pc/[rr]   &           &                  &             & \ar@{-}@/_1pc/[dd] \\
                      & s_1 &                & j_1\ar[ll]_-{\sigma_1}^-{}="b"&                  & p_1 \ar[dl] _-{\pi_1}^-{}="c"& \\
                      &     &                &           & i \ar[ul]^-{\alpha_1}^-{}="a"\ar[dl]_-{\alpha_2}^-{}="d" &             & \\
                      & s_2 &                & j_2\ar[ll]^-{\sigma_2}^-{}="e"&                  & p_2 \ar[ul] ^-{\pi_2}^-{}="f" & \\
       \ar@{-}@/_1pc/[uu] &     & \ar@{-}@/^1pc/[rr]   &           &                  &             & \ar@{-}@/^1pc/[uu]
\ar@{.}"a";"b"\ar@{.}"a";"c"\ar@{.}"d";"e"\ar@{.}"f";"d"
                  }
}
\]

We define a new algebra $V_{i}(A)=\mathrm{k}Q'/ \left<
\mathcal{P}' \right>$ with $Q_0'=Q_0$, $Q_1'$ and $\mathcal{P}'$
according to the following cases:

\begin{enumerate}
\item If $s_1\neq i \neq s_2$ let
$\alpha_m^{'},\pi_m^{'},\sigma_m^{'}$ for $m\in \{ 1,2 \}$ be
arrows such that $j_m=s(\alpha_m^{'})=e(\pi_{m}^{'})$,
$i=s(\sigma_m^{'})=e(\alpha_{m}^{'})$, $p_m=s(\pi_{3-m}^{'})$ and
$s_m=e(\sigma_m^{'})$. Define
$$Q_1':=(Q_1\setminus \{ \alpha_1 ,\alpha_2 ,\pi_1 ,\pi_2, \sigma_1 ,\sigma_2\})\cup \{\alpha_1^{'} ,\alpha_2^{'} ,\pi_1^{'} ,\pi_2^{'}, \sigma_1^{'} ,\sigma_2^{'} \}$$
and
$$\mathcal{P}':=(\mathcal{P}\setminus \{ \sigma_1 \alpha_1 ,\sigma_2 \alpha_2 ,\alpha_1 \pi_1 ,\alpha_2 \pi_2\})\cup \{\sigma_1^{'}\alpha_2^{'} ,\sigma_2^{'}\alpha_1^{'} ,\alpha_1^{'}\pi_1^{'} ,\alpha_2^{'}\pi_2^{'} \}$$
that is
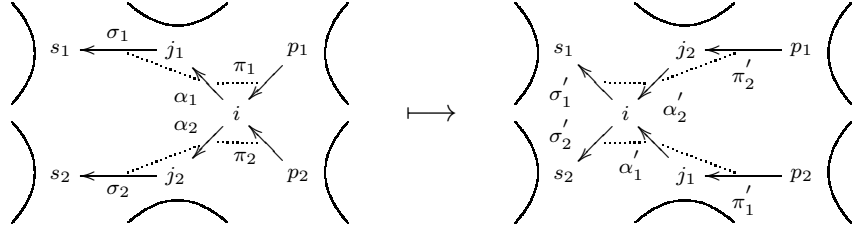
\[
\def\objectstyle{\scriptstyle}
\def\labelstyle{\scriptstyle}
\vcenter{
  \xymatrix@-1.1pc{
       \ar@{-}@/^1pc/[dd] &     & \ar@{-}@/_1pc/[rr]   &           &                  &             & \ar@{-}@/_1pc/[dd] \\
                      & s_1 &                & j_1\ar[ll]_-{\sigma_1}^-{}="b"&                  & p_1 \ar[dl] _-{\pi_1}^-{}="c"& \\
                      &     &                &           & i \ar[ul]^-{\alpha_1}^-{}="a"\ar[dl]_-{\alpha_2}^-{}="d" &             & \\
                      & s_2 &                & j_2\ar[ll]^-{\sigma_2}^-{}="e"&                  & p_2 \ar[ul] ^-{\pi_2}^-{}="f" & \\
       \ar@{-}@/_1pc/[uu] &     & \ar@{-}@/^1pc/[rr]   &           &                  &             & \ar@{-}@/^1pc/[uu]
\ar@{.}"a";"b"\ar@{.}"a";"c"\ar@{.}"d";"e"\ar@{.}"f";"d"
                  }
} \quad \longmapsto \quad
\def\objectstyle{\scriptstyle}
\def\labelstyle{\scriptstyle}
\vcenter{
  \xymatrix@-1.1pc{
       \ar@{-}@/^1pc/[dd] &                & \ar@{-}@/_1pc/[rr]             &   &                  &             & \ar@{-}@/_1pc/[dd] \\
       & s_1 &              & j_2\ar[dl]^-{\alpha_2^{'}}^-{}="a"                     &  & p_1 \ar[ll]^-{\pi_2^{'}}^-{}="e" & \\
                      &     &                          i \ar[ul]^(.7){\sigma_1^{'}}^-{}="b"\ar[dl]_(.7){\sigma_2^{'}}^-{}="d"&  &  &             & \\
                      & s_2 &              & j_1\ar[ul]^-{\alpha_1^{'}}^-{}="c"                   & & p_2 \ar[ll]^-{\pi_1^{'}}^-{}="f" & \\
       \ar@{-}@/_1pc/[uu] &                     & \ar@{-}@/^1pc/[rr]       &    &                  &             & \ar@{-}@/^1pc/[uu]
\ar@{.}"a";"b"\ar@{.}"c";"d"\ar@{.}"a";"e"\ar@{.}"c";"f"
                  }
}
\]

\item If $s_1=i$ we have
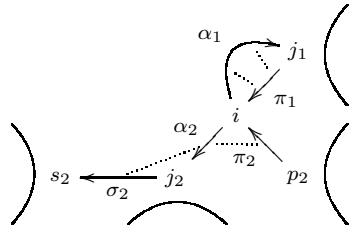
\begin{enumerate} \item
\[
\def\objectstyle{\scriptstyle}
\def\labelstyle{\scriptstyle}
\vcenter{
  \xymatrix@-1.1pc{
                      &     &                &                               &                  &                              &\ar@{-}@/_1pc/[dd] \\
                      &     &                &                               &                  & j_1 \ar[dl]^-{\pi_1}^(.6){}="c"^(.4){}="b"& \\
                      &     &                &           & i \ar@/^1.5pc/[ur]^-{\alpha_1}^(.3){}="a"^(.7){}="g"\ar[dl]_-{\alpha_2}^-{}="d" &             & \\
                      & s_2 &                & j_2\ar[ll]^-{\sigma_2}^-{}="e"&                  & p_2 \ar[ul] ^-{\pi_2}^-{}="f" & \\
       \ar@{-}@/_1pc/[uu] &     & \ar@{-}@/^1pc/[rr]   &           &                  &             & \ar@{-}@/^1pc/[uu]
\ar@{.}"a";"c"\ar@{.}"d";"e"\ar@{.}"f";"d"\ar@{.}"b";"g"
                  }
}
\]

let $\alpha_m^{'},\pi_m^{'}$ for $m\in \{ 1,2 \}$ and
$\sigma_2^{'}$ be arrows such that
$j_m=s(\alpha_m^{'})=e(\pi_{m}^{'})$,
$i=s(\sigma_2^{'})=s(\pi_{2}^{'})=e(\alpha_{m}^{'})$,
$p_2=s(\pi_{1}^{'})$ y $s_2=e(\sigma_2^{'})$. Define
$$Q_1':=(Q_1\setminus \{ \alpha_1 ,\alpha_2 ,\pi_1 ,\pi_2, \sigma_2\})\cup \{\alpha_1^{'} ,\alpha_2^{'} ,\pi_1^{'} ,\pi_2^{'}, \sigma_2^{'} \}$$
and
$$\mathcal{P}':=(\mathcal{P}\setminus \{ \pi_1 \alpha_1 ,\sigma_2 \alpha_2 ,\alpha_1 \pi_1 ,\alpha_2 \pi_2\})\cup \{\pi_2^{'}\alpha_2^{'} ,\sigma_2^{'}\alpha_1^{'} ,\alpha_1^{'}\pi_1^{'} ,\alpha_2^{'}\pi_2^{'} \}$$
that is

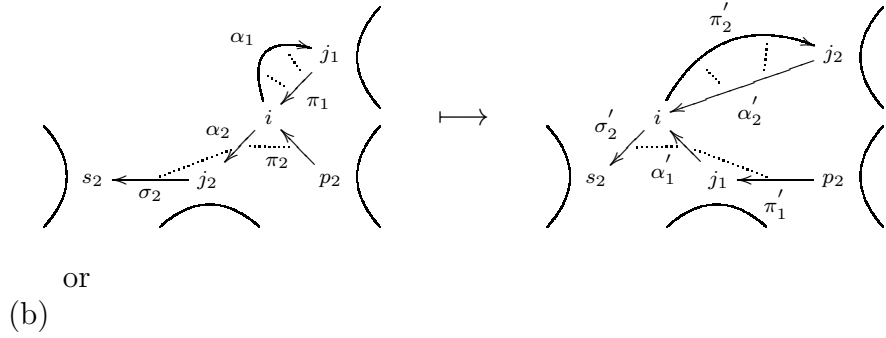
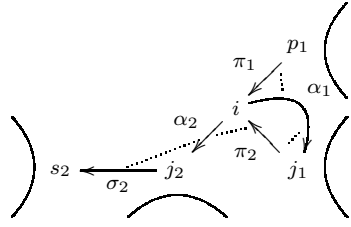
\[
\def\objectstyle{\scriptstyle}
\def\labelstyle{\scriptstyle}
\vcenter{
  \xymatrix@-1.1pc{
                      &     &                &                               &                  &                              &\ar@{-}@/_1pc/[dd] \\
                      &     &                &                               &                  & j_1 \ar[dl]^-{\pi_1}^(.6){}="c"^(.4){}="b"& \\
                      &     &                &           & i \ar@/^1.5pc/[ur]^-{\alpha_1}^(.3){}="a"^(.7){}="g"\ar[dl]_-{\alpha_2}^-{}="d" &             & \\
                      & s_2 &                & j_2\ar[ll]^-{\sigma_2}^-{}="e"&                  & p_2 \ar[ul] ^-{\pi_2}^-{}="f" & \\
       \ar@{-}@/_1pc/[uu] &     & \ar@{-}@/^1pc/[rr]   &           &                  &             & \ar@{-}@/^1pc/[uu]
\ar@{.}"a";"c"\ar@{.}"d";"e"\ar@{.}"f";"d"\ar@{.}"b";"g"
                  }
} \quad \longmapsto \quad
\def\objectstyle{\scriptstyle}
\def\labelstyle{\scriptstyle}
\vcenter{
  \xymatrix@-1.1pc{
       &         &              &                                                            &  &             & \ar@{-}@/_1pc/[dd] \\
       &         &              &                                                            &  & j_2 \ar[dlll]^-{\alpha_2^{'}}^(.6){}="c"^(.4){}="b" & \\
                      &     &           i\ar@/^1.5pc/[urrr]^-{\pi_2^{'}}^(.3){}="d"^(.7){}="a"\ar[dl]_-{\sigma_2^{'}}^-{}="f"&  &  &             & \\
                      & s_2 &              & j_1\ar[ul]^-{\alpha_1^{'}}^-{}="h"                   & & p_2 \ar[ll]^-{\pi_1^{'}}^-{}="g" & \\
       \ar@{-}@/_1pc/[uu] &                     & \ar@{-}@/^1pc/[rr]       &    &                  &             & \ar@{-}@/^1pc/[uu]
\ar@{.}"a";"b"\ar@{.}"c";"d"\ar@{.}"h";"f"\ar@{.}"g";"h"
                 }
}
\]

or  \item
\[
\def\objectstyle{\scriptstyle}
\def\labelstyle{\scriptstyle}
\vcenter{
  \xymatrix@-1.1pc{
                    &   &         &           &                  &             & \ar@{-}@/_1pc/[dd] \\
                    &   &         &           &                  & p_1 \ar[dl] _-{\pi_1}^(.3){}="e"& \\
                    &   &         &           & i \ar@/^1.5pc/[dr]^-{\alpha_1}^(.3){}="a"^(.7){}="b"\ar[dl]_-{\alpha_2}^-{}="f" &  & \\
                    &s_2&         & j_2\ar[ll]^-{\sigma_2}^-{}="g"&            & j_1 \ar[ul] ^{\pi_2}^(.3){}="c"_(.7){}="d" & \\
 \ar@{-}@/_1pc/[uu] &   & \ar@{-}@/^1pc/[rr]   &                  &            &                                            & \ar@{-}@/^1pc/[uu]
\ar@{.}"a";"e"\ar@{.}"b";"c"\ar@{.}"d";"f"\ar@{.}"f";"g"
                  }
}
\]

and we define $\alpha_m^{'},\pi_m^{'}$ for $m\in \{ 1,2 \}$ and
$\sigma_2^{'}$ arrows such that
$j_m=s(\alpha_m^{'})=e(\pi_{m}^{'})$,
$i=s(\sigma_2^{'})=s(\pi_{1}^{'})=e(\alpha_{m}^{'})$,
$p_1=s(\pi_{2}^{'})$ and $s_2=e(\sigma_2^{'})$. Define
$$Q_1':=(Q_1\setminus \{ \alpha_1 ,\alpha_2 ,\pi_1 ,\pi_2, \sigma_2\})\cup \{\alpha_1^{'} ,\alpha_2^{'} ,\pi_1^{'} ,\pi_2^{'}, \sigma_2^{'} \}$$
and
$$\mathcal{P}':=(\mathcal{P}\setminus \{ \pi_2 \alpha_1 ,\sigma_2 \alpha_2 ,\alpha_1 \pi_1 ,\alpha_2 \pi_2\})\cup \{\pi_1^{'}\alpha_2^{'} ,\sigma_2^{'}\alpha_1^{'} ,\alpha_1^{'}\pi_1^{'} ,\alpha_2^{'}\pi_2^{'} \}$$
that is

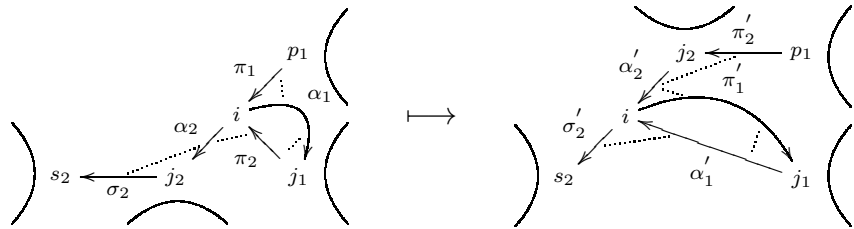
\[
\def\objectstyle{\scriptstyle}
\def\labelstyle{\scriptstyle}
\vcenter{
  \xymatrix@-1.1pc{
                    &   &         &           &                  &             & \ar@{-}@/_1pc/[dd] \\
                    &   &         &           &                  & p_1 \ar[dl] _-{\pi_1}^(.3){}="e"& \\
                    &   &         &           & i \ar@/^1.5pc/[dr]^-{\alpha_1}^(.3){}="a"^(.7){}="b"\ar[dl]_-{\alpha_2}^-{}="f" &  & \\
                    &s_2&         & j_2\ar[ll]^-{\sigma_2}^-{}="g"&            & j_1 \ar[ul] ^{\pi_2}^(.3){}="c"_(.7){}="d" & \\
 \ar@{-}@/_1pc/[uu] &   & \ar@{-}@/^1pc/[rr]   &                  &            &                                            & \ar@{-}@/^1pc/[uu]
\ar@{.}"a";"e"\ar@{.}"b";"c"\ar@{.}"d";"f"\ar@{.}"f";"g"
                  }
} \quad \longmapsto \quad
\def\objectstyle{\scriptstyle}
\def\labelstyle{\scriptstyle}
\vcenter{
  \xymatrix@-1.1pc{
       &                & \ar@{-}@/_1pc/[rr]             &   &                  &             & \ar@{-}@/_1pc/[dd] \\
       &        &              & j_2\ar[dl]_-{\alpha_2^{'}}^-{}="b"                     &  & p_1 \ar[ll]_-{\pi_2^{'}}^-{}="a" & \\
                      &     &                          i \ar[dl]_-{\sigma_2^{'}}^-{}="g" \ar@/^1.5pc/[drrr]^-{\pi_1^{'}}^(.3){}="c"^(.7){}="d"&  &  &             & \\
                      & s_2 &              & & & j_1\ar[ulll]^-{\alpha_1^{'}}^(.3){}="e"_(.7){}="f"    & \\
       \ar@{-}@/_1pc/[uu] &                     &      &    &                    &  &  \ar@{-}@/^1pc/[uu]
\ar@{.}"c";"b"\ar@{.}"d";"e"\ar@{.}"f";"g" \ar@{.}"a";"b"
                  }
}
\]
\end{enumerate}
\item If $s_1=s_2=i$ the algebra remains unchanged.
\end{enumerate}

Denote the inverse transformation by $V_i^{-1}$.

\subsection{Transformations over an arrow}\label{flecha}

Let $A=\mathrm{k}Q/ \left< \mathcal{P} \right>$ be a gentle
algebra with $Q$ connected and $\delta\in Q_1$ such that
$i:=s(\delta)$ $j:=e(\delta)$, $i\neq j$. The description of the
general situation which can be found around the arrow $\alpha$ in
the quiver with relations is the following:

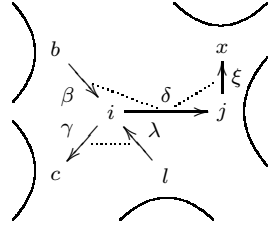
\[
\def\objectstyle{\scriptstyle}
\def\labelstyle{\scriptstyle}
\vcenter{
  \xymatrix@-1.1pc{
       \ar@{-}@/^1pc/[dd] &          &    & \ar@{-}@/_1pc/[rr] &   &              \\
                          & b \ar[dr]_-{\beta}^-{}="a"&    &                    & x &    \ar@{-}@/_1pc/[dd]                           \\
                          &          & i\ar[rr]^-{\delta}^-{}="b"\ar[dl]_-{\gamma}^-{}="e"& & j \ar[u] _-{\xi}^-{}="c"& \\
                          & c    &                &  l\ar[ul]_-{\lambda}^-{}="d"  & & \\
      \ar@{-}@/_1pc/[uu]  &  &  \ar@{-}@/^1pc/[rr]  & & &
\ar@{.}"a";"b"\ar@{.}"b";"c"\ar@{.}"d";"e"
                  }
}
\]
Let $\hat{Q}$ be the quiver obtained of $Q$ by deleting the arrows
$\beta$ and $\gamma$, that is $\hat{Q_0}=Q_0$,
$\hat{Q_1}=Q_1\setminus\{ \beta ,\gamma \}$, and
$\hat{\mathcal{P}}=\mathcal{P}\setminus \{ \delta \beta , \gamma
\lambda \}$. Consider the maximal connected subquivers of
$\hat{Q}$ which contain vertices $b$ and $c$ and let $Q^1$ be its
the union. Define also $Q^2$ as the maximal connected subquiver of
$\hat{Q}$ which contains vertex $i$. If $Q_0^1 \cap
Q_0^2=\emptyset$ define a new algebra $F_{\delta}(A)=\mathrm{k}Q'/
\left< \mathcal{P}' \right>$ taking $Q_0'=Q_0$, $Q_1'$ and
$\mathcal{P}'$ as follows. Let
$\delta^{'},\beta^{'},\gamma^{'},\lambda^{'},\xi^{'}$ be arrows
such that $j=s(\delta^{'})=e(\lambda^{'})$,
$i=s(\gamma^{'})=s(\xi^{'})=e(\beta^{'})=e(\delta^{'})$,
$b=s(\beta^{'})$, $c=e(\gamma^{'})$, $l=s(\lambda^{'})$  and
$x=e(\xi^{'})$. Define
$$Q_1':=(Q_1\setminus \{ \delta ,\beta ,\gamma ,\lambda, \xi \})\cup \{\delta^{'} ,\beta^{'} ,\gamma^{'} ,\lambda^{'}, \xi^{'} \}$$
and
$$\mathcal{P}':=(\mathcal{P}\setminus \{ \delta \beta ,\gamma \lambda ,\xi \delta \})\cup \{\gamma^{'}\delta^{'} ,\xi^{'} \beta^{'} ,\delta^{'}\lambda^{'}  \}$$
that is

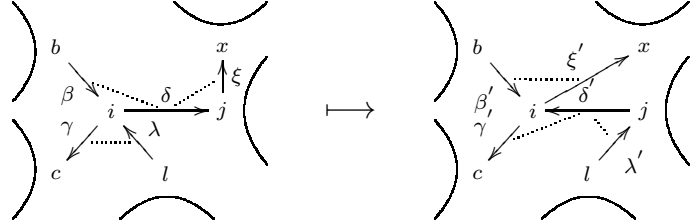
\[
\def\objectstyle{\scriptstyle}
\def\labelstyle{\scriptstyle}
\vcenter{
  \xymatrix@-1.1pc{
       \ar@{-}@/^1pc/[dd] &          &    & \ar@{-}@/_1pc/[rr] &   &              \\
                          & b \ar[dr]_-{\beta}^-{}="a"&    &                    & x &    \ar@{-}@/_1pc/[dd]                           \\
                          &          & i\ar[rr]^-{\delta}^-{}="b"\ar[dl]_-{\gamma}^-{}="e"& & j \ar[u] _-{\xi}^-{}="c"& \\
                          & c    &                &  l\ar[ul]_-{\lambda}^-{}="d"  & & \\
      \ar@{-}@/_1pc/[uu]  &  &  \ar@{-}@/^1pc/[rr]  & & &
\ar@{.}"a";"b"\ar@{.}"b";"c"\ar@{.}"d";"e"
                  }
} \quad \longmapsto \quad
\def\objectstyle{\scriptstyle}
\def\labelstyle{\scriptstyle}
\vcenter{
   \xymatrix@-1.1pc{
       \ar@{-}@/^1pc/[dd] &          &    & \ar@{-}@/_1pc/[rr] &   &             \\
                          & b \ar[dr]_-{\beta^{'}}^-{}="a"&    &                    & x &       \ar@{-}@/_1pc/[dd]                         \\
                          &          & i\ar[dl]_-{\gamma^{'}}^-{}="e"\ar[urr] ^-{\xi^{'}}^-{}="b"& & j \ar[ll]_-{\delta^{'}}^-{}="d"& \\
                          & c    &                &  l\ar[ur]_-{\lambda^{'}}^-{}="f"  & & \\
      \ar@{-}@/_1pc/[uu]  &  &  \ar@{-}@/^1pc/[rr]  & & &
\ar@{.}"a";"b"\ar@{.}"d";"f"\ar@{.}"d";"e"
                  }
}
\]
Denote the corresponding inverse transformation by
$F_{\delta'}^{-1}$.

\subsection{Transformation over a loop}\label{lazo}

Let $A=\mathrm{k}Q/ \left< \mathcal{P} \right>$ be a gentle
algebra with $Q$ connected and $\lambda\in Q_1$ be a loop of $Q$,
that is, $i:=s(\lambda)=e(\lambda)$. The description of the
general situation which can be found around the loop in the quiver
with relations is the following:

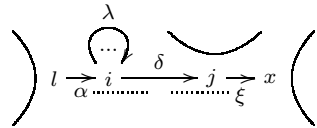
\[
\def\objectstyle{\scriptstyle}
\def\labelstyle{\scriptstyle}
\vcenter{
  \xymatrix@-1.1pc{
       \ar@{-}@/^1pc/[dd] &          &    & \ar@{-}@/_1pc/[rr] &   &  &   \ar@{-}@/_1pc/[dd]           \\
                          & l \ar[r]_-{\alpha}_-{\phantom{\cdot}}="a"&  i \ar@{}[u]|{\dots} \ar@(ul,ur)^(.2){}="c"^(.3){}="d"^-{\lambda} \ar[rr]^-{\delta}_-{\phantom{\cdot}}="b"  & & j \ar[r] _-{\xi}_-{\phantom{\cdot}}="e"  & x  &   \\
                          & & & & & &
\ar@{.}"a";"b"\ar@{.}"b";"e"
                  }
}
\]

Define a new algebra $L_{\lambda}(A)=\mathrm{k}Q'/ \left<
\mathcal{P}' \right>$ taking $Q_0'=Q_0$, $Q_1'$ y $\mathcal{P}'$
as follows. Consider $\alpha^{'},\delta^{'}, \xi^{'}$ arrows such
that $j=s(\delta^{'})=e(\alpha^{'})$,
$i=s(\xi^{'})=e(\delta^{'})$, $l=s(\alpha^{'})$ and
$x=e(\xi^{'})$. Define
$$Q_1':=(Q_1\setminus \{ \alpha ,\delta ,\xi \})\cup \{\alpha^{'} ,\delta^{'} ,\xi^{'}  \}$$
and
$$\mathcal{P}':=(\mathcal{P}\setminus \{ \delta \alpha ,\xi \delta \})\cup \{\delta^{'}\alpha^{'} ,\xi^{'}\delta^{'} \}$$
that is

\[
\def\objectstyle{\scriptstyle}
\def\labelstyle{\scriptstyle}
\vcenter{
  \xymatrix@-1.1pc{
       \ar@{-}@/^1pc/[dd] &          &    & \ar@{-}@/_1pc/[rr] &   &  &   \ar@{-}@/_1pc/[dd]           \\
                          & l \ar[r]_-{\alpha }_-{\phantom{\cdot}}="a"&  i \ar@{}[u]|{\dots} \ar@(ul,ur)^(.2){}="c"^(.3){}="d"^-{\lambda} \ar[rr]^-{\delta }_-{\phantom{\cdot}}="b"  & & j \ar[r] _-{\xi }_-{\phantom{\cdot}}="e"  & x  &   \\
                          & & & & & &
\ar@{.}"a";"b"\ar@{.}"b";"e"
                  }
} \quad \longmapsto \quad
\def\objectstyle{\scriptstyle}
\def\labelstyle{\scriptstyle}
\vcenter{
   \xymatrix@-1.1pc{
       \ar@{-}@/^1pc/[dd] &   \ar@{-}@/_1pc/[rr]  &    &  &   &  &   \ar@{-}@/_1pc/[dd]           \\
                          & l \ar[r]_-{\alpha '}_-{\phantom{\cdot}}="a"&  j \ar[rr]^-{\delta '}_-{\phantom{\cdot}}="b"  & & i \ar@{}[u]|{\dots} \ar@(ul,ur)^(.2){}="c"^(.3){}="d"^-{\lambda '} \ar[r] _-{\xi '}_-{\phantom{\cdot}}="e"  & x  &   \\
                          & & & & & &
\ar@{.}"a";"b"\ar@{.}"b";"e"
                  }
}
\]

Denote the corresponding inverse transformation by
$L_{\lambda}^{-1}$.

We say that $T$ is an {\em elementary transformation} if it is a
vertex, an arrow, a loop transformation or the inverse of one of
those. Two gentle algebra $A$ and $B$ are called {\em derived
equivalent under elementary transformations} if there exists a
finite sequence of elementary transformations $T_1$, $T_2$,\dots,
$T_r$ such that $B=T_r(\cdots T_1(A))$.

\subsection{Tilting complexes associated to the transformations}

\begin{theorem}
Let $A=\mathrm{k}Q/ \left< \mathcal{P} \right>$ be an algebra as
defined in $(1)$ or $(2)$ from Section \ref{vertice}. Consider

$$T:=T_c \oplus \bigoplus_{l\in Q_0\setminus \{i\}}P_l[1]$$

with

\[
\def\objectstyle{\scriptstyle}
\def\labelstyle{\scriptstyle}
\vcenter{
   \xymatrix@-1.1pc{
      T_c:= \ar@{.}[r] &   \ar[r]  & 0 \ar[r]   & P_{j_1}\oplus P_{j_2} \ar[r]^-{(\alpha_1,\alpha_2)} &  P_i\ar[r]  & 0 \ar[r] &  \ar@{.}[r] & .
                  }
}
\]
Then $T$ is a tilting complex and
$\operatorname{End}_{D^b(A)}(T)\simeq V_{i}(A)$.

\end{theorem}

\proof

By construction $T$ is a tilting complex and by \cite{sz} we know
that the algebra $\operatorname{End}_{D^b(A)}(T)$ is gentle, so
$\operatorname{End}_{D^b(A)}(T)=\mathrm{k}Q''/ \left<
\mathcal{P''} \right>$ with $Q_0''=Q_0$. Denote the
$\operatorname{End}_{D^b(A)}(T)$-indecomposable projective modules
for $m\in Q_0''$ by $P_m''$ , which can be identified with the
indecomposable summands of $T$. The following morphisms are
irreducible:

\begin{enumerate}
\item
\[
\alpha_1 ':= \quad
\def\objectstyle{\scriptstyle}
\def\labelstyle{\scriptstyle}
\vcenter{
   \xymatrix@-1.1pc{
       \ar@{.}[r] &   \ar[r]  & 0 \ar[r]   & P_{j_1} \ar[r] & 0 \ar[r]    & 0 \ar[r] &  \ar@{.}[r] & \\
       \ar@{.}[r] &   \ar[r]  & 0 \ar[r]   & P_{j_1}\oplus P_{j_2} \ar[u]^-{(id,0)}\ar[r]^-{(\alpha_1,\alpha_2)}^-{}="b" & P_{i}\ar[u]^-{}="a" \ar[r]  & 0\ar[u] \ar[r] &  \ar@{.}[r] &
                  }
}
\]
\[
\alpha_2 ':= \quad
\def\objectstyle{\scriptstyle}
\def\labelstyle{\scriptstyle}
\vcenter{
   \xymatrix@-1.1pc{
       \ar@{.}[r] &   \ar[r]  & 0 \ar[r]   & P_{j_2}  \ar[r] & 0 \ar[r]    & 0 \ar[r] &  \ar@{.}[r] & \\
       \ar@{.}[r] &   \ar[r]  & 0 \ar[r]   &  P_{j_1}\oplus P_{j_2} \ar[u]^-{(0,id)}\ar[r]^-{(\alpha_1,\alpha_2)}^-{}="b" & P_{i} \ar[u]^-{}="a"\ar[r]  & 0\ar[u] \ar[r] &  \ar@{.}[r] &
                  }
}
\]

\[
\pi_1 ':=\phantom{aaa:} \quad
\def\objectstyle{\scriptstyle}
\def\labelstyle{\scriptstyle}
\vcenter{
   \xymatrix@-1.1pc{
       \ar@{.}[r] &   \ar[r]   & 0 \ar[r] & P_{p_2} \ar[r]    & 0 \ar[r] & 0 \ar[r]  &  \ar@{.}[r] & \\
       \ar@{.}[r] &   \ar[r]   & 0 \ar[r]^-{}="b"     & P_{j_1} \ar[u]^-{\alpha_1\pi_2}\ar[r] & 0\ar[u] \ar[r] & 0 \ar[r]   \ar[u]     &  \ar@{.}[r] &
                  }
}
\]

\[
\pi_2 ':=\phantom{aaa: } \quad
\def\objectstyle{\scriptstyle}
\def\labelstyle{\scriptstyle}
\vcenter{
   \xymatrix@-1.1pc{
       \ar@{.}[r] &    \ar[r]   & 0 \ar[r] & P_{p_1} \ar[r]    & 0 \ar[r] & 0 \ar[r]  & \ar@{.}[r] & \\
       \ar@{.}[r] &   \ar[r]   & 0 \ar[r]^-{}="b" & P_{j_2}\ar[u]^-{\alpha_2\pi_1} \ar[r] & 0 \ar[u]\ar[r]  & 0 \ar[r]\ar[u]  & \ar@{.}[r] &
                  }
}
\]
\[
\sigma_1 ':= \quad
\def\objectstyle{\scriptstyle}
\def\labelstyle{\scriptstyle}
\vcenter{
   \xymatrix@-1.1pc{
       \ar@{.}[r] &   \ar[r]  & 0 \ar[r]   &  P_{j_1}\oplus P_{j_2} \ar[r]^-{(\alpha_1,\alpha_2)}^-{}="b" & P_i  \ar[r]  & 0 \ar[r] &  \ar@{.}[r] & \\
       \ar@{.}[r] &   \ar[r]  & 0 \ar[r]   & P_{s_1} \ar[u]^-{\binom{\sigma_1}{0}}^-{}="a"\ar[r] & 0 \ar[u]\ar[r]  & 0 \ar[r]\ar[u] &  \ar@{.}[r] &
\ar@{.}"a";"b"
                  }
}
\]

\[
\sigma_2 ':= \quad
\def\objectstyle{\scriptstyle}
\def\labelstyle{\scriptstyle}
\vcenter{
   \xymatrix@-1.1pc{
       \ar@{.}[r] &   \ar[r]  & 0 \ar[r]   &  P_{j_1}\oplus P_{j_2} \ar[r]^-{(\alpha_1,\alpha_2)}^-{}="b" & P_i  \ar[r]  & 0 \ar[r] &  \ar@{.}[r] & \\
       \ar@{.}[r] &   \ar[r]  & 0 \ar[r]   & P_{s_2}\ar[u]^-{\binom{0}{ \sigma_2}}^-{}="a" \ar[r] & 0 \ar[u]\ar[r]  & 0 \ar[u]\ar[r] &  \ar@{.}[r] &
\ar@{.}"a";"b"
                  }
}
\]
because they cannot be factorized through any other indecomposable
summand of $T$. The morphism $\sigma_1 : P_{s_1}\rightarrow
P_{j_1}$ and $\sigma_2 : P_{s_2}\rightarrow P_{j_2}$ are
irreducible in $A$, but factorize through $\alpha _1 '$ and
$\alpha _2 '$ respectively. Also, the compositions $\sigma_1 '
\alpha_2 '$ and $ \sigma_2 ' \alpha_1 '$ are zero and $ \alpha_1 '
\pi_1 '$, $\alpha_2 '\pi_2 '$ are homotopic to zero.

To see there are no further irreducible maps we use Happel's
formula, see \cite[III.1.3 and III.1.4]{ha88}:
$$\sum_{i}(-1)^i \operatorname{dim Hom}_{K^b(A)}(Q,R[i])=\sum_{r,s}(-1)^{r-s} \operatorname{dim Hom}_{A}(Q^r,R^s)$$

where $K^b(A)$ is the homotopy category of bounded complexes of
projective $A$-modules, $[.]$ is the shift operator and
$Q=(Q^r)_{r\in \mathbb{Z}}$, $R=(R^s)_{s\in \mathbb{Z}}$ are
objects in $K^b(A)$. For direct summands of tilting complexes
$\operatorname{dim Hom}_{K^b(A)}(Q,R[i])=0$ for $i\neq 0$ so in
this case
$$\operatorname{dim Hom}_{K^b(A)}(Q,R)=\sum_{r,s}(-1)^{r-s} \operatorname{dim Hom}_{A}(Q^r,R^s)$$

For $x,y,z\in Q_0$ denote by $(x,y)$ the set of paths starting in
$y$ and ending in $x$, by $(x,y;z)$ the set of paths starting in
$y$ and ending in $x$ which go through $z$ and by $(x,y;\hat{z})$
the set of paths starting in $y$ and ending in $x$ which do not go
through $z$. Also for $\alpha \in Q_1$ denote by $(x,y;\alpha)$
the set of paths starting in $y$ and ending in $x$ which pass
through $\alpha$ and by $(x,y;\hat{\alpha})$ the set of paths
starting in $y$ and ending in $x$ which do not go through
$\alpha$.

 By using Happel's formula we calculate the following:
\begin{align*}
\operatorname{dim Hom}_{K^b(A)}(T_c,P_l[1])&=\operatorname{dim
Hom}_{A}(P_{j_1}\oplus
P_{j_2},P_l)-\operatorname{dim Hom}_{A}(P_i, P_l)\\
&=\#(j_1,l)+\#(j_2,l)-\#(i,l)\\
&=\#(j_1,l;i)+\#(j_1,l;\hat{i})+\#(j_2,l;i)\\
&+\#(j_2,l;\hat{i})-\#(i,l;\pi_1)-\#(i,l;\hat{\pi_1})\\
&=\#(j_1,l;\hat{i})+\#(j_2,l;\hat{i})
\end{align*}
because $\#(j_1,l;i)=\#(i,l;\hat{\pi_1})$ and
$\#(j_2,l;i)=\#(i,l;\pi_1)$, but only the trivial paths in $j_1$
and $j_2$ produce irreducible morphisms, the rest factorize
through these.

\begin{align*} \operatorname{dim Hom}_{K^b(A)}(P_l[1],
T_c)&=\operatorname{dim Hom}_{A}(P_l,P_{j_1}\oplus
P_{j_2})-\operatorname{dim Hom}_{A}(P_l, P_i)\\
&=\#(l,j_1)+\#(l,j_2)-\#(l,i)\\
&=\#(l,j_1;\sigma_1)+\#(l,j_1;\hat{\sigma_1})+\#(l,j_2;\sigma_2)\\
&+\#(l,j_2;\hat{\sigma_2})-\#(l,i;\alpha_1)-\#(l,i;\alpha_2)\\
&=\#(l,j_1;\sigma_1)+\#(l,j_2;\sigma_2)
\end{align*}
because $\#(l,j_1;\hat{\sigma_1})=\#(l,i;\alpha_1)$ and
$\#(l,j_2;\hat{\sigma_2})=\#(l,i;\alpha_2)$, and the only
irreducible morphisms are the induced by $\sigma_1$ and $\sigma_2$
\begin{align*}
\operatorname{dim
Hom}_{K^b(A)}(P_{j_1}[1],P_l[1])&=\operatorname{dim
Hom}_{A}(P_{j_1},P_l)\\
& =\#(j_1,l)
\end{align*}

so the only irreducible morphisms in this case are the ones
corresponding to $\alpha_1 \pi_2$ for $l=p_2$ or the one induced
by an arrow connecting $l$ and $j_1$. The calculation for
$\operatorname{dim Hom}_{K^b(A)}(P_{j_2}[1],P_l[1])$ is similar.
%
\item
\begin{enumerate}
\item Let $\alpha_1 '$, $\alpha_2 '$, $\pi_1 '$ and $\sigma_2'$ be
as in $(1)$, and

\[
\pi_2 ':= \quad
\def\objectstyle{\scriptstyle}
\def\labelstyle{\scriptstyle}
\vcenter{
   \xymatrix@-1.1pc{
       \ar@{.}[r] &   \ar[r]  & 0 \ar[r]   & P_{j_1}\oplus P_{j_2}  \ar[r]^-{(\alpha_1,\alpha_2)}^-{}="b" &  P_i \ar[r]  & 0 \ar[r]&  \ar@{.}[r] & \\
       \ar@{.}[r] &   \ar[r]  &  0 \ar[r]   & P_{j_2} \ar[u]^-{\binom{\alpha_2\pi_1}{0}}^-{}="a"\ar[r] & 0\ar[r]\ar[u]  & 0 \ar[r] \ar[u] &  \ar@{.}[r] &
\ar@{.}"a";"b"
                  }
}
\]
which is also an irreducible morphism because it cannot be
factorized through any other indecomposable summand of $T$. The
morphism $\sigma_2 :P_{s_2}\rightarrow P_{j_2}$ is irreducible in
$A$, but induces a morphism which can be factorized through
$\alpha _2 '$, while $\alpha_2 \pi_1 $ induces a morphism which
can be factorized through $\alpha_1 '$. Also, the compositions $
\sigma_2 ' \alpha_1 '$ and $\pi_2 '\alpha_2 '$ are zero and $
\alpha_1 ' \pi_1 '$ is homotopic to zero as in $(1)$ such as
$\alpha_2 '\pi_2 '$.
\item Let $\alpha_1 '$, $\alpha_2 '$, $\pi_2 '$ and $\sigma_2'$ as
in $(1)$, and

\[
\pi_1 ':= \quad
\def\objectstyle{\scriptstyle}
\def\labelstyle{\scriptstyle}
\vcenter{
   \xymatrix@-1.1pc{
       \ar@{.}[r] &   \ar[r]  & 0 \ar[r]   &  P_{j_1}\oplus P_{j_2}  \ar[r]^-{(\alpha_1,\alpha_2)}^-{}="b" & P_i \ar[r]  & 0 \ar[r] &  \ar@{.}[r] & \\
       \ar@{.}[r] &   \ar[r]  & 0 \ar[r]   & P_{j_1} \ar[r]\ar[u]^-{\binom{\alpha_1\pi_2}{0}}^-{}="a" & 0\ar[u] \ar[r]  & 0 \ar[u]\ar[r] &  \ar@{.}[r] &
\ar@{.}"a";"b"
                  }
}
\]
which is also an irreducible morphism because can not be
factorized through any other indecomposable summand of $T$. The
morphism $\sigma_2 :P_{s_2}\rightarrow P_{j_2}$ which is
irreducible in $A$, induces a morphism which factorize through
$\alpha _2 '$, while $\alpha_1 \pi_2$ induces a morphism which
factorizes through $\alpha_1 '$. Also, the compositions $ \sigma_2
' \alpha_1 '$ and $\pi_1 '\alpha_2 '$ are zero and $ \alpha_2 '
\pi_2 '$ is homotopic to zero as in $(1)$ such as $\alpha_1 '\pi_1
'$.
\end{enumerate}
\end{enumerate}
For each case, the other irreducible morphisms in
$A-\operatorname{mod}$ produce irreducible morphisms in
$\operatorname{End}_{D^b(A)}(T)-\operatorname{mod}$ because the
corresponding projective indecomposable modules are direct
summands of $T$ of degree one. We verify there are no further
irreducible maps using Happel's formula:

(a)
\begin{align*}
\operatorname{dim Hom}_{K^b(A)}(T_c,P_l[1])&=\operatorname{dim
Hom}_{A}(P_{j_1}\oplus
P_{j_2},P_l)-\operatorname{dim Hom}_{A}(P_i, P_l)\\
&=\#(j_1,l)+\#(j_2,l)-\#(i,l)\\
&=\#(j_1,l;i)+\#(j_1,l;\hat{i})+\#(j_2,l;i)\\
&+\#(j_2,l;\hat{i})-\#(i,l;\pi_1)-\#(i,l;\hat{\pi_1})\\
&=\#(j_1,l;\hat{i})+\#(j_2,l;\hat{i})
\end{align*}
because $\#(j_1,l;i)=\#(i,l;\hat{\pi_1})$ and
$\#(j_2,l;i)=\#(i,l;\pi_1)$, but only the trivial paths in $j_1$
and $j_2$ produce irreducible morphisms, the rest factorize
through these.

\begin{align*} \operatorname{dim Hom}_{K^b(A)}(P_l[1],
T_c)&=\operatorname{dim Hom}_{A}(P_l,P_{j_1}\oplus
P_{j_2})-\operatorname{dim Hom}_{A}(P_l, P_i)\\
&=\#(l,j_1)+\#(l,j_2)-\#(l,i)\\
&=\#(l,j_1;i)+\#(l,j_1;\hat{i})+\#(l,j_2;\sigma_2)\\
&+\#(l,j_2;\hat{\sigma_2})-\#(l,i;j_2)-\#(l,i;\hat{j_2})\\
&=\#(l,j_1;i)+\#(l,j_2;\sigma_2)
\end{align*}
because $\#(l,j_1;\hat{i})=\#(l,i;\hat{j_2})$ and
$\#(l,j_2;\hat{\sigma_2})=\#(l,i;j_2)$, and the only irreducible
morphisms are the induced by $\alpha_2\pi_1$ (in this case
$l=j_2$) and $\sigma_2$ (in this case $l=s_2$).
\begin{align*}
\operatorname{dim
Hom}_{K^b(A)}(P_{j_1}[1],P_l[1])&=\operatorname{dim
Hom}_{A}(P_{j_1},P_l)\\
& =\#(j_1,l)
\end{align*}

so this are the only irreducible morphisms in this case are the
one corresponding to $\alpha_1 \pi_2$ for $l=p_2$ or the one
induced by an arrow connecting $l$ and $j_1$.

The calculation for $\operatorname{dim
Hom}_{K^b(A)}(P_{j_2}[1],P_l[1])$ is similar.

(b)
\begin{align*}
\operatorname{dim Hom}_{K^b(A)}(T_c,P_l[1])&=\operatorname{dim
Hom}_{A}(P_{j_1}\oplus
P_{j_2},P_l)-\operatorname{dim Hom}_{A}(P_i, P_l)\\
&=\#(j_1,l)+\#(j_2,l)-\#(i,l)\\
&=\#(j_1,l;i)+\#(j_1,l;\hat{i})+\#(j_2,l;i)\\
&+\#(j_2,l;\hat{i})-\#(i,l;\pi_1)-\#(i,l;\hat{\pi_1})\\
&=\#(j_1,l;\hat{i})+\#(j_2,l;\hat{i})
\end{align*}
because $\#(j_1,l;i)=\#(i,l;\hat{\pi_1})$ and
$\#(j_2,l;i)=\#(i,l;\pi_1)$ but only the trivial paths in $j_1$
and $j_2$ produce irreducible morphisms, the rest factorize
through these.

\begin{align*} \operatorname{dim Hom}_{K^b(A)}(P_l[1],
T_c)&=\operatorname{dim Hom}_{A}(P_l,P_{j_1}\oplus
P_{j_2})-\operatorname{dim Hom}_{A}(P_l, P_i)\\
&=\#(l,j_1)+\#(l,j_2)-\#(l,i)\\
&=\#(l,j_1;i)+\#(l,j_1;\hat{i})+\#(l,j_2;\sigma_2)\\
&+\#(l,j_2;\hat{\sigma_2})-\#(l,i;\alpha_2)-\#(l,i;\hat{\alpha_2})\\
&=\#(l,j_1;i)+\#(l,j_2;\sigma_2)
\end{align*}
because $\#(l,j_1;\hat{i})=\#(l,i;\hat{\alpha_2})$ and
$\#(l,j_2;\hat{\sigma_2})=\#(l,i;\alpha_2)$, and the only
irreducible morphisms are the induced by $\alpha_1\pi_2$ (in this
case $l=j_1$) and $\sigma_2$ (in this case $l=s_2$).
\begin{align*}
\operatorname{dim
Hom}_{K^b(A)}(P_{j_1}[1],P_l[1])&=\operatorname{dim
Hom}_{A}(P_{j_1},P_l)\\
& =\#(j_1,l)
\end{align*}

so this are the only irreducible morphisms in this case are the
one corresponding to $\alpha_1 \pi_2$ for $l=p_2$ or the one
induced by an arrow connecting $l$ and $j_1$.

The calculation for $\operatorname{dim
Hom}_{K^b(A)}(P_{j_2}[1],P_l[1])$ is similar.

Then $\operatorname{End}_{D^b(A)}(T)$ can be identified with
$V_i(A)$.

\endproof

\begin{remark}{\normalfont
In case $(3)$ of Section \ref{vertice}, using the previous tilting
complex, $\operatorname{End}_{D^b(A)}(T)$ identifies with an
algebra defined by the same quiver with relations which defines
$A$. }\end{remark}

\begin{theorem}
Let $A=\mathrm{k}Q/ \left< \mathcal{P} \right>$ as in Section
\ref{flecha}. Define

$$T:=T_c\oplus \bigoplus_{m\in Q_0^1}P_m \oplus \bigoplus_{r\in Q_0^2\setminus \{i\}}P_r[1]$$

with

\[
\def\objectstyle{\scriptstyle}
\def\labelstyle{\scriptstyle}
\vcenter{
   \xymatrix@-1.1pc{
      T_c:= \ar@{.}[r] &   \ar[r]  & 0 \ar[r]   & P_j \ar[r]^-{\delta} & P_i \ar[r]  & 0 \ar[r] &  \ar@{.}[r] & .
                  }
}
\]
Then $T$ is a tilting complex and
$\operatorname{End}_{D^b(A)}(T)\simeq F_{\delta}(A)$.

\end{theorem}

\proof

By construction $T$ is a tilting complex. By \cite{sz} we know the
algebra $\operatorname{End}_{D^b(A)}(T)$ is gentle, we write then
$\operatorname{End}_{D^b(A)}(T)=\mathrm{k}Q''/ \left<
\mathcal{P''} \right>$ with $Q_0''=Q_0$. Denote the
$\operatorname{End}_{D^b(A)}(T)$-indecomposable projective modules
for $m\in Q_0''$ by $P_m''$, which are identified with the
indecomposable summands of $T$. We analyze the irreducible
morphisms in $\operatorname{End}_{D^b(A)}(T)-\operatorname{mod}$.
Some of them are the following:

\[
\beta ':= \quad
\def\objectstyle{\scriptstyle}
\def\labelstyle{\scriptstyle}
\vcenter{
   \xymatrix@-1.1pc{
       \ar@{.}[r] &   \ar[r]  & 0 \ar[r]   & 0\ar[r] &  P_b \ar[r]    & 0 \ar[r] &  \ar@{.}[r] & \\
       \ar@{.}[r] &   \ar[r]  & 0 \ar[r]   & P_j \ar[u]\ar[r]^-{\delta}^-{}="b" & P_i \ar[u]_-{\beta}^-{}="a"\ar[r]  & 0 \ar[u]\ar[r] &  \ar@{.}[r] &
\ar@{.}"a";"b"
                  }
}
\]

\[
\gamma ':= \quad
\def\objectstyle{\scriptstyle}
\def\labelstyle{\scriptstyle}
\vcenter{
   \xymatrix@-1.1pc{
       \ar@{.}[r] &   \ar[r]  & 0 \ar[r]   & P_j \ar[r]^-{\delta} & P_i \ar[r]  & 0 \ar[r] &  \ar@{.}[r] & \\
       \ar@{.}[r] &   \ar[r]  & 0 \ar[r]   & 0\ar[u]  \ar[r]^-{}="b" & P_c  \ar[r]\ar[u]_-{\gamma}^-{}="a"  & 0 \ar[u]\ar[r] &  \ar@{.}[r] &
                  }
}
\]

\[
\delta ':= \quad
\def\objectstyle{\scriptstyle}
\def\labelstyle{\scriptstyle}
\vcenter{
   \xymatrix@-1.1pc{
       \ar@{.}[r] &   \ar[r]  & 0 \ar[r]   & P_j \ar[r] & 0   \ar[r]  & 0 \ar[r] &  \ar@{.}[r] & \\
       \ar@{.}[r] &   \ar[r]  & 0 \ar[r]   & P_j \ar[u]^-{id}^-{}="a"\ar[r]^-{\delta}^-{}="b" & P_i \ar[u]\ar[r]  & 0 \ar[u]\ar[r] &  \ar@{.}[r] &
                  }
}
\]

\[
\xi ':= \quad
\def\objectstyle{\scriptstyle}
\def\labelstyle{\scriptstyle}
\vcenter{
   \xymatrix@-1.1pc{
       \ar@{.}[r] & \ar[r]& 0 \ar[r]   & P_j\ar[r]^-{\delta}="a" & P_i \ar[r]   & 0 \ar[r]\ar[d] &  \ar@{.}[r] & \\
       \ar@{.}[r] & \ar[r]& 0 \ar[r]   &P_x  \ar[u]^-{\xi}="b" \ar[r] & 0\ar[u]\ar[r]  & 0 \ar[r] &  \ar@{.}[r] &
\ar@{.}"a";"b"
                  }
}
\]

\[
\phantom{a:}\lambda ':= \quad
\def\objectstyle{\scriptstyle}
\def\labelstyle{\scriptstyle}
\vcenter{
   \xymatrix@-1.1pc{
       \ar@{.}[r] & \ar[r]    & 0 \ar[r] & 0 \ar[r]& P_l \ar[r]& 0 \ar[r]  &  \ar@{.}[r] & \\
       \ar@{.}[r] & \ar[r]    & 0 \ar[r] & 0 \ar[r]\ar[u]   & P_j \ar[u]^-{\delta\lambda}^-{}="a" \ar[r]^-{}="b" & 0 \ar[u]\ar[r]  &  \ar@{.}[r] &
                  }
}
\]
they do not factorize through any other indecomposable summand of
$T$ and are therefore irreducible. The morphism $\xi
:P_x\rightarrow P_j$ is irreducible in $A$, but produces a
morphism which factorizes through $\delta '$. Also, compositions
$\xi ' \beta '$ and $ \gamma ' \delta '$ are zero and  $ \delta '
\lambda '$ is homotopic to zero.

The rest of the irreducible morphisms in $A-\operatorname{mod}$
produce irreducible morphisms in
$\operatorname{End}_{D^b(A)}(T)-\operatorname{mod}$ because the
corresponding indecomposable projective modules are direct
summands of $T$ in the suitable degrees.

We verify there are no further irreducible maps using Happel's
formula. We calculate the following:
\begin{align*}
\operatorname{dim Hom}_{K^b(A)}(T_c,P_m)&=\operatorname{dim
Hom}_{A}(P_{i},P_m)-\operatorname{dim Hom}_{A}(P_j, P_m)\\
&=\#(i,m)
\end{align*}
but only $\beta$ produce an irreducible morphism, for $m=b$.

\begin{align*} \operatorname{dim Hom}_{K^b(A)}(P_r[1],
T_c)&=\operatorname{dim Hom}_{A}(P_r,P_{j})-\operatorname{dim Hom}_{A}(P_r, P_i)\\
&=\#(r,j)-\#(r,i)
=\#(r,j;x)+\#(r,j;\hat{x})-\#(r,i;j)\\
&=\#(r,j;x)
\end{align*}
because $\#(r,j;\hat{x})=\#(r,i;j)$, and the only irreducible
morphism is the one induced by $\xi$ for $r=x$.
\begin{align*}
\operatorname{dim Hom}_{K^b(A)}(P_m,T_c)&=\operatorname{dim
Hom}_{A}(P_{m},P_i)-\operatorname{dim Hom}_{A}(P_m, P_j)\\
&=\#(m,i)-\#(m,j)=\#(m,i)
\end{align*}
but only $\gamma$ produce an irreducible morphism, for $m=c$.
\begin{align*}
\operatorname{dim Hom}_{K^b(A)}(T_c,P_r[1])&=\operatorname{dim
Hom}_{A}(P_{j},P_r)-\operatorname{dim Hom}_{A}(P_i, P_r)\\
&=\#(j,r)-\#(i,r)=\#(j,r;i)+\#(j,r;\hat{i})-\#(i,r)\\
&=\#(j,r;\hat{i})
\end{align*}
because $\#(j,r;i)=\#(i,r)$ but the only irreducible morphism is
the one induced by the trivial path in $j$ for $r=j$, the others
factorize through this one.

Then $\operatorname{End}_{D^b(A)}(T)$ identifies with
$F_{\delta}(A)$.

%

So these are the only irreducible morphisms.
\endproof

\begin{theorem}
Let $A=\mathrm{k}Q/ \left< \mathcal{P} \right>$ as in Section
\ref{lazo}. Define

$$T:=T_c\oplus \bigoplus_{m\in Q_0\setminus \{i\}}P_m[1]$$

with

\[
\def\objectstyle{\scriptstyle}
\def\labelstyle{\scriptstyle}
\vcenter{
   \xymatrix@-1.1pc{
      T_c:= \ar@{.}[r] &   \ar[r]  & 0 \ar[r]   & P_j\oplus P_j  \ar[r]^-{(\delta,\delta \lambda)} &P_i \ar[r]  & 0 \ar[r] &  \ar@{.}[r] & .
                   }
}
\]
Then $T$ is a tilting complex and
$\operatorname{End}_{D^b(A)}(T)\simeq L_{\lambda}(A)$.

\end{theorem}

\proof

By construction $T$ is a tilting complex and by \cite{sz} the
algebra $\operatorname{End}_{D^b(A)}(T)$ is gentle, then
$\operatorname{End}_{D^b(A)}(T)=\mathrm{k}Q''/ \left<
\mathcal{P''} \right>$ with $Q_0''=Q_0$. Denote by $P_m''$ the
$\operatorname{End}_{D^b(A)}(T)$-indecomposable projective modules
for $m\in Q_0''$, which are identified with the indecomposable
summands of $T$. The following morphisms:

\[
\alpha ':=\phantom{aa:} \quad
\def\objectstyle{\scriptstyle}
\def\labelstyle{\scriptstyle}
\vcenter{
   \xymatrix@-1.1pc{
       \ar@{.}[r] & \ar[r]  & 0 \ar[r]  & 0 \ar[r] & P_l \ar[r]    & 0 \ar[r] &  \ar@{.}[r] & \\
       \ar@{.}[r] & \ar[r]  & 0 \ar[r]  & 0 \ar[r]\ar[u]^-{}="a"^-{}="b" & P_j \ar[u]^-{\delta\lambda\alpha}\ar[r]  & 0 \ar[u]\ar[r] &  \ar@{.}[r] &
                  }
}
\]

\[
\delta ':= \quad
\def\objectstyle{\scriptstyle}
\def\labelstyle{\scriptstyle}
\vcenter{
   \xymatrix@-1.1pc{
       \ar@{.}[r] &   \ar[r]  & 0 \ar[r]   & P_j \ar[r] & 0 \ar[r] & 0 \ar[r] &  \ar@{.}[r] & \\
       \ar@{.}[r] &   \ar[r]  & 0 \ar[r]   & P_j\oplus P_j \ar[u]^-{(0,id)}="a"\ar[r]^-{(\delta,\delta\lambda)}^-{}="b" &P_i  \ar[u]\ar[r]  & 0 \ar[r]\ar[u] &  \ar@{.}[r] &
                  }
}
\]

\[
\xi ':= \quad
\def\objectstyle{\scriptstyle}
\def\labelstyle{\scriptstyle}
\vcenter{
   \xymatrix@-1.1pc{
       \ar@{.}[r] & \ar[r]& 0 \ar[r] & P_j\oplus P_j \ar[r]^-{(\delta,\delta\lambda)}^-{}="a" & P_i \ar[r]& 0 \ar[r] &  \ar@{.}[r] & \\
       \ar@{.}[r] & \ar[r]& 0 \ar[r] & P_x \ar[u]^-{\binom{\xi}{0}}="b"\ar[r]^-{} & 0 \ar[u]_-{}\ar[r]  & 0 \ar[u]\ar[r] &  \ar@{.}[r] &
\ar@{.}"a";"b"
                  }
}
\]
are irreducible because they cannot be factorized through any
other indecomposable summand of $T$. The morphism $\xi
:P_x\rightarrow P_j$ is irreducible in $A$, but the corresponding
morphism in $\operatorname{End}_{D^b(A)}(T)$ factorizes through
$\delta '$ and the morphism induced by the projection of
$P_j\oplus P_j$ over the first component, factorizes through

\[
\lambda':= \quad
\def\objectstyle{\scriptstyle}
\def\labelstyle{\scriptstyle}
\vcenter{
   \xymatrix@-1.1pc{
       \ar@{.}[r] & \ar[r]& 0 \ar[r] & P_j\oplus P_j \ar[r]^-{(\delta,\delta\lambda)}^-{}="c" & P_i \ar[r]& 0 \ar[r] &  \ar@{.}[r] & \\
       \ar@{.}[r] & \ar[r]& 0 \ar[r] & P_j\oplus P_j \ar[r]^-{(\delta,\delta\lambda)}^-{}="b"\ar[u]^-{\binom{0\phantom{0}0}{1\phantom{0}0}}="d" & P_i\ar[r]\ar[u]_-{\lambda}^-{}="a"  & 0 \ar[u]\ar[r] &  \ar@{.}[r] &
                  }
}
\]

Also, notice that the compositions $ \xi ' \delta '$ and $\lambda
'\lambda '$ are zero and $ \delta ' \alpha'$ is homotopic to zero.

The other morphisms which are irreducible in
$A-\operatorname{mod}$ produce irreducible morphisms in
$\operatorname{End}_{D^b(A)}(T)-\operatorname{mod}$ because the
corresponding indecomposable projective modules are direct
summands of $T$ of degree one.

By using Happel's formula we calculate the following:
\begin{align*}
\operatorname{dim Hom}_{K^b(A)}(T_c,P_m[1])&=\operatorname{dim
Hom}_{A}(P_{j}\oplus P_j,P_m)-\operatorname{dim Hom}_{A}(P_i, P_m)\\
&=\#2(j,m)-\#(i,m)\\
&=2\#(j,m;i)+2\#(j,m;\hat{i})-\#(i,m) =2\#(j,m;\hat{i})
\end{align*}
because $\#(i,m)=2\#(j,m;i)$ but only the trivial path in $j$
induce an irreducible morphism (applied in the second coordinate)
for $m=j$, the rest factorize through this one.

\begin{align*} \operatorname{dim Hom}_{K^b(A)}(P_m[1],
T_c)&=\operatorname{dim Hom}_{A}(P_m,P_{j}\oplus P_j)-\operatorname{dim Hom}_{A}(P_m, P_i)\\
&=2\#(m,j)-\#(m,i) \\
&=2\#(m,j;x)+2\#(m,j;\hat{x})-\#(m,i) \\
&=2\#(m,j;x)
\end{align*}
because $2\#(m,j;\hat{x})=\#(m,i)$, and the only irreducible
morphism is the one induced by $\xi$ for $m=x$.

\begin{align*}
\operatorname{dim Hom}_{K^b(A)}(P_m[1],P_t[1])&=\operatorname{dim
Hom}_{A}(P_{m},P_t)=\#(m,t)
\end{align*}
but the only irreducible morphism are the ones induced by arrows
starting in $l$ and ending in $m$ and the one induced by the path
$\delta\lambda\alpha$ for $t=l$ and $m=j$.

Finally
\begin{align*}
\operatorname{dim Hom}_{K^b(A)}(T_c,T_c)&=\operatorname{dim
Hom}_{A}(P_{j}\oplus P_j,P_j\oplus P_j)-\operatorname{dim Hom}_{A}(P_j\oplus P_j, P_i)\\
&-\operatorname{dim
Hom}_{A}(P_i,P_{j}\oplus P_j)+\operatorname{dim Hom}_{A}(P_i, P_i)\\
&=4\#(j,j)-2\#(j,i)-2\#(i,j)+\#(i,i)
\end{align*}
If $\#(i,j)\neq 0$ there must be only one path $\Pi$ from $j$ to
$i$ not involving $\lambda$ because the algebra is gentle. Notice
the only from $j$ to $j$ not passing through $i$ must be the
trivial path because the algebra is finite dimensional. Then
$$(j,j)=\{ 1_j, \delta\lambda\Pi \},(j,i)=\{ \delta, \delta\lambda \},(i,j)=\{ \Pi, \lambda\Pi \},(i,i)=\{ 1_i, \lambda \}, $$
so
$$\operatorname{dim Hom}_{K^b(A)}(T_c,T_c)=8-4-4+2=2.$$
If $\#(i,j)=0$ for each path $\Pi$ from $j$ to $j$ there are
exactly two paths from $i$ to $j$ $\Pi\delta\lambda$ and
$\Pi\delta$ so $\#(j,i)=2\#(j,j)$. Therefore $2\#(j,i)=4\#(j,j)$
and
$$\operatorname{dim Hom}_{K^b(A)}(T_c,T_c)=\#(i,i)=2.$$

In both cases there is only one irreducible morphism from $T_c$ to
itself. This must be $\lambda'$ as it does not factorizes through
any of the previous irreducible morphisms we already analyzed.

We verified that there are no further irreducible morphisms and
this completes the proof.
%

\endproof

\section{Acknowledgement}
The research was supported by DGAPA and DGEP, UNAM and CONACyT.

I should like to thank my advisor Christof Geiss without whom this
work would have not been possible and also Corina Saenz and Bertha
Tom\'e for useful conversations. I would also like to thank
Thorsten Holm, Jan Sch\"roer and Alexander Zimmermann for their
comments and for giving their consent to display the results of
their unpublished manuscript in the Appendix.

\end{document}